\documentclass[author-year, preprint]{elsarticle}
\usepackage[margin=2.1cm]{geometry}
\journal{Transportation Research Part E: Logistics and Transportation Review}
\usepackage{lineno}
\usepackage{hyperref}

\usepackage{setspace}
\usepackage{color,graphicx,shortvrb}
\usepackage{multirow}
\usepackage{verbatim}
\usepackage{amsmath,amsfonts,amssymb}
\usepackage{rotating}
\usepackage{lscape}
\usepackage{pdflscape}
\usepackage{booktabs}
\usepackage[titletoc,title]{appendix}
\usepackage{hhline}
\usepackage{eurosym}
\usepackage{tikz}
\usepackage{tikz-qtree}
\usepackage{float}
\usepackage{graphics}
\usepackage{subcaption}
\usepackage{changepage}
\usepackage{soul}
\soulregister\cite7
\soulregister\ref7
\soulregister\eqref7
\soulregister\pageref7
\usepackage{cancel}
\usepackage{amsmath,amsfonts,amssymb}
\usepackage{enumitem}
\usetikzlibrary{arrows.meta}
\usepackage{soul}
\usepackage{booktabs}
\usepackage{makecell}
\usepackage{mathtools}

\usepackage{algorithm}
\usepackage[noend]{algpseudocode}

\usepackage{tikz-cd}
\usepackage{cases}
\usetikzlibrary{matrix,arrows,decorations.pathmorphing}
\usepackage{forest}
\usetikzlibrary{shadows,arrows.meta}

\usepackage{numcompress}\biboptions{authoryear}

\makeatletter
\def\ps@pprintTitle{%
 \let\@oddhead\@empty
 \let\@evenhead\@empty
 \def\@oddfoot{}%
 \let\@evenfoot\@oddfoot}
\makeatother

\newcommand{\be}{\begin{equation}}
\newcommand{\ee}{\end{equation}}
\newcommand{\beenum}{\begin{enumerate}}
\newcommand{\eenum}{\end{enumerate}}
\newcommand{\bi}{\begin{itemize}}
\newcommand{\ei}{\end{itemize}}

\definecolor{darkorange}{rgb}{1.0, 0.55, 0.0}

\usepackage{amsmath}
\tikzset{
dot/.style={circle,draw,inner sep=1.2,fill=black},
}
\DeclareMathOperator*{\argmin}{argmin}

\newtheorem{prop}{Proposition}

\begin{document}
  \begin{frontmatter}
	\title{Stochastic Fleet Size, Fleet Mix and Consistent Vehicle Routing Problem for Last Mile Delivery}
		
		\author[BeatriciMaggioniBiroliniMalighetti_address]{Paolo Beatrici}
		\ead{paolo.beatrici@unibg.it}

        \author[BeatriciMaggioniBiroliniMalighetti_address]{Sebastian Birolini}\ead{sebastian.birolini@unibg.it}
		\author[BeatriciMaggioniBiroliniMalighetti_address]{Francesca Maggioni\corref{mycorrespondingauthor}}
	\ead{francesca.maggioni@unibg.it}
	\cortext[mycorrespondingauthor]{Corresponding author}
		
	\author[BeatriciMaggioniBiroliniMalighetti_address]{Paolo Malighetti}\ead{paolo.malighetti@unibg.it}

	\address[BeatriciMaggioniBiroliniMalighetti_address]{{Department of Management, Information and Production Engineering, University of Bergamo, Viale G. Marconi 5, Dalmine 24044, Italy}}
		
		\begin{abstract}
      In this paper, we address the joint optimization of fleet size and mix, along with vehicle routing, under uncertain customer demand. We propose a two-stage stochastic mixed-integer programming model, where first-stage decisions concern the composition of the delivery fleet and the design of consistent baseline routes. In the second stage, approximate recourse actions are introduced to adapt the initial routes in response to realized customer demands. The objective is to minimize the total delivery cost, including vehicle acquisition, travel distance, and penalty costs for unserved demand.
To tackle the computational challenges arising in realistic problem instances, we develop a path-based reformulation of the model and design a Kernel Search-based heuristic to enhance scalability. Computational experiments on small synthetic instances, generated through a population-density-based sampling approach, are conducted to validate the formulation and assess the effects of demand stochasticity through standard stochastic measures, after applying a scenario reduction technique. Additional tests on large-scale real-world instances, based on data from the Italian postal company, demonstrate the effectiveness of the proposed approach and provide managerial and practical insights.
		\end{abstract}
		
		
	\end{frontmatter}
   \section{Introduction}\label{sec:introduction}

In recent years, urban logistics systems operating in dense metropolitan areas have undergone profound changes driven by rapid urbanization, the growth of e-commerce, and increasing sustainability concerns. These dynamics are particularly critical in city centers characterized by high customer density, limited road capacity, traffic congestion, and strict environmental regulations. In such contexts, the final stage of distribution, commonly referred to as the Last Mile, is widely recognized as one of the most costly and operationally challenging segments of supply chain management \citep{mohammad2023innovative}. At the same time, the continuous expansion of urban transportation systems \citep{shahmohammadi2020comparative} and the widespread use of internal combustion engine vehicles \citep{anderluh2017synchronizing} have exacerbated congestion, air pollution, and energy consumption, making sustainable Last Mile logistics an increasingly important research and practical challenge.

To address these issues, logistics operators are increasingly adopting heterogeneous delivery fleets that combine conventional vehicles with Light Freight Vehicles (LFVs), such as Electric Cargo Bikes (ECBs). Mixed fleets exploit the complementary characteristics of different vehicle types: while cargo vans remain essential for serving large delivery volumes and longer travel distances, ECBs are particularly suitable for dense urban environments because of their maneuverability, low environmental impact, and ability to access restricted traffic areas. Consequently, optimizing fleet composition has emerged as a key strategic decision for improving both operational efficiency and environmental sustainability \citep{fikar2018decision,Beatrici2026}.

Another major challenge in Last Mile logistics is the presence of uncertainty. Customer demand, travel times, customer availability, and other operational factors are typically unknown when tactical decisions are made and may substantially affect delivery performance. Consequently, uncertainty has received increasing attention in urban logistics, motivating the development of stochastic, robust, and distributionally robust optimization approaches \citep{perboli2017progressive,maggioni2019stochastic,cavagnini2024two,MAGGIONI2025,oyola2018stochastic, spinelli2026stochastic}. Among these sources of uncertainty, customer demand is particularly relevant because it directly influences fleet utilization, routing decisions, and overall system costs.

Although the literature has extensively investigated fleet size and mix decisions, stochastic vehicle routing, and route consistency, these topics have largely been studied independently. To the best of our knowledge, limited attention has been devoted to integrating fleet sizing and composition decisions with consistent routing strategies under stochastic customer demand within a unified optimization framework. Such an integrated perspective is particularly relevant for urban logistics operators, where tactical decisions regarding fleet composition interact closely with the need to maintain efficient and stable delivery operations despite daily demand fluctuations.

Motivated by these challenges, this paper introduces the Stochastic Fleet Size and Mix Consistent Vehicle Routing Problem (SFSMConVRP). The proposed model simultaneously determines the composition of a heterogeneous delivery fleet and the design of delivery routes while explicitly accounting for uncertain customer demand through a two-stage stochastic optimization framework.

A distinctive feature of the proposed approach is the incorporation of route consistency. Rather than redesigning delivery routes from scratch every day after customer requests become available, the model determines a set of stable baseline routes that can be repeatedly adopted over time. Daily operational decisions are then limited to modest adjustments required by the realized demand, while preserving the underlying routing structure. Following \citet{kovacs2014vehicle}, route consistency refers to maintaining stable customer sequences and travel paths across successive operational days. This planning strategy simplifies operational management, increases drivers' familiarity with their assigned service areas, and reduces unnecessary route variability, while preserving the flexibility required to accommodate stochastic demand.

Specifically, we formulate the SFSMConVRP as a two-stage stochastic mixed-integer linear programming model. In the first stage, tactical decisions determine the fleet composition together with the baseline consistent routes. After customer demand is revealed, approximate recourse actions modify these routes only when necessary, without disrupting their tactical structure.
Given the computational complexity of the resulting optimization problem, we derive a path-based reformulation and develop a tailored Kernel Search heuristic \citep{angelelli2010kernel} to improve computational scalability on large-scale instances.
    
The proposed methodology is first validated on a set of synthetically generated instances inspired by a real Last Mile delivery application of the Italian postal operator Poste Italiane. A scenario reduction procedure is employed to efficiently manage large scenario sets, and the value of explicitly modeling uncertainty is assessed through classical stochastic programming measures. Finally, the approach is evaluated on large-scale real-world instances derived from the operational network of Poste Italiane, providing managerial insights into fleet planning and consistent routing decisions under uncertainty.

The main contributions of this paper can be summarized as follows:

\begin{enumerate}
\item[(1)] We develop a novel two-stage stochastic optimization model for the fleet size and mix consistent vehicle routing problem with approximate recourse actions;
\item[(2)] We characterize the proposed recourse actions and derive valid inequalities that improve model scalability;
\item[(3)] We propose a path-based reformulation that reduces computational complexity and facilitates solution of large-scale instances;
\item[(4)] We design a Kernel Search-based heuristic that further enhances scalability;
\item[(5)] We perform numerical experiments on small instances to validate the theoretical properties and compare deterministic versus stochastic solutions using classical stochastic measures (\cite{maggioni2012analyzing});
\item[(6)] We test the effectiveness of the proposed methodology on a real-world case study.
\end{enumerate}

The remainder of the paper is organized as follows. Section \ref{sec:lit_rev} reviews the relevant literature. Section \ref{sec:problem_description_and_formulation} describes the problem and presents the mathematical formulations, along with valid inequality and characterization of recourse actions. Section \ref{sec:sol_meth} details the proposed solution approaches, while Section \ref{sec:real_case_study} discusses the experimental setup, computational results, and managerial implications. Finally, Section \ref{sec:conclusion} concludes the paper.

\section{Literature Review}\label{sec:lit_rev}

Last Mile delivery problems are commonly modeled as Vehicle Routing Problems (VRPs), in which a fleet of vehicles must serve a set of geographically dispersed customers while optimizing one or more performance criteria, such as transportation costs, travel distance, travel time, or service quality. Owing to their practical relevance, VRPs have been extensively studied over the past decades and extended to incorporate a wide variety of operational constraints and application-specific requirements (see \cite{toth2014vehicle} for a comprehensive overview).

Among the numerous VRP variants, the Fleet Size and Mix Vehicle Routing Problem (FSMVRP) and the Consistent Vehicle Routing Problem (ConVRP) are particularly relevant for Last Mile logistics. The FSMVRP addresses tactical decisions concerning the optimal size and composition of a heterogeneous fleet together with the associated routing decisions, whereas the ConVRP focuses on maintaining stable customer–driver assignments and routing patterns over time in order to improve service quality and operational regularity. Since the proposed problem combines these two research directions under demand uncertainty, we review them separately in the following.

The main contributions on the FSMVRP and its variants are summarized in Table \ref{tab:LR_overview_FSMVRP}. The reviewed studies differ with respect to the considered fleet characteristics, the presence of uncertainty, the mathematical formulation, and the adopted solution methodology.

\begin{table}[h!]
    \resizebox{\textwidth}{!}{
    \begin{tabular}{ccccccccc}
    \toprule
    Reference & \makecell{Problem\\type} & \makecell{Fleet with\\LFVs}  & \makecell{Presence of\\uncertainty} & Objective &\makecell{Model\\type\textsuperscript{a}} & \makecell{Solution\\Approach\textsuperscript{b}} & \makecell{Solution\\Method\textsuperscript{c}}\\
    \hline\\ [-1ex]
        \cite{golden1984fleet} & FSMVRP & - & - & \makecell{Minimize\\ fleet costs} & MILP & H & Savings heuristic\\  
         \\[0.5ex]
	\cite{gheysens1984comparison} & FSMVRP & - & - & \makecell{Minimize\\ fleet costs} & MILP & H & Assignment-based heuristic\\ 
	\\[0.5ex] 
	\cite{desrochers1991new} & FSMVRP & - & - & \makecell{Minimize\\ fleet costs} & - & H & Matching-based 		heuristic\\  
         \\[0.5ex]
	\cite{osman1996local} & FSMVRP & - & - & \makecell{Minimize\\ fleet costs} & LP & H & Tabu-search-based 		LS\\  
	\\[0.5ex]
        \cite{liu1999fleet} & \makecell{FSMVRP with\\Time Windows} & - & - & \makecell{Minimize\\fleet costs} & MILP & H & Insertion-based heuristic\\
        \\[0.5ex]
	\cite{dullaert2002new} & \makecell{FSMVRP with\\Time Windows} & - & - & \makecell{Minimize\\fleet costs} & - & H & Sequential insertion techniques\\
        \\[0.5ex]
	\cite{renaud2002sweep} & FSMVRP & - & - & \makecell{Minimize\\ fleet costs} & - & H & Sweep-based heuristic\\  
	\\[0.5ex]
        \cite{dell2007heuristic} & \makecell{FSMVRP with\\Time Windows} & - & - &\makecell{Minimize\\fleet costs} & MILP & H & CI heuristic, Metaheuristic\\
        \\[0.5ex]
	\cite{baldacci2009valid} & \makecell{FSMVRP with\\Time Windows} & - & - &\makecell{Minimize\\fleet costs} & MILP & E & CS, valid inequalities\\
        \\[0.5ex]
        \cite{salhi2013fleet} & \makecell{FSMVRP with\\Backhauls} & - & - & \makecell{Minimize\\fleet costs} & ILP & H & SP-based heuristic\\
        \\[0.5ex]

        \cite{hiermann2016electric} & \makecell{Electric\\FSMVRP} & $\tiny{\surd}$ & - & \makecell{Minimize\\fleet costs}& MILP & E,H & B\&P, ALNS, LS\\
        \\[0.5ex]

        \cite{kilby2016fleet} & FSMVRP & - & - & \makecell{Minimize\\delivery costs} & CP & H & LNS\\
        \\[0.5ex]

	\cite{pasha2016simple} & \makecell{Multi-period\\FSMVRP}  & - & Demand &\makecell{Minimize\\delivery costs} & MILP & H & TS\\
        \\[0.5ex]

        \cite{bertoli2020column} & FSMVRP & - & - & \makecell{Minimize\\fleet costs} &  - & H & CG-based heuristic\\
        \\[0.5ex]

        \cite{malladi2022stochastic} & \makecell{Stochastic\\FSMVRP} & - & Customers &\makecell{Minimize\\fleet costs}& MILP & H & SAA-based heuristic, ALNS\\
        \\[0.5ex]

        Current work & \makecell{Stochastic\\FSMConVRP} & $\tiny{\surd}$ &Demand & \makecell{Minimize\\ delivery costs} & MILP & E, H & CS, KS-based heuristic\\
        \\[0.5ex]

        \bottomrule
    \end{tabular}}
\caption{Classification of the references on the FSMVRP and its variants for Last Mile delivery problems.\\
\textsuperscript{a} MILP: Mixed Integer Linear program, LP: Linear program, ILP: Integer Linear program, CP: Constraint Programming.\\
    \textsuperscript{b} H: Heuristic, E: Exact.\\
    \textsuperscript{c} CI: Constructive Insertion, SP: Set Partitioning, B\&P: Branch\&Price, ALNS: Adaptive Large Neighbourhood Search, LS: Local Search, LNS: Local Neighbourhood Search, TS: Tabu Search, CG: Column Generation, SAA: Sample Average Approximation, CS: Commercial Solver, KS: Kernel Search.}
    \label{tab:LR_overview_FSMVRP}
\end{table}

The FSMVRP was originally introduced by \cite{golden1984fleet}, who formulated the problem under deterministic customer demand while jointly considering fixed vehicle ownership costs and variable routing costs. Since simultaneously determining fleet composition and vehicle routes leads to a computationally challenging optimization problem, early research primarily focused on developing efficient heuristic solution methods. Representative examples include the savings-based and route-first--cluster-second heuristics proposed in \cite{golden1984fleet}, the generalized assignment-based approach of \cite{gheysens1984comparison}, the matching-based tour fusion heuristic developed by \cite{desrochers1991new}, and the tabu-search local search procedure proposed by \cite{osman1996local}.

Subsequent research enriched the original formulation by incorporating additional operational constraints frequently encountered in distribution systems. Time-window constraints gave rise to the Fleet Size and Mix Vehicle Routing Problem with Time Windows (FSMVRPTW), first introduced by \cite{liu1999fleet}. This variant stimulated the development of increasingly sophisticated solution approaches, including insertion-based heuristics \cite{liu1999fleet}, sequential insertion procedures \cite{dullaert2002new}, and hybrid ruin-and-recreate metaheuristics \cite{dell2007heuristic}. In parallel, alternative mathematical formulations were proposed to improve computational performance. For example, \cite{renaud2002sweep} adopted a set-partitioning formulation based on routes generated through a sweep heuristic, whereas \cite{baldacci2009valid} developed a strengthened two-commodity network flow formulation incorporating valid inequalities.

More recent studies have focused on improving the realism of the FSMVRP by considering additional operational and sustainability-related aspects. Pickup operations after deliveries are incorporated in the FSMVRP with backhauls proposed by \cite{salhi2013fleet}, together with a set-partitioning-based heuristic solution approach. Constraint Programming has also been explored by \cite{kilby2016fleet}, who embedded a Large Neighbourhood Search (LNS) procedure within a CP framework to jointly optimize fleet composition and routing decisions.

The increasing interest in sustainable urban logistics has further motivated the integration of alternative vehicle technologies into fleet design. In particular, \cite{hiermann2016electric} introduced the Electric FSMVRP, explicitly accounting for battery limitations, recharging times, and charging infrastructure when determining fleet composition and delivery routes. Their solution framework combines an exact Branch-and-Price algorithm with an Adaptive Large Neighbourhood Search heuristic to efficiently solve large-scale instances. Fleet design decisions have also been investigated from a long-term planning perspective by \cite{bertoli2020column}, who considered the possibility of combining owned and hired vehicles and proposed a Column Generation-based heuristic to minimize acquisition, maintenance, hiring, and routing costs.

Overall, the FSMVRP literature has progressively evolved from the original deterministic formulation towards richer operational settings involving heterogeneous fleets, sustainability considerations, and increasingly effective exact and heuristic solution approaches. Nevertheless, the vast majority of existing contributions still assume deterministic delivery environments, despite the significant impact that uncertainty may have on tactical fleet planning decisions.

Modeling uncertainty has become increasingly important in Last Mile logistics, where customer demand, travel times, and service conditions are typically unknown when tactical planning decisions are made. Stochastic optimization provides a natural framework for addressing these sources of uncertainty (see, e.g., \cite{kall1994stochastic,birge2011introduction,shapiro2021lectures}). Despite its practical relevance, uncertainty has received only limited attention within the FSMVRP literature.

Among the few available contributions, \cite{pasha2016simple} consider demand variability over multiple planning periods and develop a Tabu Search heuristic to determine a common fleet composition across different demand realizations. More recently, \cite{malladi2022stochastic} introduced a two-stage stochastic FSMVRP in which uncertain customer requests are explicitly modeled through scenario-based optimization, jointly determining fleet composition and routing decisions while accounting for electric vehicles.

Although these contributions represent important advances towards stochastic fleet planning, they do not consider route consistency across operational days. Delivery routes are redesigned after demand realization, and no mechanism is introduced to preserve stable routing structures over time. Consequently, the interaction between fleet sizing and composition, stochastic customer demand, and route consistency remains largely unexplored.

Similarly to the FSMVRP, the Consistent Vehicle Routing Problem (ConVRP) has been progressively extended to address increasingly realistic operational settings while preserving route consistency over time. The main contributions on the ConVRP and its variants are summarized in Table \ref{tab:LR_overview_ConVRP}. Existing studies differ with respect to fleet characteristics, uncertainty modeling, mathematical formulation, and solution methodology.

The ConVRP was originally introduced by \cite{groer2009consistent} to improve customer satisfaction in the small-package delivery industry through stable customer--driver assignments over a multi-period planning horizon. The underlying motivation is that maintaining consistent service increases drivers' familiarity with their delivery areas and customers' preferences, while also reducing variability in delivery schedules.

Building on this idea, subsequent research has progressively relaxed the assumptions of the original formulation and incorporated additional operational features. For instance, \cite{kovacs2015generalized} proposed the Generalized ConVRP, allowing each customer to be served by a limited number of drivers rather than enforcing a strict one-driver-per-customer policy. Their formulation explicitly models the trade-off between routing efficiency and service consistency by penalizing variations in customer arrival times.

Several studies have further enriched the ConVRP by considering operational characteristics frequently encountered in practice. Simultaneous pickup and delivery operations are incorporated by \cite{zhen2020consistent}, whereas \cite{mancini2021collaborative} investigate a collaborative logistics setting in which customers may be exchanged among carriers while preserving time consistency, workload balance, and minimum market-share requirements. More recently, \cite{stavropoulou2022consistent} extended the ConVRP to heterogeneous vehicle fleets, accounting for different vehicle capacities and operating costs while maintaining both driver and arrival-time consistency. However, the fleet composition is assumed to be known in advance, and no tactical decisions regarding fleet sizing or fleet selection are considered.

\begin{table}[h!]
    \resizebox{\textwidth}{!}{
    \begin{tabular}{ccccccccc}
    \toprule
    Reference & \makecell{Problem\\type} & \makecell{Mixed\\fleet} & \makecell{Presence of\\uncertainty}& Objective & \makecell{Model\\type\textsuperscript{a}} & \makecell{Solution\\Approach\textsuperscript{b}} & \makecell{Solution\\Method\textsuperscript{c}}\\
    \hline\\ [-1ex]
        
        \cite{groer2009consistent} & ConVRP& - & - &\makecell{Minimize\\travel time} & MIP & H & Problem-based heuristic\\
        \\[0.5ex]

        \cite{kovacs2015generalized} & \makecell{Generalized\\ConVRP} & - & - & \makecell{Minimize travel\\and arrival time} & MIP & H & Flexible LNS\\
        \\[0.5ex]

        \cite{goeke2019exact} & ConVRP & - & - &\makecell{Minimize\\operating time} & MIP & E & CG-based approach\\
        \\[0.5ex]

        \cite{barros2020exact} & ConVRP & - & - & \makecell{Minimize\\ travel time} & MILP & E & CS \\
        \\[0.5ex]

        \cite{zhen2020consistent} & \makecell{ConVRP with\\ distribution\\and collection} & - & - &\makecell{Minimize\\operating time} & MILP & H & RTR-travel algorithm, LSVNS, TS\\
        \\[0.5ex]

        \cite{mancini2021collaborative} & \makecell{Collaborative\\ConVRP}& - & - & \makecell{Maximize\\profit} & MILP & H & Matheuristic, ILS-based approach\\
        \\[0.5ex]

        \cite{yang2021consistent} & \makecell{Uncertain\\ConVRP} & - &\makecell{Demand, Travel\\ and service time} &\makecell{Minimize\\travel time} & MILP & H & ABC heuristic\\
        \\[0.5ex]
        
        \cite{stavropoulou2022consistent} & \makecell{Heterogeneous\\ConVRP} & $\small{\surd}$ & - & \makecell{Minimize\\transportation\\costs} & MILP & H & TS, VND algorithm\\
        \\[0.5ex]

        \cite{alvarez2024consistent} & \makecell{Stochastic\\ConVRP} & - & Demand & \makecell{Minimize\\consistency\\violations} & MILP & E & SAA, Benders decomposition\\
        \\[0.5ex]

        Current work & \makecell{Stocahstic\\FSMConVRP} & $\small{\surd}$ & Demand & \makecell{Minimize\\delivery costs} &  MILP & E, H & CS, KS-based heuristic\\
        \\[0.5ex]
        
        \bottomrule
    \end{tabular}}
\caption{Classification of the references on the ConVRP and its variants for Last Mile delivery problems.\\
\textsuperscript{a} MIP: Mixed Integer Program, MILP: Mixed Integer Linear Program. \\
    \textsuperscript{b} H: Heuristic, E: Exact.\\
    \textsuperscript{c} LNS: Large Neighbourhood Search, CG: Column Generation, MILP: Mixed Integer Linear Program, RTR: Record-To-Record, LSVNS: Local Search with Variable Neibourhood Search, TS: Tabu Search, ILS: Iterated Local Search, ABC: Artificial Bee Colony, VND: Variable Neighbourhood Descent, SAA: Sample Average Approximation, CS: Commercial Solver, KS: Kernel Search.}
    \label{tab:LR_overview_ConVRP}
\end{table}

As with the FSMVRP, the computational complexity of the ConVRP rapidly increases with problem size, making exact solution methods impractical for many real-world applications. Consequently, most research has focused on the development of efficient heuristic and matheuristic approaches capable of producing high-quality solutions within reasonable computational times.

Representative examples include the flexible Large Neighbourhood Search proposed by \cite{kovacs2015generalized}, which iteratively destroys and repairs routes to improve consistency, the Record-to-Record, Local Search with Variable Neighbourhood Search, and Tabu Search heuristics developed by \cite{zhen2020consistent}, the matheuristic and Iterated Local Search framework introduced by \cite{mancini2021collaborative}, and the Tabu Search combined with Variable Neighbourhood Descent proposed by \cite{stavropoulou2022consistent}. Although heuristic methods dominate the literature, several exact approaches have also been developed. In particular, \cite{goeke2019exact} proposed a Column Generation-based exact algorithm, while \cite{barros2020exact} introduced a strengthened mixed-integer linear programming formulation to improve computational tractability.

Despite the practical relevance of uncertainty in Last Mile logistics, only a limited number of studies have incorporated stochastic elements into the ConVRP. Most existing contributions assume deterministic customer demands and travel conditions, focusing primarily on service consistency.

Among the few exceptions, \cite{yang2021consistent} introduced an Uncertain ConVRP that simultaneously considers stochastic customer demand, travel times, and service durations. More recently, \cite{alvarez2024consistent} proposed a two-stage stochastic programming formulation in which customer-to-driver assignments are determined before demand realization, while routing decisions are optimized afterwards using Sample Average Approximation and Benders decomposition.

Although these contributions represent important advances in stochastic consistent routing, they assume a predefined fleet and do not optimize fleet size or fleet composition. Consequently, the interaction between tactical fleet planning and route consistency under uncertain customer demand remains largely unexplored.

Overall, the above review highlights three main research gaps. First, the FSMVRP literature has extensively investigated heterogeneous fleet design and routing decisions, but only a few studies explicitly account for demand uncertainty, and none of them incorporates route consistency into the tactical planning process. Second, although the ConVRP literature has introduced stochastic formulations to preserve service consistency under uncertainty, fleet size and composition are generally assumed to be fixed and are not optimized jointly with routing decisions. Third, despite the increasing adoption of Light Freight Vehicles (LFVs), including Electric Cargo Bikes (ECBs), in sustainable Last Mile logistics, their integration within a stochastic framework that simultaneously addresses fleet sizing, fleet composition, and route consistency remains unexplored.

To bridge these gaps, this paper introduces the \emph{Stochastic Fleet Size and Mix Consistent Vehicle Routing Problem} (SFSMConVRP), a novel optimization framework that jointly determines the composition of a heterogeneous delivery fleet and the design of consistent delivery routes under uncertain customer demand. The proposed model is formulated as a two-stage stochastic mixed-integer linear program, where tactical decisions define fleet composition and baseline consistent routes, while approximate recourse actions adapt route execution after demand realization without disrupting their underlying structure. To solve realistic instances efficiently, we further develop a path-based reformulation and a tailored Kernel Search heuristic. The effectiveness of the proposed methodology is finally assessed through extensive computational experiments on both synthetic benchmark instances and a real-world Last Mile delivery case study based on data provided by the Italian postal operator \emph{Poste Italiane}.

\section{Problem Description and Formulation}\label{sec:problem_description_and_formulation}

We consider a company responsible for a daily parcel delivery service that aims to determine the optimal fleet size and composition, along with consistent routing plans for its distribution operations. The objective is to design a delivery system that avoids ad-hoc routing decisions typically made by drivers at the beginning of each working day, once actual delivery requests are known. Instead, the company seeks to establish consistent, repeatable routes that can be operated regularly, while maintaining flexibility to adapt to daily variations in customer demand. This setting reflects the strategic and tactical nature of the problem, where long-term decisions on fleet acquisition and route design must account for operational adjustments driven by uncertain customer demands.
 
Let $\mathcal{G}=(\mathcal{I}^*,\mathcal{A})$ be a complete directed graph, where $\mathcal{I}^*=\{0,\ldots,N\}$ denotes the set of nodes and $\mathcal{A}=\{(i,j): i,j\in\mathcal{I}^*\,, i\neq j\}$ the set of directed arcs. Node $0$ represents the single depot, referred to as Distribution Center (DC), from which all delivery operations start and end. The set of customers is denoted by $\mathcal{I}=\{1,\ldots,N\}\subset\mathcal{I}^*$. The distance between nodes $i$ and $j$ is denoted by $\delta_{ij}$.

Customer demand is uncertain. For each node $i\in\mathcal{I}$, the random variable $d_i$ represents its demand, and these random variables are assumed to be independent. Let $\mathcal{S}=\{1,\ldots,S\}$ denote the finite set of demand scenarios. Accordingly, $d_{is}$ indicates the realized non-negative demand of customer $i$ under scenario $s\in\mathcal{S}$, and $\pi_s\in[0,1]$ denotes the probability of scenario $s$, with $\sum_{s\in\mathcal{S}}\pi_s=1$.

Deliveries are carried out by a heterogeneous fleet of vehicles. Let $\mathcal{P}=\{1,\dots,P\}$ be the set of vehicle types. For each type $p\in\mathcal{P}$, $l_p$ and $f_p$ denote the load capacity and the fixed ownership cost, respectively. The maximum capacity among all vehicle types is indicated by $l^{\text{max}}$. Each vehicle type $p$ is characterized by a fixed deployment cost $c_{ip}^{\text{fix}}$ and a variable travel cost $c_{ijp}^{\text{var}}$. The fixed deployment cost $c^{\text{fix}}_{ip}$ is associated with activating a vehicle of type $p$ whose consistent route terminates at customer $i$. It includes the fixed ownership cost $f_p$ of vehicle type $p$ and the travel cost required to return from the last visited customer $i$ to the depot. Hence, $c^{\text{fix}}_{ip}=f_p+\omega_{\delta_p}\delta_{i0}$, for $i\in\mathcal{I}$ and $p\in\mathcal{P}$.
The variable travel cost $c^{\text{var}}_{ijp}$ is incurred whenever a vehicle of type $p$ travels from node $i$ to node $j$ along a first-stage route, and is defined as $c^{\text{var}}_{ijp}=\omega_{\delta_p}\delta_{ij}$ for $i,j\in\mathcal{I}^*$ with $i\neq j$ and $p\in\mathcal{P}$.
  Since all deliveries must be completed within a working day, a large time window is imposed to guarantee that operations are completed by the end of the shift. Parameter $\bar{t}$ represents this time limit, and $t_{ijp}$ denotes the travel time for a vehicle of type $p$ moving from node $i$ to node $j$.

To model fleet sizing, composition, and consistent routing decisions, we define the following binary variables:
\begin{itemize}
\item $z_{ip}=1$ if there exists a route operated by a vehicle of type $p$ that terminates at node $i$ (i.e., $i$ is the last customer visited on that route); 0 otherwise. Summing $z_{ip}$ over all $i$ provides the number of vehicles of type $p$ included in the delivery fleet;
\item $v_{ip}=1$ if a route operated by a vehicle of type $p$ serves customer $i$ without terminating there; 0 otherwise; 
\item $x_{ijp}=1$ if a vehicle of type $p$ travels directly from node $i\in\mathcal{I}^*$ to node $j\in\mathcal{I}$ ($i\neq j$); 0 otherwise.
\end{itemize}

These first-stage variables define the tactical decisions associated with the fleet configuration and the design of consistent delivery routes. Since routes are not indexed by individual vehicles, each active route associated with vehicle type $p$ can be interpreted as being operated by one vehicle of that type.

After the realization of customer demand in each scenario $s \in \mathcal{S}$, the first-stage routes remain unchanged in order to preserve tactical consistency. Nevertheless, operational recourse actions may be required to ensure feasibility with respect to the realized delivery requests.

Customer demand is stochastic for all nodes. If a node $i$ is not included in the first-stage routing plan, its realized demand $d_{is}$ under scenario $s$ is treated as \emph{additional} demand, and $i$ is referred to as an \emph{additional customer}. Such demand can be accommodated through recourse actions performed by the vehicles in the fleet, subject to residual load capacity and operating time constraints, or alternatively outsourced to an external delivery service to guarantee full demand coverage. For customers already assigned to first-stage routes, realized demand may differ from its expected value. Excess demand is handled through the same recourse mechanism, whereas customers with zero realized demand are still visited to preserve route consistency, although no delivery is performed.

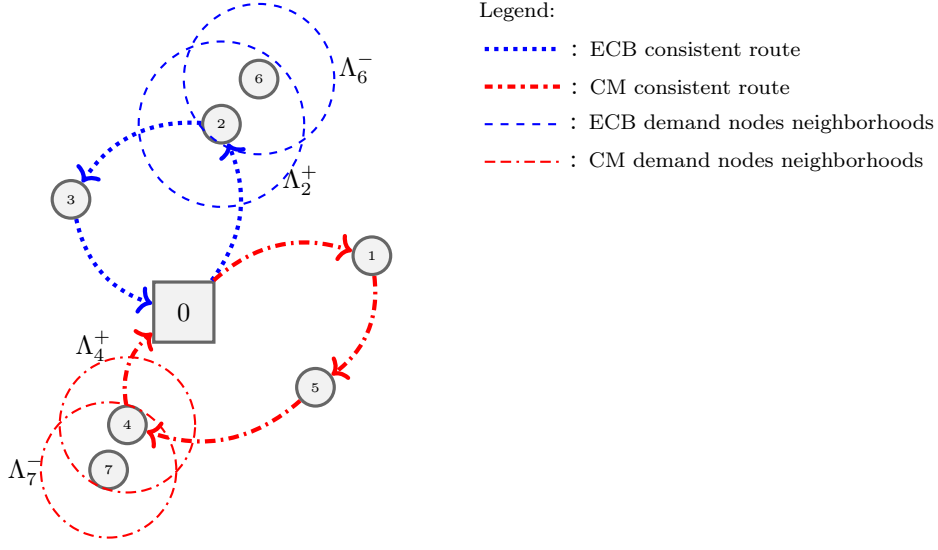
\begin{figure}[h!]
    \centering
    \begin{tikzpicture}[roundnode/.style={circle, draw=black!60, fill=black!5, very thick, minimum size=5mm},squarednode/.style={rectangle, draw=black!60, fill=black!5, very thick, minimum size=8mm}]
    \begin{scope}[scale=0.5]
        \node[squarednode] (0) at (0,0) {0};
        \node[roundnode] (1) at (5,1.5) {\footnotesize \tiny 1};
        \node[roundnode] (2) at (1,5) {\tiny 2};
        \node[roundnode] (3) at (-3,3) {\tiny 3};
        \node[roundnode] (4) at (-1.5,-3) {\tiny 4};
        \node[roundnode] (5) at (3.5,-2) {\tiny 5};
        \node[roundnode] (6) at (2,6.2) {\tiny 6};
        \node[roundnode] (7) at (-2,-4.2) {\tiny 7};

        \path (0) edge[->,ultra thick, dotted, bend right, blue]node[right]{} (2);
        \path (2) edge[->,ultra thick, dotted, bend right,blue]node[above]{}  (3);
        \path (3) edge[->,ultra thick, dotted, bend right,blue]node[right]{} (0);

        \path (0) edge[->,ultra thick, dash pattern=on 4pt off 2pt on 1pt off 2pt, bend left,red]node[below]{}(1);
        \path (1) edge[->,ultra thick, dash pattern=on 4pt off 2pt on 1pt off 2pt, bend left,red]node[right]{} (5);
        \path (5) edge[->,ultra thick, dash pattern=on 4pt off 2pt on 1pt off 2pt, bend left,red]node[below]{} (4);
        \path (4) edge[->,ultra thick, dash pattern=on 4pt off 2pt on 1pt off 2pt, bend left,red] node[above left]{} (0);

        \draw[color=black!60, thick, dashed,blue](1,5) circle (2.2);
        \node at (3.1,3.5) {$\Lambda_2^{+}$};
        \draw[color=black!60, thick, dashed,blue](2,6.2) circle (2);
        \node at (4.6,6.5) {$\Lambda_6^{-}$};
        \draw [color=black!60, thick, dash pattern=on 4pt off 2pt on 4pt off 2pt on 1pt off 2pt,red](-1.5,-3) circle (1.8);
        \node at (-2.4,-0.9) {$\Lambda_4^{+}$};
        \draw [color=black!60, thick, dash pattern=on 4pt off 2pt on 4pt off 2pt on 1pt off 2pt,red](-2,-4.2) circle (1.8);
        \node at (-4.2,-4.2) {$\Lambda_7^{-}$};
        \end{scope}
\begin{scope}[shift={(4,3.5)}] 
        \draw[](-0.2,0.5) node[right, black]{\footnotesize Legend:};
        \draw[blue, ultra thick, dotted] (0,0) -- (1,0) node[right, black]{: \footnotesize ECB consistent route};
        \draw[red, ultra thick, dash dot] (0,-0.5) -- (1,-0.5) node[right, black]{: \footnotesize CM consistent route};
        \draw[blue, thick, dashed] (0,-1) -- (1,-1) node[right, black]{: \footnotesize ECB demand nodes neighborhoods};
        \draw[red, thick, dash pattern=on 4pt off 2pt on 4pt off 2pt on 1pt off 2pt] (0,-1.5) -- (1,-1.5) node[right, black]{: \footnotesize CM demand nodes neighborhoods};
        \end{scope}
    
    \end{tikzpicture}
    
    \caption{Example delivery network with consistent routes assigned to fleet vehicles and neighborhoods for selected demand nodes. Two vehicle types are considered: Conventional Motorcycles (CMs) and Electric Cargo Bikes (ECBs).}

    \label{different_size_neighbourhood}
\end{figure}

To capture the local flexibility allowed by the proposed recourse mechanism, we define for each customer $i\in\mathcal{I}$ the forward neighborhood
$\Lambda_i^{+}=\{j\in\mathcal{I}\mid \delta_{ij}\leq \bar{\delta}\},
$ where $\bar{\delta}$ is a predefined distance threshold. The set $\Lambda_i^{+}$ contains all customers that can be served through a recourse action originating from customer $i$, and, by definition, $i\in\Lambda_i^{+}$. Figure~\ref{different_size_neighbourhood} illustrates an example of these neighborhoods in a simplified delivery network, showing how local recourse actions may extend consistent routes to accommodate stochastic demand realizations.

For completeness, we also define the backward neighborhood of customer $j\in\mathcal{I}$ as
$\Lambda_j^-=\{i\in\mathcal{I}\mid \delta_{ij}\leq\bar{\delta}\},$
which contains all customers from which $j$ can be served through a feasible recourse action. By construction,
$\Lambda_j^-=\{\,i\in\mathcal{I}\mid j\in\Lambda_i^+\}$, for all $j\in\mathcal{I}$. The proposed framework can be readily generalized by allowing customer-dependent neighborhood sizes, thus providing additional flexibility in modeling the recourse actions.

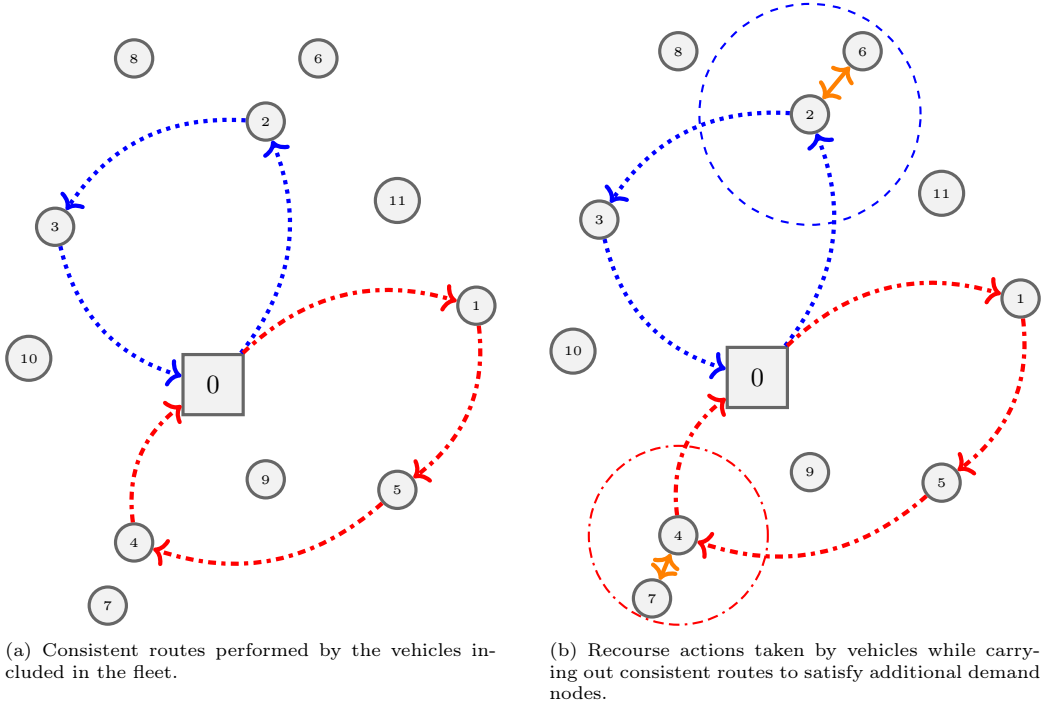
\begin{figure}[h!]
    \centering
    \subfloat[Consistent routes performed by the vehicles included in the fleet.\label{(2a)}]{\begin{tikzpicture}[roundnode/.style={circle, draw=black!60, fill=black!5, very thick, minimum size=0.3mm},squarednode/.style={rectangle, draw=black!60, fill=black!5, very thick, minimum size=8mm}]
    \begin{scope}[scale=0.7]
        \node[squarednode] (0) at (0,0) {0};
        \node[roundnode] (1) at (5,1.5) {\tiny 1};
        \node[roundnode] (2) at (1,5) {\tiny 2};
        \node[roundnode] (3) at (-3,3) {\tiny 3};
        \node[roundnode] (4) at (-1.5,-3) {\tiny 4};
        \node[roundnode] (5) at (3.5,-2) {\tiny 5};
        \node[roundnode] (6) at (2,6.2) {\tiny 6};
        \node[roundnode] (7) at (-2,-4.2) {\tiny 7};
        \node[roundnode] (8) at (-1.5,6.2) {\tiny 8};
        \node[roundnode] (9) at (1,-1.8) {\tiny 9};
        \node[roundnode] (10) at (-3.5,0.5) {\tiny 10};
        \node[roundnode] (11) at (3.5,3.5) {\tiny 11};
	\node[roundnode] (12) at (-3, -2.5) {\tiny 12};

        \path (0) edge[->,ultra thick, dotted, bend right,blue]node[right]{} (2);
        \path (2) edge[->,ultra thick, dotted, bend right,blue]node[above]{}  (3);
        \path (3) edge[->,ultra thick, dotted, bend right,blue]node[right]{} (0);

        \path (0) edge[->,ultra thick, dash dot, bend left,red]node[below]{}(1);
        \path (1) edge[->,ultra thick, dash dot, bend left,red]node[right]{} (5);
        \path (5) edge[->,ultra thick, dash dot, bend left,red]node[below]{} (4);
        \path (4) edge[->,ultra thick, dash dot, bend left,red] node[above left]{} (0);

    \end{scope}
    \end{tikzpicture}}\quad\quad
    \subfloat[Recourse actions taken by vehicles while carrying out consistent routes to satisfy additional demand nodes.\label{(2b)}]{\begin{tikzpicture}[roundnode/.style={circle, draw=black!60, fill=black!5, very thick, minimum size=0.3mm},squarednode/.style={rectangle, draw=black!60, fill=black!5, very thick, minimum size=8mm}]
    \begin{scope}[scale=0.7]
        \node[squarednode] (0) at (0,0) {0};
        \node[roundnode] (1) at (5,1.5) {\tiny 1};
        \node[roundnode] (2) at (1,5) {\tiny 2};
        \node[roundnode] (3) at (-3,3) {\tiny 3};
        \node[roundnode] (4) at (-1.5,-3) {\tiny 4};
        \node[roundnode] (5) at (3.5,-2) {\tiny 5};
        \node[roundnode] (6) at (2,6.2) {\tiny 6};
        \node[roundnode] (7) at (-2,-4.2) {\tiny 7};
        \node[roundnode] (8) at (-1.5,6.2) {\tiny 8};
        \node[roundnode] (9) at (1,-1.8) {\tiny 9};
        \node[roundnode] (10) at (-3.5,0.5) {\tiny 10};
        \node[roundnode] (11) at (3.5,3.5) {\tiny 11};
	\node[roundnode] (12) at (-3, -2.5) {\tiny 12};

        \path (0) edge[->,ultra thick, dotted, bend right,blue]node[right]{} (2);
        \path (2) edge[->,ultra thick, dotted, bend right,blue]node[above]{}  (3);
        \path (3) edge[->,ultra thick, dotted, bend right,blue]node[right]{} (0);
        \path (2) edge[<->,ultra thick,orange]node[right]{} (6);

        \path (0) edge[->,ultra thick, dash dot, bend left,red]node[below]{}(1);
        \path (1) edge[->,ultra thick, dash dot, bend left,red]node[right]{} (5);
        \path (5) edge[->,ultra thick, dash dot, bend left,red]node[below]{} (4);
        \path (4) edge[->,ultra thick, dash dot, bend left,red] node[above left]{} (0);
        \path (4) edge[<->,ultra thick,orange]node[right]{} (7);
	\path (4) edge[<->,ultra thick,orange]node[right]{} (12);

        \draw[color=black!60, thick, dashed,blue](1,5) circle (2.1);
        \draw[color=black=!60, thick, dash pattern=on 4pt off 2pt on 4pt off 2pt on 1pt off 2pt,red](-1.5,-3) circle (2);
    \end{scope}
    \end{tikzpicture}}

 \caption{Example of consistent routes and corresponding approximate recourse actions (shown as orange arrows) performed by the vehicles while serving the delivery network. The same vehicle types as in Figure~\ref{different_size_neighbourhood} are considered.}
\label{fig:mean_and_adjusted_route}
\end{figure}

We are now ready to describe the approximate recourse actions adopted in the second stage to adapt the first-stage consistent routes after demand realization. Figure~\ref{fig:mean_and_adjusted_route} illustrates these recourse actions on a representative delivery network.

Recourse actions are represented by non-negative continuous variables $y_{ijs}$, which denote the proportion of the realized demand of customer $j\in\Lambda_i^{+}$ that is served through a recourse action originating from customer $i\in\mathcal{I}$ under scenario $s\in\mathcal{S}$. These variables quantify additional deliveries performed from customer $i$ to neighboring customers $j$ in response to realized demand deviations.

The proposed recourse mechanism is deliberately modeled in an approximate manner. In particular, the variables $y_{ijs}$ do not explicitly represent the routing decisions required to travel from customer $i$ to customer $j$, nor do they allow for sequential recourse operations, such as serving additional customers starting from $j$.

The adoption of approximate recourse actions is motivated by both computational and modeling considerations. From a computational perspective, modeling these second-stage routing decisions explicitly would require a substantial number of additional binary variables to represent scenario-dependent routes. This would significantly increase the size and complexity of the formulation and rendering the model computationally intractable for the real-world instances considered in this work. From a modeling perspective, it is consistent with the tactical nature of the problem, whose primary objective is to design stable and repeatable delivery routes. Fleet vehicles operate along fixed first-stage routes that are intended to be executed on a daily basis, while recourse actions provide the flexibility required to accommodate realized demand without modifying the underlying routing structure. In this way, the proposed approach balances operational adaptability and route consistency, ensuring robustness with respect to demand variability while maintaining a tractable optimization model.

Each recourse action incurs an additional cost equal to $\beta c^{\text{rec}}_{ij}y_{ijs}$. Parameter $c^{\text{rec}}_{ij}$ denotes the unit cost of performing a recourse action from customer $i$ to customer $j$ and depends on the recourse operating cost $\omega_{\text{rec}}$ and on the additional travel distance $\delta_{ij}$. Specifically, $c^{\text{rec}}_{ij}=\omega_{\text{rec}}\delta_{ij},
\
i\in\mathcal{I},j\in\Lambda_i^+.$ 
Parameter $\beta$ is a tunable penalty coefficient controlling the level of conservatism of the solution. For example, setting $\beta=2$ assumes that serving customer $j$ from customer $i$ requires a round trip, thereby providing a conservative approximation of the corresponding routing cost.

Given the variability of customer demand and the limited load capacity of fleet vehicles, part of the realized demand may remain unserved by the fleet and must therefore be outsourced to an external provider or postponed to a subsequent delivery round (see \cite{malladi2022stochastic}). This situation is modeled through non-negative continuous variables $w_{is}$, defined for each customer $i\in\mathcal{I}$ and scenario $s\in\mathcal{S}$, representing the proportion of the realized demand at customer $i$ that is not served by the fleet. Consequently, $w_{is}=0$ indicates that the entire realized demand is satisfied by the company's fleet, whereas $w_{is}=1$ corresponds to complete outsourcing of customer $i$'s demand under scenario $s$. To discourage outsourced or postponed deliveries, a penalty coefficient $\gamma$ is introduced and typically chosen larger than $\beta$. Accordingly, each unit of unserved demand incurs an additional penalty equal to $\gamma d_{is}w_{is}$.

Because of vehicle capacity limitations, it is necessary to track the cumulative quantity delivered along each route. To this end, non-negative continuous variables $u_{is}$ are introduced, representing the cumulative quantity delivered up to and including customer $i$ under scenario $s$.

Furthermore, since all delivery operations must be completed within a prescribed working-time window, the temporal evolution of each route is explicitly modeled. Non-negative continuous variables $\tau_{is}$ denote the cumulative operating time required to reach and serve customer $i$ under scenario $s$. This quantity is computed as the sum of travel times between consecutive customers together with the service times accumulated up to and including customer $i$. Parameter $\hat{t}_i$ denotes the service time required at customer $i$, which is assumed to be independent of the vehicle type.

The decision variables introduced above define the first- and second-stage decisions of the proposed optimization framework and constitute the basis of the two-stage stochastic mixed-integer programming formulation presented in the following subsection.

\subsection{A Stochastic Two-Stage Mixed Integer Linear Programming Formulation}\label{sec:integrated_MILP}

In this section, we present a two-stage stochastic mixed-integer linear programming formulation for the proposed Stochastic Fleet Size and Mix Consistent Vehicle Routing Problem (SFSMConVRP). For ease of presentation, the sets, parameters, and decision variables used throughout the formulation are summarized below.

\bigskip

\noindent\emph{\underline{Sets}}\\
$\mathcal{I}=\{1,\ldots,N\}$: set of customers;\\
$\mathcal{I}^{\ast}=\{0,1,\ldots,N\}$: set of nodes, including the Distribution Center (DC), denoted by node $0$;\\
$\mathcal{P}=\{1,\ldots,P\}$: set of vehicle types;\\
$\mathcal{S}=\{1,\ldots,S\}$: set of demand scenarios;\\
$\Lambda_i^{+}$: forward neighborhood of customer $i\in\mathcal{I}$, containing the customers that can be served through a recourse action originating from $i$;\\
$\Lambda_i^{-}$: backward neighborhood of customer $i\in\mathcal{I}$, containing the customers from which $i$ can be served through a recourse action.

\bigskip

\noindent\emph{\underline{Deterministic parameters}}\\
$f_p$: fixed ownership cost of a vehicle of type $p\in\mathcal{P}$;\\
$\omega_{\delta_p}$: operating cost per unit of travel distance for vehicles of type $p\in\mathcal{P}$;\\
$l_p$: load capacity of vehicles of type $p\in\mathcal{P}$;\\
$l^{\max}$: maximum load capacity among all vehicle types;\\
$\Delta_p$: driving range of vehicles of type $p\in\mathcal{P}$;\\
$\delta_{ij}$: travel distance from node $i$ to node $j$, with $i,j\in\mathcal{I}^{\ast}$ and $i\neq j$;\\
$c_{ijp}^{\mathrm{var}}$: variable travel cost incurred by a vehicle of type $p\in\mathcal{P}$ traveling from node $i$ to node $j$;\\
$c_{ip}^{\mathrm{fix}}$: fixed deployment cost of a vehicle of type $p\in\mathcal{P}$ whose route terminates at customer $i$;\\
$c_{ij}^{\mathrm{rec}}$: unit cost associated with a recourse action from customer $i$ to customer $j$;\\
$t_{ijp}$: travel time of a vehicle of type $p\in\mathcal{P}$ from node $i$ to node $j$;\\
$\hat{t}_i$: service time required at customer $i\in\mathcal{I}$;\\
$\omega_{\mathrm{rec}}$: operating cost per unit of travel distance associated with recourse actions;\\
$\bar{t}$: maximum duration of a working day;\\
$\bar{\delta}$: neighborhood radius defining the sets $\Lambda_i^{+}$ and $\Lambda_i^{-}$;\\
$\beta$: penalty coefficient associated with recourse actions;\\
$\gamma$: penalty coefficient associated with unserved demand.

\bigskip

\noindent\emph{\underline{Stochastic parameters}}\\
$d_{is}$: realized demand of customer $i\in\mathcal{I}$ under scenario $s\in\mathcal{S}$;\\
$\pi_s$: probability of occurrence of scenario $s\in\mathcal{S}$.

\bigskip

\noindent\emph{\underline{Decision variables}}\\
$z_{ip}$: binary variable equal to 1 if a route operated by a vehicle of type $p\in\mathcal{P}$ terminates at customer $i\in\mathcal{I}$, and 0 otherwise;\\
$v_{ip}$: binary variable equal to 1 if customer $i\in\mathcal{I}$ is visited by a route operated by a vehicle of type $p\in\mathcal{P}$ without being its terminal customer, and 0 otherwise;\\
$x_{ijp}$: binary variable equal to 1 if a vehicle of type $p\in\mathcal{P}$ travels directly from node $i\in\mathcal{I}^{\ast}$ to customer $j\in\mathcal{I}$, with $i\neq j$, and 0 otherwise;\\
$y_{ijs}$: non-negative variable representing the proportion of the realized demand of customer $j\in\Lambda_i^{+}$ served through a recourse action originating from customer $i$ under scenario $s\in\mathcal{S}$;\\
$w_{is}$: non-negative variable representing the proportion of the realized demand of customer $i\in\mathcal{I}$ that is not served by the fleet under scenario $s\in\mathcal{S}$;\\
$u_{is}$: non-negative variable representing the cumulative quantity delivered up to and including customer $i$ under scenario $s\in\mathcal{S}$;\\
$\tau_{is}$: non-negative variable representing the cumulative operating time up to and including customer $i$ under scenario $s\in\mathcal{S}$.

We now present the proposed two-stage stochastic mixed-integer linear programming model, denoted by $\mathcal{M}$, for the SFSMConVRP.

\allowdisplaybreaks

\begin{align}
\min  \quad & \sum_{i\in\mathcal{I}}\sum_{p\in \mathcal{P}} c_{ip}^{fix} z_{ip} + \sum_{i\in\mathcal{I}^\ast}\sum_{j\in\mathcal{I}\setminus\{i\}}\sum_{p\in\mathcal{P}}c_{ijp}^{var} x_{ijp} +  \sum_{s\in\mathcal{S}}\pi_s\Bigg( \sum_{i\in\mathcal{I}}\sum_{j\in\Lambda_{i}^{+}}\beta c_{ij}^{rec} y_{ijs} +\sum_{i\in\mathcal{I}}\gamma d_{is}w_{is}\Bigg) \label{stochastic_objective_function} \\
\text{s.t.}\quad & \sum_{j\in\mathcal{I}} x_{0jp}=\sum_{i\in\mathcal{I}}z_{ip}\quad\quad p\in\mathcal{P}\label{stochastic_arcs=trips}\\
& \sum_{i\in\mathcal{I}^\ast\setminus\{j\}} x_{ijp}= z_{jp}+v_{jp} \quad\quad p\in\mathcal{P},\, j \in\mathcal{I}\label{stochastic_nr_vehicles_equal_trips_from_depot}\\
& \sum_{j\in\mathcal{I}\setminus\{i\}} x_{ijp}= v_{ip} \quad\quad p\in\mathcal{P},\, i\in\mathcal{I}\label{stochastic_outgoing_arc_from_not_end_customers}\\
& \sum_{p\in\mathcal{P}} (z_{ip}+ v_{ip}) \leq 1 \quad\quad i\in\mathcal{I}\label{stochastic_one_customer_one_type}\\
&y_{ijs} \leq \sum_{p\in\mathcal{P}}(z_{ip}+v_{ip}) \quad\quad i \in \mathcal{I},\, j \in \Lambda_i^{+},\, s \in \mathcal{S}\label{stochastic_demand_servable only_if_i_is_served}\\
&\sum_{i \in \Lambda_j^{-}}  y_{ijs} + w_{js} = 1 \quad\quad j \in \mathcal{I},\, s \in \mathcal{S}\label{stochastic_demand_satisfaction}\\
&u_{is}\leq \sum_{p \in \mathcal{P}} (z_{ip} + v_{ip}) l_p+ \Big[1-\sum_{p \in \mathcal{P}} (z_{ip} + v_{ip}) \Big] l^{max} \quad\quad  i \in\mathcal{I},\, s \in \mathcal{S}\label{stochastic_UB_load_capacity}\\
&u_{is} \geq d_{is} y_{iis}  + \sum_{j\in \Lambda_i^{+}\setminus\{i\}} d_{js} y_{ijs}  \quad\quad i \in\mathcal{I},\, s \in \mathcal{S} \label{stochastic_LB_demand_served}\\
&u_{js} \geq u_{is} + \sum_{h \in \Lambda_j^{+}} d_{hs} y_{jhs}  - l^{max} \Big(1 - \sum_{p\in\mathcal{P}} x_{ijp} \Big) \quad\quad i,j\in\mathcal{I},\,i\ne j,\, s\in\mathcal{S}\label{stochastic_track_of_u}\\
& \tau_{js} \geq \tau_{is}+t_{ijp}x_{ijp} -\Bar{t}(1-x_{ijp})+ \notag\\&\quad\quad +\sum_{h\in\Lambda_j^{+}} (\hat{t}_{j}+\beta t_{jhp})y_{jhs}\quad\quad p\in\mathcal{P},\, i,j\in\mathcal{I},\,i\ne j,\,s\in\mathcal{S}\label{stochastic_large_time_windows1}\\
& \tau_{is} \geq \sum_{p\in\mathcal{P}}\sum_{j\in\mathcal{I}^*} t_{jip}x_{jip}\quad\quad i\in\mathcal{I},\,s\in\mathcal{S}\label{stochastic_large_time_windows2}\\
&\tau_{is} + \sum_{p\in\mathcal{P}} t_{i0p}z_{ip} \leq \bar{t} \quad\quad i\in\mathcal{I},\,s\in\mathcal{S}\label{stochastic_large_time_windows3}\\
& z_{ip},\,v_{ip}\in\{0,1\} \quad\quad p\in\mathcal{P},\,i\in\mathcal{I}\label{stochastic_variables_1}\\
&x_{ijp}\in\{0,1\} \quad\quad p\in\mathcal{P},\,i\in\mathcal{I}^\ast,\,j\in\mathcal{I}\setminus\{i\}\label{stochastic_variables_2}\\
&y_{ijs},w_{is},\,u_{is},\,\tau_{is}\in\mathbb{R}^+ \quad\quad i\in\mathcal{I}, j\in\mathcal{I}, s\in\mathcal{S}.\label{stochastic_variables_3}
\end{align}

The objective function \eqref{stochastic_objective_function} minimizes the expected total delivery cost. The first two terms account for the first-stage decisions, namely the fixed deployment costs associated with fleet composition and the variable travel costs of the consistent routes. The last term represents the expected second-stage cost, including the cost of approximate recourse actions and the penalties associated with outsourced or postponed deliveries.

Constraints \eqref{stochastic_arcs=trips}-\eqref{stochastic_one_customer_one_type} define the first-stage fleet composition and consistent routing decisions. Specifically, constraints \eqref{stochastic_arcs=trips} ensure that the number of routes departing from the Distribution Center equals the number of vehicles deployed for each vehicle type. Constraints \eqref{stochastic_nr_vehicles_equal_trips_from_depot} and \eqref{stochastic_outgoing_arc_from_not_end_customers} enforce flow conservation along each consistent route, whereas constraints \eqref{stochastic_one_customer_one_type} guarantee that each customer is assigned to at most one first-stage route.

Constraints \eqref{stochastic_demand_servable only_if_i_is_served}-\eqref{stochastic_demand_satisfaction} model the second-stage recourse actions and demand satisfaction. Constraints \eqref{stochastic_demand_servable only_if_i_is_served} allow recourse actions to originate only from customers belonging to the first-stage routing plan. Since $i\in\Lambda_i^+$ for every $i\in\mathcal{I}$, these constraints also apply to customers directly served along the consistent routes. Consequently, the recourse variables $y_{ijs}$ are naturally bounded between 0 and 1. Constraints \eqref{stochastic_demand_satisfaction} ensure that the realized demand of every customer is completely satisfied under each scenario, either through the planned delivery route, through recourse actions performed from neighboring customers, or by outsourcing. As a consequence, the outsourcing variables $w_{is}$ are also bounded between 0 and 1.

Constraints \eqref{stochastic_UB_load_capacity}-\eqref{stochastic_track_of_u} regulate vehicle load capacities and track the cumulative delivered quantity along each route. Constraints \eqref{stochastic_UB_load_capacity} impose the load-capacity limits of the vehicles assigned to the first-stage routes. Constraints \eqref{stochastic_LB_demand_served} establish a lower bound on the quantity delivered at each visited customer, accounting for both planned deliveries and recourse deliveries performed from that customer. Finally, constraints \eqref{stochastic_track_of_u} recursively update the cumulative delivered quantity along each route. Whenever customer $j$ is visited immediately after customer $i$, the cumulative quantity delivered up to customer $j$ equals the cumulative quantity delivered up to customer $i$, plus the demand served at customer $j$, including any recourse deliveries originating from $j$. Otherwise, the big-$M$ term deactivates the constraint.

Constraints \eqref{stochastic_large_time_windows1}-\eqref{stochastic_large_time_windows3} enforce the working-time limitations of the delivery routes. Constraints \eqref{stochastic_large_time_windows1} recursively compute the cumulative operating time along each route by accounting for travel times, service times, and the additional time required to perform recourse actions. As in the previous capacity constraints, the big-$M$ term deactivates the constraint whenever customer $j$ does not immediately follow customer $i$. Constraints \eqref{stochastic_large_time_windows2} initialize the cumulative operating time at each visited customer, while constraints \eqref{stochastic_large_time_windows3} ensure that every vehicle completes its route and returns to the Distribution Center before the end of the working day.

Finally, constraints \eqref{stochastic_variables_1}-\eqref{stochastic_variables_3} define the domains of the first-stage and second-stage decision variables, thereby completing the formulation.

\subsubsection{Valid Inequalities}\label{sec:valid_inequalities}

In this section, the mathematical formulation of $\mathcal{M}$ is strengthened by introducing valid inequalities that tighten the feasible region, reduce the solution space, and accelerate convergence. 
These inequalities do not eliminate any feasible integer solution when added to the original formulation.

The following valid inequalities are derived for model $\mathcal{M}$.

\begin{prop}
For each scenario $s \in \mathcal{S}$, the inequality
\begin{linenomath}
\begin{equation}\label{valid_in}
    \sum_{i \in \mathcal{I}} \sum_{p \in \mathcal{P}} l_p z_{ip} 
    \geq 
    \sum_{i \in \mathcal{I}} \sum_{j \in \Lambda_i^+} d_{js} y_{ijs}
\end{equation}
\end{linenomath}
is valid for the feasible region of $\mathcal{M}$.
\end{prop}

\noindent\textit{Proof.}
Consider any feasible solution of model $\mathcal{M}$ and fix a scenario $s\in\mathcal{S}$. 
For each vehicle route selected in the first stage, let $p$ be the associated vehicle type and let $i$ be the last customer visited by that route. By definition, this route contributes $l_p z_{ip}$ to the left-hand side of \eqref{valid_in}. Hence, the left-hand side represents the total load capacity of the vehicles deployed in the first-stage routing plan.

In the second stage, the quantity
\[
\sum_{i \in \mathcal{I}} \sum_{j \in \Lambda_i^+} d_{js} y_{ijs}
\]
represents the total amount of realized demand served by the company fleet under scenario $s$. Indeed, for each customer $j$, the term $d_{js}y_{ijs}$ is the portion of the realized demand of $j$ served through customer $i$, either directly when $i=j$ or through a recourse action when $j\in\Lambda_i^+$ and $i\neq j$. Demand not served by the fleet is instead captured by the outsourcing variable $w_{js}$ and therefore does not require vehicle capacity.

Since constraints \eqref{stochastic_UB_load_capacity}-\eqref{stochastic_track_of_u} ensure that the cumulative quantity delivered along each selected route does not exceed the capacity of the vehicle operating that route, the total amount of demand served by all fleet vehicles in scenario $s$ cannot exceed the total capacity of the selected fleet. Therefore,
\[
\sum_{i \in \mathcal{I}} \sum_{j \in \Lambda_i^+} d_{js} y_{ijs}
\leq
\sum_{i \in \mathcal{I}} \sum_{p \in \mathcal{P}} l_p z_{ip},
\]
which proves the validity of \eqref{valid_in} for every $s\in\mathcal{S}$.
\hfill\qed

The above inequalities are included in the mathematical formulation of $\mathcal{M}$, where they significantly improve the computational performance of the model, as discussed in Section~\ref{sec:small_results}.

\subsection{A Stochastic Two-Stage Mixed Integer Linear Program - path based reformulation}\label{sec:PB_ref}

In this section, we present a Path-Based ($PB$) reformulation designed to address the computational challenges encountered when solving $\mathcal{M}$ on large-scale instances (see Section \ref{sec:large_results}). We denote this reformulation as $\mathcal{M}_{PB}$. The approach follows a well-established modeling strategy based on composite path variables, which has been shown to offer a more computationally tractable representation of routing decisions (see, e.g., \cite{bruglieri2019path}). A known challenge of path-based formulations is that the number of feasible routes may be exponentially large. To overcome this issue, several authors rely on Column Generation procedures that iteratively generate promising routes (see, e.g., \cite{dayarian2015column}).

In contrast, in this work we preprocess a pool of high-quality candidate routes, generated using a state-of-the-art deterministic CVRP model, as described in Section \ref{sec:path_generations}, and provide them as input to the optimization model.

Let $\mathcal{R}$ denote the set of preprocessed feasible routes. For each $r\in\mathcal{R}$, let $\mathcal{I}_r\subseteq \mathcal{I}^*$ be the set of demand nodes visited by route $r$, and for each $i\in\mathcal{I}^*$ let $\mathcal{R}_i\subseteq\mathcal{R}$ denote the subset of routes that include $i$. To account for fleet heterogeneity, for each vehicle type $p\in\mathcal{P}$ we define $\mathcal{R}_p\subseteq\mathcal{R}$ as the set of routes whose total length does not exceed the driving range $\Delta_p$.

We introduce binary variables $\psi_{rp}$ that take value 1 if route $r\in\mathcal{R}_p$ is assigned to a vehicle of type $p$, and 0 otherwise. Since each vehicle performs exactly one route, $\psi_{rp}$ also determines the number of vehicles of type $p$ used in the fleet. If route $r$ is operated by a vehicle of type $p$, a cost $c_{rp}$ is incurred, computed as the product of the route length and the operating cost $\omega_{\delta_p}$, i.e., $c_{rp}=\sum_{i,j\in\mathcal{I}r}\omega_{\delta_p}\delta_{ij}$ for all $p\in\mathcal{P}$ and $r\in\mathcal{R}_p$.

Similarly, for each $r$ and $p$, we define $\tilde{t}_{rp}$ as the travel time of route $r$ when performed by a vehicle of type $p$, computed as $\tilde{t}_{rp}=\sum_{i,j\in\mathcal{I}r}\delta_{ij}/v_p$, where $v_p$ is its average speed. In the second stage, recourse variables $y_{ijrs}$ assign recourse deliveries to the route from which they are performed.

We observe that since routes are directly represented, load, time, and fleet composition variables from $\mathcal{M}$ are no longer required, substantially simplifying the formulation.

All new sets, parameters, and decision variables introduced for this reformulation, are resumed below for clarity.\\

\noindent\emph{\underline{Sets:}}\\
$\mathcal{R}=\{1,\ldots,R\}:\,$ set of preprocessed feasible routes;\\
$\mathcal{R}_p:\,$ set of routes suitable for vehicles of type $p\in\mathcal{P}$;\\
$\mathcal{R}_i:\,$ set of routes that includes node $i\in\mathcal{I}^*$;\\
$\mathcal{I}_r:\,$ set of nodes in $\mathcal{I}^*$ included in route $r\in\mathcal{R}$.\\

\noindent\emph{\underline{Deterministic parameters:}}\\
$c_{rp}:\,$ operating cost of route $r\in\mathcal{R}_p$ when traveled by vehicles of type $p\in\mathcal{P}$;\\
$\tilde{t}_{rp}:\,$ time spent by vehicles of type $p\in\mathcal{P}$ on route $r\in\mathcal{R}_p$.\\

\noindent\emph{\underline{Decision variables:}}\\
$\psi_{rp}:\,$ binary variables indicating if route $r\in\mathcal{R}$ is operated by a vehicle of type $p\in\mathcal{P}$;\\
$y_{ijrs}:\,$ non-negative variables indicating the percentage amount of demand of $j\in\mathcal{I}$ that is served from $i\in\mathcal{I}$, included in route $r\in\mathcal{R}$, under scenario $s\in\mathcal{S}$.\\

The resulting stochastic two-stage mixed-integer linear model $\mathcal{M}_{PB}$ is:
\begin{align}
\min  \quad & \sum_{p\in \mathcal{P}}\sum_{r\in\mathcal{R}_p} f_p\psi_{rp} + \sum_{p\in\mathcal{P}}\sum_{r\in\mathcal{R}_p} c_{rp}\psi_{rp} + \sum_{s\in\mathcal{S}}\pi_s \Bigg(\sum_{r\in\mathcal{R}}\sum_{i\in\mathcal{I}}\sum_{j\in\Lambda_{i}^{+}}\beta c^{rec}_{ij} y_{ijrs} + \sum_{i \in \mathcal{I}} \gamma d_{is} w_{is}\Bigg) \label{ref_objective_function} \\
\text{s.t.}\quad & 
\psi_{rp}=0 \quad\quad p\in\mathcal{P},\, r\in\mathcal{R}\setminus\mathcal{R}_p\label{ref_r_not_activated_if_not_linked_to_p}\\
& \sum_{p\in\mathcal{P}} \psi_{rp}\leq 1 \quad\quad r\in\mathcal{R}\label{ref_one_routes_atmost}\\
& \sum_{p\in\mathcal{P}}\sum_{r\in\mathcal{R}_{i}}\psi_{rp}\leq 1 \quad\quad i\in\mathcal{I}\label{ref_arc_included_only_one}\\
&\sum_{r\in\mathcal{R}}\sum_{i \in\mathcal{I}_r\cap \Lambda_j^{-}\setminus{\{0\}}}  y_{ijrs} + w_{js} = 1 \quad\quad j \in \mathcal{I},\, s \in \mathcal{S}\label{ref_demand_satisfaction}\\
&y_{ijrs} \leq \sum_{p\in\mathcal{P}}\psi_{rp} \quad\quad r\in\mathcal{R},\,i \in \mathcal{I}_r\setminus\{0\},\, j \in \Lambda_i^{+},\, s \in \mathcal{S}\label{ref_demand_servable only_if_i_is_in_route}\\
& \sum_{i\in\mathcal{I}_r\setminus\{0\}}\sum_{j\in\Lambda_i^+}d_{js}y_{ijrs} \leq \sum_{p\in\mathcal{P}} l_p\psi_{rp} \quad\quad r\in\mathcal{R},\, s\in\mathcal{S}\label{ref:capacity_constraints}\\
& \tilde{t}_{rp}+\sum_{i\in\mathcal{I}_r\setminus\{0\}}\sum_{j\in\Lambda_i^+\setminus{\{i\}}}(\hat{t}_{j}+\beta t_{ijp}y_{ijrs})\leq \bar{t} \quad\quad p\in\mathcal{P},\,r\in\mathcal{R}_p,\, s\in\mathcal{S}\label{ref:time_windows_contraint}\\
& \psi_{rp}\in\{0,1\} \quad\quad p\in\mathcal{P},\,r\in\mathcal{R}\label{ref_variables_1}\\
&y_{ijrs},\,w_{is}\in\mathbb{R}^+ \quad\quad i\in\mathcal{I},\, j\in\mathcal{I},\,r\in\mathcal{R},\,s\in\mathcal{S}.\label{ref_variables_3}
\end{align}

In $\mathcal{M}_{PB}$, the objective function \eqref{ref_objective_function} minimizes the total expected service cost. The first-stage components capture fleet-related costs: the first term represents the fixed cost of using vehicles of different types, while the second term accounts for route operating costs, which depend on the vehicle type assigned to each route. The third term represents the expected second-stage recourse costs, weighted by scenario probabilities, analogous to the cost structure in $\mathcal{M}$.

Constraints \eqref{ref_r_not_activated_if_not_linked_to_p} ensure that a route can be assigned only to vehicle types capable of operating it. Constraints \eqref{ref_one_routes_atmost} impose that each route is selected at most once, while constraints \eqref{ref_arc_included_only_one} guarantee that each demand node is visited by at most one route. Constraints \eqref{ref_demand_satisfaction} are equivalent to \eqref{stochastic_demand_satisfaction} in $\mathcal{M}$, ensuring that the demand of each node is either satisfied along a selected route or covered through recourse. Constraints \eqref{ref_demand_servable only_if_i_is_in_route} restrict recourse deliveries to originate only from nodes belonging to an activated route.

Capacity feasibility is enforced by \eqref{ref:capacity_constraints}, which limit, for each scenario, the total amount of demand served along a route (both directly and via recourse) to the load capacity of the vehicle assigned to that route. Similarly, constraints \eqref{ref:time_windows_contraint} ensure temporal feasibility by requiring that the total travel and service time (including additional time induced by recourse actions) does not exceed the maximum allowed route duration. Finally, \eqref{ref_variables_1}-\eqref{ref_variables_3} specify the domain of the decision variables.

\subsection{Characterization of Approximate Recourse Actions}\label{sec:recourse_actions_properties}

This section characterizes the approximate recourse mechanism embedded in the proposed stochastic optimization models $\mathcal{M}$ and $\mathcal{M}_{PB}$.

Let $\mathcal{R}^{\ast}=\{R^{\ast,1},\ldots,R^{\ast,\bar{N}}\}$ denote the set of first-stage consistent routes obtained from the optimal solution of either $\mathcal{M}$ or $\mathcal{M}_{PB}$, where $\bar{N}$ is the number of vehicles in the fleet. Each route $R^{\ast,n}\in\mathcal{R}^{\ast}$ is a sequence of nodes visited by the corresponding vehicle, starting and ending at the central depot $0$. Then, we write $i \in R^{\ast,n}$ to indicate that customer $i$ is included in route $R^{\ast,n}$.

Depending on the realized customer demand, the residual load capacity, and the available operating time of the fleet, approximate recourse actions may be implemented in different ways, as illustrated in Figure~\ref{fig:types_of_approximate_ra}.

\begin{figure}[h!]
    \centering
    \subfloat[\footnotesize Adjusted consistent routes with additional demand nodes served through complete approximate recourse actions.\label{(a)}]{\begin{tikzpicture}[roundnode/.style={circle, draw=black!60, fill=black!5, very thick, minimum size=0.3mm},squarednode/.style={rectangle, draw=black!60, fill=black!5, very thick, minimum size=8mm}]
    \begin{scope}[scale=0.45]
        \node[squarednode] (0) at (0,0) {0};
        \node[roundnode] (1) at (5,1.5) {\tiny 1};
        \node[roundnode] (2) at (1,5) {\tiny 2};
        \node[roundnode] (3) at (-3,3) {\tiny 3};
        \node[roundnode] (4) at (-1.5,-3) {\tiny 4};
        \node[roundnode] (5) at (3.5,-2) {\tiny 5};
        \node[roundnode] (6) at (2,6.2) {\tiny 6};
        \node[roundnode] (7) at (-2,-4.2) {\tiny 7};
        \node[roundnode] (8) at (-1.5,6.2) {\tiny 8};
        \node[roundnode] (9) at (1,-1.8) {\tiny 9};
        \node[roundnode] (10) at (-3.5,0.5) {\tiny 10};
        \node[roundnode] (11) at (3.5,3.5) {\tiny 11};

        \path (0) edge[->,ultra thick, dotted, bend right,blue]node[right]{} (2);
        \path (2) edge[->,ultra thick, dotted, bend right,blue]node[above]{}  (3);
        \path (3) edge[->,ultra thick, dotted, bend right,blue]node[right]{} (0);
        \path (2) edge[<->,ultra thick,orange] node[above right]{} (11);
        \path (3) edge[<->,ultra thick,orange] node[above left]{} (10);

        \path (0) edge[->,ultra thick, dash dot, bend left,red]node[below]{}(1);
        \path (1) edge[->,ultra thick, dash dot, bend left,red]node[right]{} (5);
        \path (5) edge[->,ultra thick, dash dot, bend left,red]node[below]{} (4);
        \path (4) edge[->,ultra thick, dash dot, bend left,red] node[above left]{} (0);
        \path (4) edge[<->,ultra thick,orange]node[above left]{} (9);


    \end{scope}
    \end{tikzpicture}}\quad\subfloat[\footnotesize Adjusted consistent routes with additional demand nodes served through complete and split approximate recourse actions.\label{(b)}]{\begin{tikzpicture}[roundnode/.style={circle, draw=black!60, fill=black!5, very thick, minimum size=0.3mm},squarednode/.style={rectangle, draw=black!60, fill=black!5, very thick, minimum size=8mm}]
    \begin{scope}[scale=0.45]
        \node[squarednode] (0) at (0,0) {0};
        \node[roundnode] (1) at (5,1.5) {\tiny 1};
        \node[roundnode] (2) at (1,5) {\tiny 2};
        \node[roundnode] (3) at (-3,3) {\tiny 3};
        \node[roundnode] (4) at (-1.5,-3) {\tiny 4};
        \node[roundnode] (5) at (3.5,-2) {\tiny 5};
        \node[roundnode] (6) at (2,6.2) {\tiny 6};
        \node[roundnode] (7) at (-2,-4.2) {\tiny 7};
        \node[roundnode] (8) at (-1.5,6.2) {\tiny 8};
        \node[roundnode] (9) at (1,-1.8) {\tiny 9};
        \node[roundnode] (10) at (-3.5,0.5) {\tiny 10};
        \node[roundnode] (11) at (3.5,3.5) {\tiny 11};

        \path (0) edge[->,ultra thick, dotted, bend right,blue]node[right]{} (2);
        \path (2) edge[->,ultra thick, dotted, bend right,blue]node[above]{}  (3);
        \path (3) edge[->,ultra thick, dotted, bend right,blue]node[right]{} (0);
        \path (2) edge[<->,ultra thick,orange] node[above right]{} (11);
        \path (1) edge[<->,ultra thick,orange] node[above right]{} (11);
        \path (3) edge[<->,ultra thick,orange] node[above left]{} (10);

        \path (0) edge[->,ultra thick, dash dot, bend left,red]node[below]{}(1);
        \path (1) edge[->,ultra thick, dash dot, bend left,red]node[right]{} (5);
        \path (5) edge[->,ultra thick, dash dot, bend left,red]node[below]{} (4);
        \path (4) edge[->,ultra thick, dash dot, bend left,red] node[above left]{} (0);
        \path (4) edge[<->,ultra thick,orange]node[above left]{} (9);
        \path (5) edge[<->,ultra thick,orange]node[above right]{} (9);


    \end{scope}
    \end{tikzpicture}}\quad
    \subfloat[\footnotesize Adjusted consistent routes with additional demand nodes served through split approximate recourse actions or externalization.\label{(c)}]{\begin{tikzpicture}[roundnode/.style={circle, draw=black!60, fill=black!5, very thick, minimum size=0.3mm},squarednode/.style={rectangle, draw=black!60, fill=black!5, very thick, minimum size=8mm}]
    \begin{scope}[scale=0.45]
        \node[squarednode] (0) at (0,0) {0};
        \node[roundnode] (1) at (5,1.5) {\tiny 1};
        \node[roundnode] (2) at (1,5) {\tiny 2};
        \node[roundnode] (3) at (-3,3) {\tiny 3};
        \node[roundnode] (4) at (-1.5,-3) {\tiny 4};
        \node[roundnode] (5) at (3.5,-2) {\tiny 5};
        \node[roundnode] (6) at (2,6.2) {\tiny 6};
        \node[roundnode] (7) at (-2,-4.2) {\tiny 7};
        \node[roundnode] (8) at (-1.5,6.2) {\tiny 8};
        \node[roundnode] (9) at (1,-1.8) {\tiny 9};
        \node[roundnode] (10) at (-3.5,0.5) {\tiny 10};
        \node[roundnode] (11) at (3.5,3.5) {\tiny 11};

        \path (0) edge[->,ultra thick, dotted, bend right,blue]node[right]{} (2);
        \path (2) edge[->,ultra thick, dotted, bend right,blue]node[above]{}  (3);
        \path (3) edge[->,ultra thick, dotted, bend right,blue]node[right]{} (0);
        \path (2) edge[<->,ultra thick,orange] node[above right]{} (11);
        \path (3) edge[<->,ultra thick,orange] node[above left]{} (10);

        \path (0) edge[->,ultra thick, dash dot, bend left,red]node[below]{}(1);
        \path (1) edge[->,ultra thick, dash dot, bend left,red]node[right]{} (5);
        \path (5) edge[->,ultra thick, dash dot, bend left,red]node[below]{} (4);
        \path (4) edge[->,ultra thick, dash dot, bend left,red] node[above left]{} (0);
        \path (0) edge[<->,ultra thick,green]node[left]{} (10);
        \path (0) edge[<->,ultra thick,green]node[left]{} (5);
        \path (1) edge[<->,ultra thick,orange]node[above right]{} (11);


    \end{scope}
    \end{tikzpicture}}\\
    \caption{\small Examples of approximate recourse actions on consistent routes for a delivery fleet comprising two vehicle types. Double arrows denote approximate recourse actions with $\beta = 2$. Orange arrows indicate complete or split approximate recourse actions, whereas green arrows represent demand externalization.}   \label{fig:types_of_approximate_ra}
\end{figure}

Consider a scenario $s\in\mathcal{S}$ and suppose that additional demand arises at customer $h\in\mathcal{I}$. Without loss of generality, the following discussion focuses on additional customers, although the same mechanism also applies whenever the realized demand of a customer already assigned to a consistent route exceeds the residual capacity of its associated vehicle. First, we assume that the demand of customer $h$ can be served only by vehicles operating routes that contain at least one customer $i$ such that $h\in\Lambda_i^{+}$. Among all such feasible routes, let $R^{\ast,\bar n}_{h}\in\mathcal{R}^*$ denote the one that is closest to customer $h$, where
$\displaystyle \bar n=\argmin_{n\in\{1,\ldots,\bar N\},\, i\in R^{\ast,n}}\delta_{ih}$. For notational simplicity, we assume that $\bar n$ is unique; if multiple routes attain the same minimum distance, the results stated below apply to each of them.

The proposed recourse mechanism distinguishes the following three operational situations:

\begin{itemize}
    \item[-]
     \textbf{Complete approximate recourse action:} the realized demand of customer $h$ is entirely satisfied by a single vehicle. This occurs when the vehicle operating route $R^{\ast,\bar n}_{h}$, or more generally any route $R^{\ast,n}\in\mathcal{R}^*$ containing a customer $i$ such that $h\in\Lambda_i^{+}$, has sufficient residual load capacity to accommodate the realized demand and the resulting operating time does not exceed the working-time limit $\bar t$ (see Figure~\ref{(a)}).
    \item[-] \textbf{Split approximate recourse action:} if no single vehicle can satisfy the entire realized demand, but multiple vehicles operating routes in $\mathcal{R}^{\ast}$ satisfy the neighborhood condition, residual capacity, and operating-time requirements, then the demand of customer $h$ is split among these vehicles through independent approximate recourse actions (see Figure~\ref{(b)}). A detailed example of split approximate recourse actions is also reported in \hyperref[sec:appendix_a]{Appendix A}.
    \item[-] \textbf{Externalized demand:} if neither complete nor split approximate recourse actions are feasible because the available fleet resources are insufficient, the remaining realized demand is fulfilled through an external delivery service (see Figure~\ref{(c)}). In this case, $0<w_{hs}\leq1$, where $w_{hs}$ denotes the proportion of the realized demand of customer $h$ that is outsourced under scenario $s$.
\end{itemize}

The above classification summarizes the approximate recourse mechanism embedded in both stochastic optimization models. Whenever feasible, realized demand is accommodated through complete or split approximate recourse actions, whereas demand externalization is used only when fleet resources are insufficient. Importantly, these recourse actions do not alter the first-stage routing plan, thereby preserving the consistency of the routes determined at the tactical level.

\section{Solution Methods}\label{sec:sol_meth}

In this section, we describe the solution methodology adopted for the proposed problem. Section \ref{sec:scen_red} details the scenario reduction procedure adopted to limit the number of scenarios. Section \ref{sec:KS_approach} introduces a Kernel Search (KS)-based approach aimed at improving the computational performance of the formulation when applied to large-scale real-world instances.

\subsection{Scenario Reduction}\label{sec:scen_red}

The performance of model $\mathcal{M}$ is strongly influenced by the cardinality of the scenario set $\mathcal{S}$, as discussed in Section \ref{sec:small_results}. Larger scenario sets increase the representativeness of the stochastic demand but also lead to significantly higher computational effort, even for small-size instances. To balance accuracy and tractability, we adopt a scenario reduction technique adapted from \cite{dupavcova2003scenario} and \cite{heitsch2003scenario}, with the aim of retaining the most informative scenarios while keeping the size of $\mathcal{S}$ manageable.

The main concepts of the implemented scenario reduction procedure is summarized below (for extended details, refer to \cite{dupavcova2003scenario} and \cite{heitsch2003scenario}).

Let $\Pi$ be a discrete probability distribution carried by finitely many scenarios $d_{is}$, with $i\in\mathcal{I}^*$, $s\in\mathcal{S}$, whose probabilities are $\pi_{s}\geq0$ with $\sum_{s\in\mathcal{S}}\pi_s=1$. Hence, $\Pi=\sum_{s\in\mathcal{S}}\pi_s \delta_{d_s}$ where $d_s$ represents the random vector of customers demand, that is $d_s=\left[d_{is}\right]_{i\in\mathcal{I}^*}$, for $s\in\mathcal{S}$, and $\delta_{d_s}$ represents the Dirac measure at $d_s$.
Let $\mathcal{S}^{RED}\subset\mathcal{S}$ be the set of the retained scenarios and $\mathcal{S}^{DEL}=\mathcal{S}\setminus\mathcal{S}^{RED}$ the set of the deleted scenarios. Once that $\mathcal{S}^{RED}$ is identified, a new probability measure $\Pi^{RED}=\sum_{s\in\mathcal{S}^{RED}}\pi_s^{RED} \delta_{d_s}$ can be defined having scenarios $d_s$ whose probabilities are $\pi_s^{RED}\geq 0$ with $s\in\mathcal{S}^{RED}$ and $\sum_{s\in\mathcal{S}^{RED}}\pi_s^{RED}=1$. In addition, let $\pi^{RED}=\left[\pi_{s}^{RED}\right]_{s\in\mathcal{S}^{RED}}$. 
The new probability measure will be defined by deleting scenarios $d_s$ for all $s\in\mathcal{S}^{DEL}$ and by assigning new probabilities $\pi_s^{RED}$ to each scenario $d_s$ with $s\in\mathcal{S}^{RED}$ without modifying the values of $d_{is}$ in the reduced scenarios. The aim is to determine the set $\mathcal{S}^{RED}\subset\mathcal{S}$ of given cardinality $|\mathcal{S}^{RED}|$ and the probability measure $\Pi^{RED}$ such that the probability distance $D(\mathcal{S}^{RED};\pi^{RED})$ is minimized, where
{\small\begin{equation*}\label{distance_prob}
    D(\mathcal{S}^{RED};\pi^{RED})\coloneqq\min\Bigg\{\sum_{\substack{s\in\mathcal{S}\\\tilde{s}\in\mathcal{S}^{RED}}}\|d_s-d_{\tilde{s}}\|_2 \cdot\eta_{s\tilde{s}}: \eta_{s\tilde{s}}\geq 0, \sum_{s\in\mathcal{S}}\eta_{s\tilde{s}}=\pi_{\tilde{s}}^{RED}, \sum_{\tilde{s}\in\mathcal{S}^{RED}}\eta_{s\tilde{s}}=\pi_s\Bigg\}
\end{equation*}}
represents a linear transportation problem, and $\|\cdot\|_2$ is the Euclidean norm.
\noindent If $\mathcal{S}^{RED}$ is given, 
according to \cite{dupavcova2003scenario}, then $D(\mathcal{S}^{RED};\pi^{RED})$ can be reformulated as
\begin{equation}\label{red_formula}
    D^{RED}=\sum_{\Tilde{s}\not\in\mathcal{S}^{RED}} \pi_{\Tilde{s}}\min_{s\in\mathcal{S}^{RED}}\|d_{\Tilde{s}}- d_s\|_2,
\end{equation}
where $D^{RED}:=\min \Big\{D(\mathcal{S}^{RED}; \pi^{RED}) : \sum_{s\in\mathcal{S}^{RED}} \pi_s^{RED} =1, \pi_s^{RED}\geq 0,\,s\in\mathcal{S}^{RED}\Big\}$. In addition, the minimum is attained at
\begin{equation}\label{new_probabilities}
    \pi_s^{RED} = \pi_s + \sum_{\Tilde{s}\in\mathcal{S}_s^{DEL}}\pi_{\Tilde{s}}, \quad \forall s\in \mathcal{S}^{RED},
\end{equation}
with $\mathcal{S}_s^{DEL}=\{\Tilde{s}\in\mathcal{S}^{DEL} : s=s(\Tilde{s})\}$ and $s(\Tilde{s})\in\argmin_{s\in\mathcal{S}^{RED}} \|d_{\Tilde{s}}- d_s\|_2$ for all $\Tilde{s}\not\in\mathcal{S}^{RED}$. Equation (\ref{new_probabilities}), called \textit{optimal redistribution rule}, specifies that the new probabilities of each preserved scenario is computed by adding to its former probability the sum of all the probabilities of the deleted scenarios that are closest to it.

The problem is then how to optimally choose the index set $\mathcal{S}^{RED}$ for the scenario reduction with $|\mathcal{S}^{RED}|$ fixed. This leads to consider (\ref{red_formula}) as a set-covering problem that can be formulated as a 0-1 integer program, known to be $\mathcal{NP}$-hard. Some efficient algorithms can be taken into account to solve it when $|\mathcal{S}^{RED}| = 1$ or $|\mathcal{S}^{RED}| = |\mathcal{S}|-1$. Among them, we take into account the \textit{Forward Selection (FS)} algorithm and, in the specific, we consider the \textit{Fast Forward Selection (FFS)} algorithm, provided and described in details in \cite{heitsch2003scenario} and \cite{Gioia2023}.

\begin{figure}[h!]
    \centering
    \subfloat[Original scenario tree with $|\mathcal{S}|=50$.\label{fig:original_scen_tree}]{\includegraphics[width=.45\textwidth]{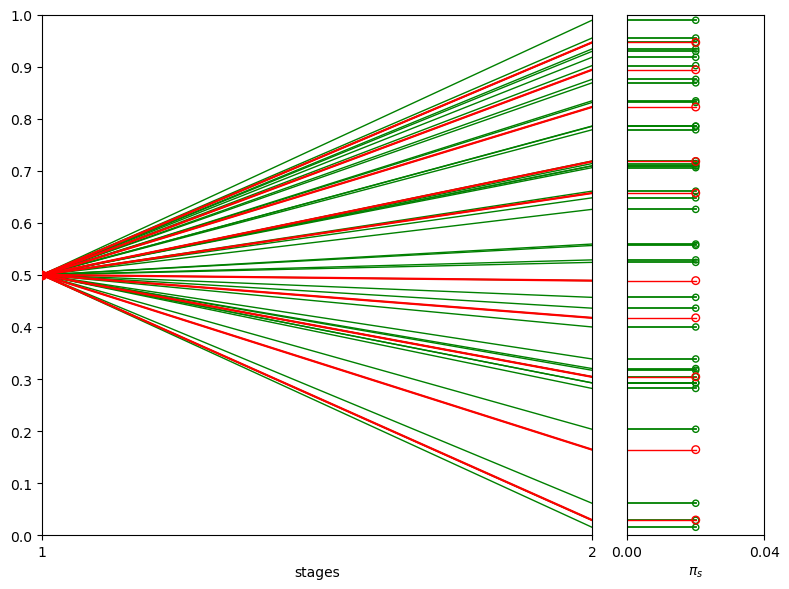}} \hspace{0.05\textwidth}
\subfloat[Reduced scenario tree with $|\mathcal{S}^{RED}|=10$.\label{fiog:red_scen_tree}]{\includegraphics[width=.45\textwidth]{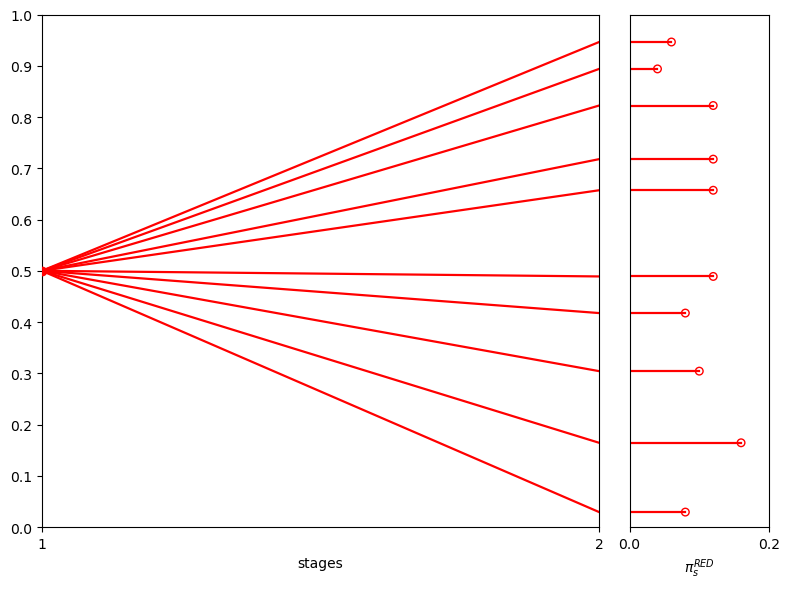}}
\caption{\small{Example of the \emph{FFS} algorithm application to a scenario tree with 50 scenarios, reduced to a scenario tree with 10 scenarios. Selected scenarios are highlighted in red, together with their new probabilities.}}
    \label{fig:scen_reduction}
\end{figure}

Figure \ref{fig:scen_reduction} illustrates the reduction process. Starting from an initial one-dimensional scenario tree with $|\mathcal{S}|=50$ equally likely scenarios, the \emph{FFS} algorithm selects a subset of $|\mathcal{S}^{RED}|=10$ scenarios (highlighted in red in Figure \ref{fig:original_scen_tree}) and assigns updated probabilities to them (Figure \ref{fiog:red_scen_tree}). The resulting reduced set provides a compact yet informative representation of the original uncertainty.

As shown in Section \ref{sec:small_results}, this scenario reduction significantly decreases computational effort while preserving solution quality, with negligible deviation from the results obtained using the full scenario set.

\subsection{Kernel Search-based Approach}\label{sec:KS_approach}

In this section, we develop a KS-based heuristic to enhance the tractability of the proposed formulations. Although the $PB$ reformulation improves performance, it remains computationally challenging for the real-world case study (see Section \ref{sec:real_case_study}), due to the large number of candidate routes generated in the preprocessing phase. To address this issue, we design a tailored Kernel Search approach that exploits the structure of $\mathcal{M}_{PB}$ by iteratively solving reduced subproblems over carefully selected subsets of routes (see \cite{angelelli2010kernel} for the general KS framework and \cite{kinene2023electric} for a recent application in large-scale network design problems). In line with these contributions, the proposed heuristic alternates between intensification, by refining a restricted promising region of the solution space, and diversification, by progressively incorporating new variables into the search. The complete iterative process is summarized below.

We initialize the kernel $\mathcal{K}$ with all elementary routes, i.e., routes serving a single customer, to ensure feasibility and to provide a structurally valid starting solution. The remaining routes in $\mathcal{R}\setminus\mathcal{K}$ are then randomly partitioned into buckets $B$ of fixed cardinality $|B|$. This yields $N=\lfloor |\mathcal{R}\setminus\mathcal{K}| / |B| \rfloor$ buckets of equal size, while the remaining routes are grouped into a final bucket. The resulting list of buckets is denoted by $\mathcal{B}$. Such a partitioning strategy allows the method to control the size of the subproblems while systematically exploring the pool of candidate routes.

The algorithm proceeds in two cycles. In each iteration of the first cycle, a bucket $B\in\mathcal{B}$ is temporarily added to $\mathcal{K}$ and $\mathcal{M}_{PB}$ is solved on the reduced set $\mathcal{R}=\mathcal{K}\cup B$. Any routes selected in the optimal solution but not yet in $\mathcal{K}$ are then permanently added to the kernel, thus expanding the promising region of the search space. This process continues until all buckets in $\mathcal{B}$ have been examined, progressively enriching the kernel with high-quality routes. The second cycle aims to reassess the routes excluded in the first phase and to mitigate potential biases introduced by the initial random partition. To this end, the set $\mathcal{R}\setminus\mathcal{K}$ is re-partitioned into new buckets $\bar{B}$, forming a new list $\bar{\mathcal{B}}$, and the same iterative update mechanism is applied. A stopping rule based on a maximum computation time $T_{max}$ and an Optimality Threshold $OT$, computed as the gap between the incumbent solution and the best lower bound obtained so far, is used to terminate the search.

In summary, the proposed KS-based heuristic systematically explores the solution space by repeatedly solving smaller and more tractable subproblems derived from the original path-based model. The kernel is progressively enriched with promising routes, ensuring monotonic improvement of the incumbent solution and fostering convergence, while maintaining computational efficiency. As shown in Section \ref{sec:large_results}, this approach delivers high-quality solutions for large-scale instances within practical computational limits. In Algorithm \ref{alg:TKS}, the full pseudocode of the proposed KS-based procedure is reported.

\begin{algorithm}[h!]
    \caption{KS-based algorithm.}
    \textbf{Input:} a set $\mathcal{R}$ of non elementary routes, a set $\mathcal{R}_{el}$ of elementary routes, bucket size $N$, maximum run time $T^{max}$, optimality threshold $OT$.\newline
	\textbf{Output:} the objective function value $M_{PB}^*$ and the solution $z_{PB}^*, \psi_{PB}^*, y_{PB}^*, w_{PB}^*$ of $\mathcal{M}_{PB}$
      \begin{algorithmic}[1]
	\Function{KS-based algorithm}{}
        \Statex \hspace*{0.5cm}\small{\texttt{\# Initialization}}
    \State Build the initial kernel: $\mathcal{K} \gets \mathcal{R}_{el}$
    \State Iteration counter: $k \gets 0$
        \While{$\textit{runtime}\leq T^{max}$ and $\textit{optimality gap}\geq OT$}
            \Statex \hspace*{0.5cm}\small{\texttt{\# First cycle}}
        \State create a set of non-kernel routes: $\mathcal{Q}=\mathcal{R}\setminus\mathcal{K}$
        \State create a buckets list: $\mathcal{B}\gets\emptyset$
                \While{$\mathcal{Q}$ is not empty}
            \State randomly select up to $N$ routes from $\mathcal{Q}$ and add them to a bucket $B$
            \State add bucket $B$ to the buckets list: $\mathcal{B}\gets\mathcal{B}\cup\{B\}$
            \State remove bucket routes from $\mathcal{Q}$: $\mathcal{Q}\gets\mathcal{Q}\setminus B$
        \EndWhile
                \For{$B\in\mathcal{B}$}
            \State solve restricted model $\mathcal{M}_{PB}(\mathcal{K}\cup B)$
            \State get the activated routes: $\mathcal{A}=\{r\in B: \psi_{rp}=1,\, p\in\mathcal{P}\}$
            \State update the kernel: $\mathcal{K}\gets \mathcal{K}\cup\mathcal{A}$
            \State store incumbent solution $(M^*_{PB},z_{PB}^*, \psi_{PB}^*, y_{PB}^*, w_{PB}^*)$
        \EndFor
                \Statex \hspace*{0.5cm}\small{\texttt{\# Second cycle}}
        \State define remaining routes: $\bar{\mathcal{Q}}=\mathcal{R}\setminus\mathcal{K}$
        \State create a new buckets list: $\bar{\mathcal{B}}\gets\emptyset$
                \While{$\bar{\mathcal{Q}}$ is not empty}
            \State randomly select up to $N$ routes from $\bar{\mathcal{Q}}$ and add them to a bucket $\bar{B}$
            \State add bucket $\bar{B}$ to the list: $\bar{\mathcal{B}}\gets\bar{\mathcal{B}}\cup\{\bar{B}\}$
            \State remove bucket routes from $\bar{\mathcal{Q}}$: $\bar{\mathcal{Q}}\gets\bar{\mathcal{Q}}\setminus \bar{B}$
        \EndWhile
                \For{$\bar{B}\in\bar{\mathcal{B}}$}
            \State solve restricted model $\mathcal{M}_{PB}(\mathcal{K}\cup \bar{B})$
            \State get the activated routes: $\bar{\mathcal{A}}=\{r\in \bar{B}: \psi_{rp}=1,\, p\in\mathcal{P}\}$
            \State update the kernel: $\mathcal{K}\gets \mathcal{K}\cup\bar{\mathcal{A}}$
            \State update incumbent solution if improved
        \EndFor
                \State update iteration counter: $k\gets k+1$
            \EndWhile	
    	\EndFunction
    \end{algorithmic}
    \label{alg:TKS}
\end{algorithm}

\section{Computational Experiments}\label{sec:real_case_study}

This section describes the empirical setting used in the numerical experiments (Section \ref{sec:empirical_settings}), which are based on a real-world dataset provided by the Italian postal company. Computational results are presented in Section \ref{sec:comput_results}, followed by managerial insights in Section \ref{sec:practical_insights}.

\subsection{Empirical Setting}\label{sec:empirical_settings}

We first introduce the real-world operational context and data (Section \ref{sec:data_and_context}). Then, we describe the procedures used to generate small and large instances (Section \ref{sec:scenario_gen_red}), and finally the route generation process (Section \ref{sec:path_generations}).

\subsubsection{Data and Context}\label{sec:data_and_context}

The case study concerns a Last Mile delivery service operated by the Italian postal company. The delivery fleet is based at a distribution center located in Trescore Balneario (IT) and serves 18 surrounding municipalities. For large-scale instances, we use a comprehensive operational dataset containing the daily volume and spatial distribution of parcel deliveries for each working day from January 2024 to April 2024. This dataset provides both the number of delivered parcels and the geographic coordinates of delivery locations.
The fleet operating at the distribution center includes two vehicle types: conventional motor vehicles (CMs), referred to as type $a$, and electric cargo bikes (ECBs), referred to as type $b$. These two vehicle types differ in load capacity, operating cost, and driving range, and are therefore modeled explicitly in the numerical experiments.

\subsubsection{Instance Generation}\label{sec:scenario_gen_red}

\paragraph{Small Instances}

A set of small synthetic instances is generated to validate the formulation $\mathcal{M}$ and to assess the impact of demand uncertainty. In line with the strategic–tactical planning scope of the problem, customers are aggregated into delivery areas (i.e., demand nodes) by partitioning the region into a hexagonal grid with cell size 500m × 500m. This approach is consistent with standard aggregation strategies commonly adopted in Last Mile routing (e.g., zip-code or arc-based aggregation; see \cite{sankaran2007solving}, \cite{wang2015vehicle}).

\begin{figure}[h!]
\centering
\subfloat[Generated customers demand with $100$ delivery requests. \label{fig:gendem1}]{\includegraphics[width=6.5 cm]
{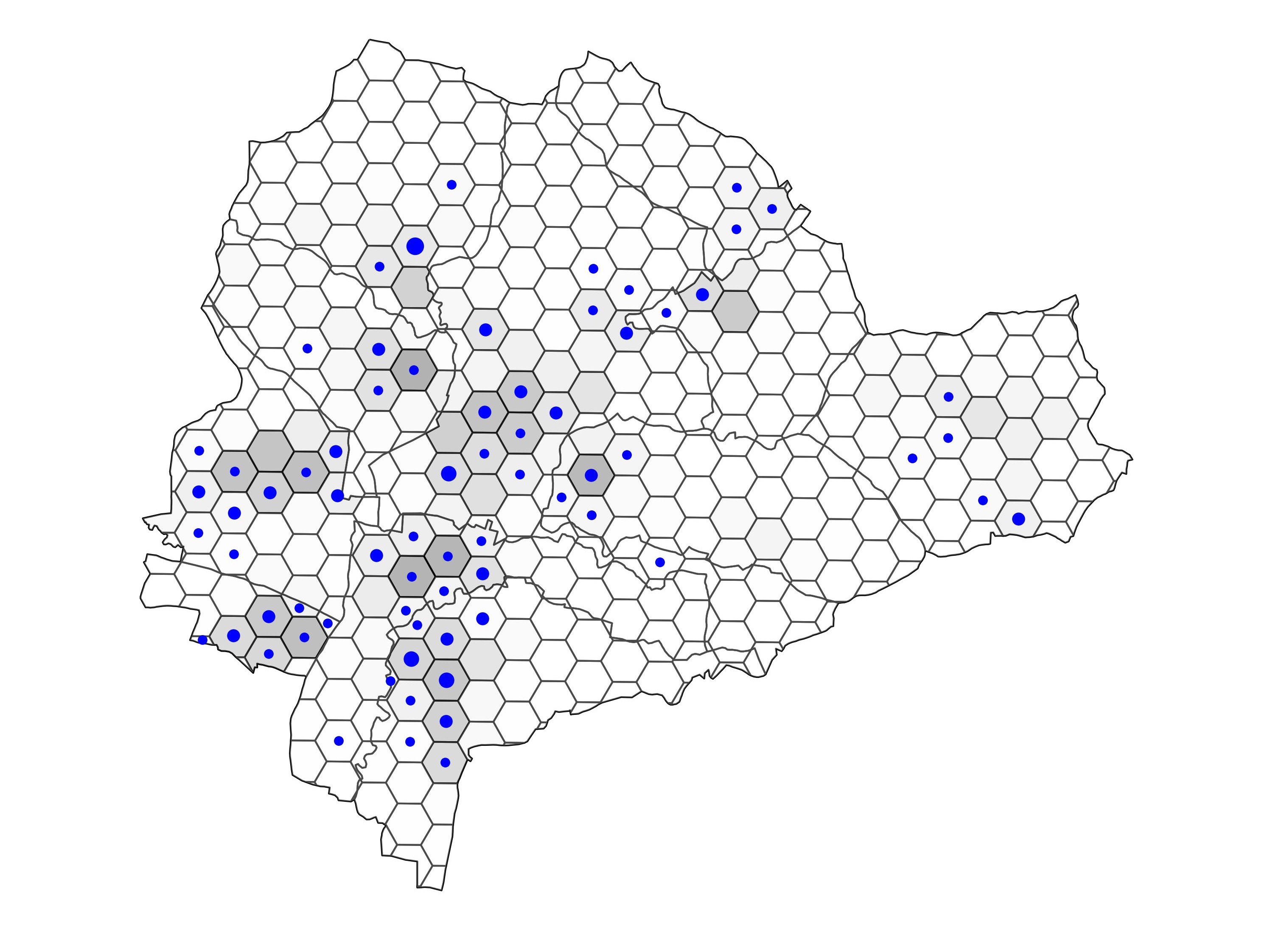}}\hspace{0.25cm}
\subfloat[Generated customers demand with $1000$ delivery requests.   \label{fig:gendem2}]{\includegraphics[width=6.5 cm]{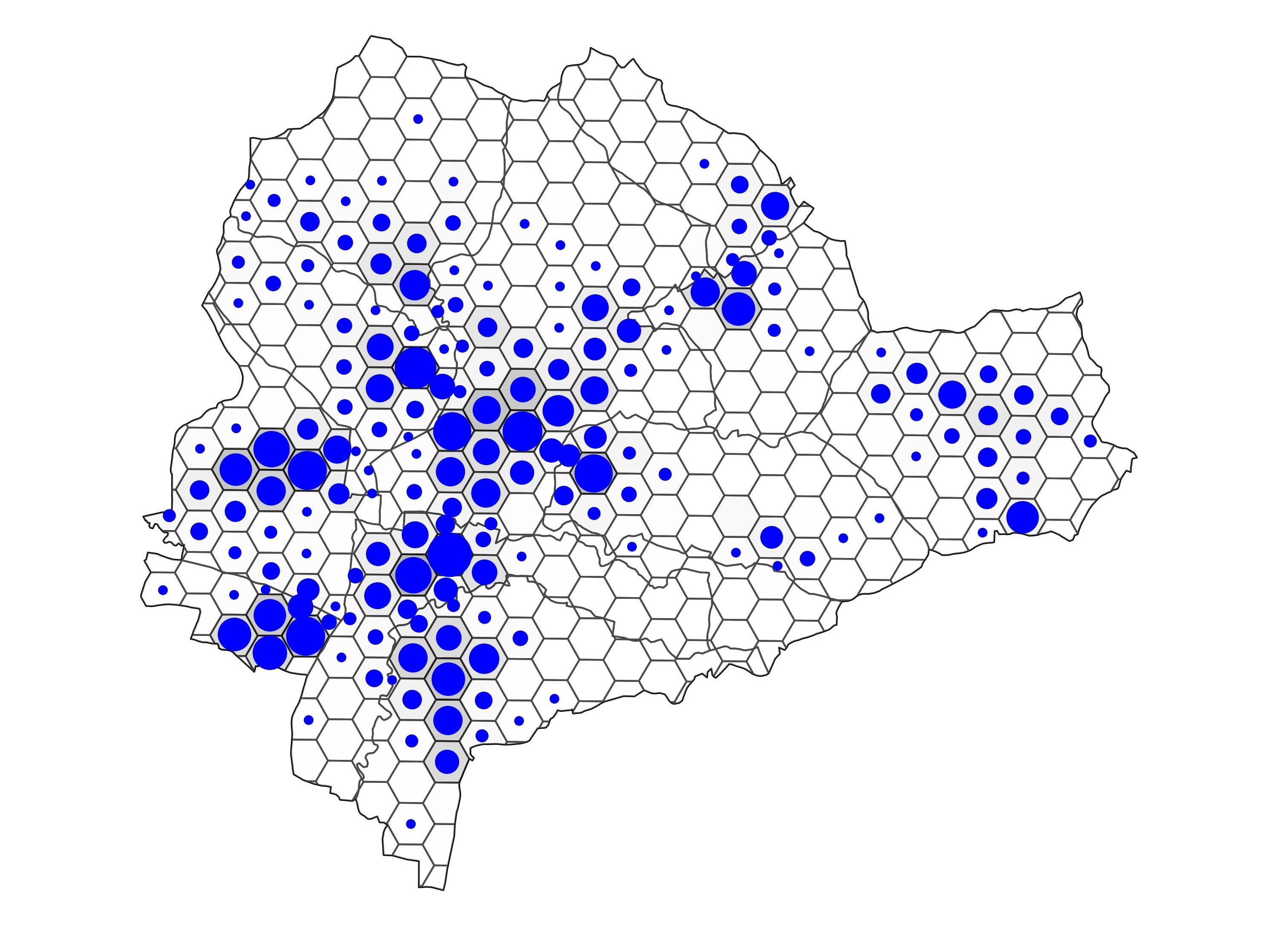}}
\caption{Demand generation based on population density. Gray shading indicates population density; marker size reflects the number of generated delivery requests in each cell.}
\label{fig:validation_plots} 
\end{figure}

Delivery requests are generated by sampling from a population density-based kernel distribution calibrated on the study area. Figure \ref{fig:validation_plots} illustrates the resulting spatial patterns of demand for two different total request volumes (100 and 1000). Darker cells indicate higher resident population density, while the size of the blue markers reflects the number of sampled delivery requests in each grid cell.

Five deterministic small instances, each containing 20 delivery requests, are produced following this procedure. Stochastic instances are derived from the deterministic ones by perturbing node demands. Let $\bar{d}_i$ denote the deterministic demand at node $i\in\mathcal{I}$, and let $\rho_i \sim \mathcal{U}(0,4)$ be a multiplicative noise factor. Then, for each scenario $s\in\mathcal{S}$ and node $i\in\mathcal{I}$, the scenario-dependent demand is defined as $d_{is} = \bar{d}_i \cdot \rho_i$. Consistent with the strategic planning level, service times are assumed deterministic and negligible compared to travel times.
Table \ref{values_of_parameters} reports the parameter values adopted for the small-size instances.

\begin{table}[h!]
    \centering
    \resizebox{\textwidth}{!}{\begin{tabular}{l|c|c|c|c|c|c|c|c|c|c|c|c}
    \toprule
         Parameter & \makecell{$l_a$\\(units)} & \makecell{$l_b$\\(units)} & \makecell{$v_a$\\(km/h)} & \makecell{$v_b$\\(km/h)} & \makecell{$f_a$\\(\euro/d)} & \makecell{$f_b$\\(\euro/d)} & \makecell{$\omega_{\delta_a}$\\(\euro/km)} & \makecell{$\omega_{\delta_b}$\\(\euro/km)} & \makecell{$\Bar{t}$\\(hours)} & \makecell{$\bar{\delta}$\\(km)} & $\beta$ & \makecell{$\gamma$\\ (\euro)} \\
\midrule
         Value  &  15 &
         5 & 45 &
          15 & 7  &
           3  &
          0,20 & 
          0,15 & 5 & 2 &
          2 & 100\\
         \bottomrule
    \end{tabular}}
    \caption{Values of the parameters of the model for the tests on small size instances (adapted from \cite{carracedo2022electric} and \cite{rudolph2017cargo}).}
    \label{values_of_parameters}
\end{table}

\paragraph{Large Instances} To generate large instances, elementary demand units have been first created by aggregating customers, according to their geographical location. Thus, $67$ demand nodes have been considered, covering the 95\% of the total demand included in the real world dataset, used to generate large instances. To model the stochasticity in the daily demand, we exploit the information included in the dataset described in Section \ref{sec:data_and_context}. Therefore, each working day (Saturdays and public holidays excluded) is considered a different scenario for demand node realization. Hence, we derive 84 different scenarios for customers demand. As for the values of the parameters, they are the same as those included in Table \ref{values_of_parameters}. The only differences are in the values of $l_p$, that are $l_a=170$ units and $l_b=100$ units, and of $f_p$, that are $f_a=13$\euro/d and $f_a=6.5$\euro/d, reflecting realistic parameters. 
Finally, routes in $\mathcal{R}$ are grouped together according to the vehicle type $p$ potentially assigned to them, to define sets $R_a$ and $R_b$. Specifically, we suppose that type $a$ vehicles can be assigned to every $r\in\mathcal{R}$ (that is $\mathcal{R}_a=\mathcal{R}$). On the other hand, $R_b$ includes routes that are feasible for deliveries with ECBs, which are characterized by a driving range $\Delta_b=15$km.

\subsubsection{Path Generation}\label{sec:path_generations}

Candidate routes are generated using a parallelized Adaptive Large Neighborhood Search (ALNS) heuristic for the Capacitated Vehicle Routing Problem (CVRP), following \cite{pisinger2007general} and \cite{kocc2015hybrid}. The heuristic employs classical removal and insertion operators, including random removal, worst removal, Shaw removal based on spatial proximity, greedy insertion, and regret insertion.

To enhance diversification, $N=10$ initial solutions are constructed in parallel using a giant-tour representation followed by a standard splitting procedure. Each solution is then refined through the ALNS search process, and the best solution from each run is retained. The procedure is repeated under varying demand realizations and fleet configurations, using the same computational conditions as those adopted in the optimization phase.

The final set of candidate routes $\mathcal{R}$ is obtained by selecting the 900 most frequently activated routes across all runs. This threshold provides a good balance between solution diversity and model tractability as tests showed that increasing this number yields negligible improvement in solution quality while significantly increasing computational time.

\subsection{Computational Results}\label{sec:comput_results}

The computational results are summarized below. Section \ref{sec:small_results} reports the outcomes obtained on small instances and assesses the impact of demand stochasticity once the scenario reduction procedure has been applied. Section \ref{sec:large_results} presents the numerical results for large instances, demonstrating the effectiveness of the proposed formulations on the real-world case study.

All experiments were performed in Python 3.11.4 using Gurobi 11.0.3 as the optimization solver. Tests were executed on an Intel(R) Core(TM) i7-8565U 64-bit processor running at 1.80 GHz with 8 GB of RAM.

\subsubsection{Small Instances Results}\label{sec:small_results}

This section discusses the results obtained from numerical experiments on small instances. For these tests, the computational time limit was set to one hour for deterministic instances and one day for stochastic ones.

First, model $\mathcal{M}$ was tested on deterministic instances to evaluate the effect of valid inequalities (\ref{valid_in}) on computational performance. Table \ref{results_deterministic_model} reports the results obtained with and without valid inequalities (first and second block, respectively). The superscript $^\dagger$ on the CPU time indicates that the time limit was reached.

 \begin{table}[h!]
  \centering
  \resizebox{0.8\textwidth}{!}{
  \begin{tabular}{l|c|c|c|c|c|c|c}
  \toprule
    Instance & CMs & ECBs & \makecell{CMs distance\\ (km)} & \makecell{ECBs distance\\(km)}  & \makecell{Objective\\ Function (\euro/d)} & \makecell{CPU\\(s)} & Gap\\
\midrule
    $20\_1$ & 1     & 1     & 24,1720 & 14,3710 & 17,5612 & $3600,0^\dagger$ & 53,5\% \\
    $20\_2$ & 1     & 1     & 16,8351 & 11,6069 & 15,1081 & $3600,0^\dagger$  & 64,1\% \\
    $20\_3$ & 1 &	1	& 5,2565 &	0,9105	& 16,3926 &	$3600,0^\dagger$	&58,7\% \\
    $20\_4$ &1 &	1 &	32,8315 &	10,3896 &	18,3933 &	$3600,0^\dagger$ &	56,8\% \\
    $20\_5$ & 1 &	1	&28,4913	& 8,0007 &	16,8984 &	$3600,0^\dagger$	& 60,5\% \\
    \midrule
    20\_1 & 1     & 1     & 24,1720 & 18,1789 & 17,5612 & 1198,1 & 0\% \\
    20\_2 & 1     & 1     & 16,8351 & 11,6069 & 15,1081 & 44,8  & 0\% \\
    20\_3 & 1     & 1     & 17,1980 & 17,8762 & 16,1210 & 1016,1 & 0\% \\
    20\_4 & 1     & 1     & 19,3828 & 21,1815 & 17,0538 & 122,4 & 0\% \\
    20\_5 & 1     & 1     & 16,4940 & 15,8684 & 15,9540 & 147,5  & 0\% \\
    \bottomrule
    
    \end{tabular}}
    \caption{Computational results of the model $\mathcal{M}$ without valid inequalities (first block) and with valid inequalities, respectively (second block).}
  \label{results_deterministic_model}
\end{table}

According to these results, an optimal solution is found within the time limit only when valid inequalities are included in the formulation. Therefore, in the following analyses, we refer exclusively to model $\mathcal{M}$ with the valid inequalities (\ref{valid_in}) activated.
 
Subsequently, $\mathcal{M}$ was tested on small stochastic instances. These experiments were performed on stochastic versions generated from instance $20\_2$. Specifically, ten runs were executed for each case with an increasing number of scenarios, i.e., $|\mathcal{S}| \in  \{5, 10, 20, 50, 100, 150\}$. Table \ref{avg_stochastic_computational_results} reports the average results, while extended outcomes are provided in the \hyperref[sec:appendix_b]{Appendix B}. Besides fleet size and composition (columns CMs and ECBs) and traveled distances (columns CMs distances and ECBs distances), the table reports the number of complete recourse actions, i.e. additional demand nodes fully served by a single vehicle, and split recourse actions, i.e. demand nodes visited by multiple vehicles. The number of unserved customers (either fully or partially) is also provided to indicate outsourced deliveries. Finally, for each instance, we report the objective function value ($RP$), CPU time, and average optimality gap. 

\begin{table}[h!]
  \centering
  \resizebox{\textwidth}{!}{
    \begin{tabular}
    {l|c|c|c|c|c|c|c|c|c|c|c|c}
    \toprule
    Instance & $\mid \mathcal{S} \mid $ & CMs & ECBs  & \makecell{CMs distance\\ (km)} & \makecell{ECBs distance\\(km)} & \makecell{Complete\\ recourse \\actions} & \makecell{Split \\recourse \\actions} &  \makecell{Complete\\ unserved \\customers}& \makecell{Partial\\ unserved \\customers}& \makecell{$RP$\\ (\euro/d)} & \makecell{CPU\\(s)} & Gap\\
    \midrule
    \text{20\_2\_5\_AVG} & 5 & 3 & 1 & 31,7644 & 3,9527
 & 0 & 4 & 0 & 0  & 30,7698 & 1073,6 & 0\% \\
    \text{20\_2\_10\_AVG} & 10 & 3 & 1 & 32,8784 & 4,3069
 & 1 & 6 & 0 & 0 &  32,2151 & 2231,9 & 0\% \\
    \text{20\_2\_20\_AVG} & 20 & 3 & 1 & 33,3863 &	2,5766 & 3 & 18 &  0 & 0 & 33,0727 & 6943,6 & 0\% \\
    \text{20\_2\_50\_AVG} & 50 & 3 & 1 & 34,0111 &	2,6061
 & 6 & 36 & 0 & 2 &  33,4368 & 32443,2 & 0\% \\
    \text{20\_2\_100\_AVG} & 100 & 3& 1 & 34,6007 &	2,5415 & 26 & 84 & 0 & 4 & 33,4370 & 56934,9 & 0\%\\  
    \text{20\_2\_150\_AVG}  & 150 & 3 & 1 & 34,6346 & 2,8002 & 28 & 91 & 0 & 6 & 33,4372 & 76470,6 & 0,58\% \\
    \bottomrule
    \end{tabular}}
  \caption{Average computational results provided by $\mathcal{M}$ for stochastic instances with increasing number of scenarios.The notation $20\_2\_|\mathcal{S}|\_{\text{AVG}}$ refers to the average results of the instances with $|\mathcal{S}|$ scenarios generated from $20\_2$).}
\label{avg_stochastic_computational_results}
\end{table}

As shown in Table \ref{avg_stochastic_computational_results}, the average fleet size and composition remain constant over different cardinalities of the scenario trees. However, as indicated by the extended results in \hyperref[sec:appendix_b]{Appendix B}, these values depend on the realization of customer demand and become more stable as the cardinality of $\mathcal{S}$ increases. Average traveled distances also vary with $|\mathcal{S}|$, as different consistent routes are designed to serve more customers and reduce recourse actions, stabilizing when $|\mathcal{S}|\geq 100$. Regarding approximate recourse actions, their number generally increases as $|\mathcal{S}|$ grows, since a larger scenario set incorporates a broader spectrum of demand realizations. Nevertheless, after in-sample stability is achieved, the number of recourse actions tends to converge, suggesting that the scenario set is sufficiently representative of the underlying demand uncertainty.

\begin{figure}[h!]
    \centering
    \includegraphics[width=0.55\textwidth]{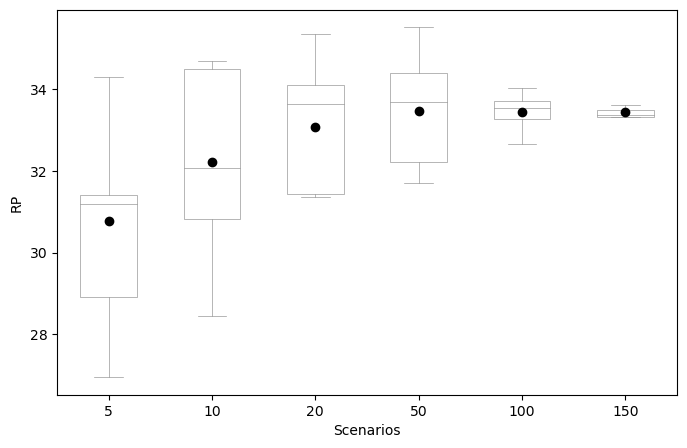}
    \caption{\small In-sample stability analysis for the stochastic instances generated from the 20\_2 deterministic instance. The box-plots represent the $RP$ value over 10 runs with increasing number of scenario.}
    \label{fig:in_sample_stability}
\end{figure}

From results in Table \ref{avg_stochastic_computational_results}, and box-plots of $RP$ values included in Figure \ref{fig:in_sample_stability} , we conclude that $RP$ reaches stability when $|\mathcal{S}| = 100$. 

Moreover, since the average CPU time increases with $|\mathcal{S}|$, and $\mathcal{M}$ fails to reach optimality within the time limit for some instances with $|\mathcal{S}|=150$, we set $|\mathcal{S}|=100$ as the benchmark for applying the scenario reduction procedure. The resulting reduced set of scenarios is then used to compute the classical stochastic measures for analyzing the impact of demand uncertainty.

The scenario reduction procedure is implemented by generating a scenario tree with $|\mathcal{S}| = 1000$ scenarios and retaining $|\mathcal{S}^{RED}| = 100$ of them. Table \ref{stochastic_computational_results_scen_red} reports the extended computational results obtained with $\mathcal{M}$ on the small stochastic instances after applying scenario reduction described in Section \ref{sec:scen_red}. The notation $20\_2\_100\_n\_RED$ denotes the instance resulting from the reduction procedure applied to the corresponding stochastic instance $20\_2\_100\_n$, for $n = 1,2,\dots,10$.

\begin{table}[h!]
  \centering
  \resizebox{\textwidth}{!}{
    \begin{tabular}
    {l|c|c|c|c|c|c|c|c|c|c|c}
    \toprule
    Instance & CMs & ECBs &  \makecell{CMs distance\\ (km)} & \makecell{ECBs distance\\(km)}  
 & \makecell{Complete\\ recourse \\actions} & \makecell{Split \\recourse \\actions} &  \makecell{Complete\\ unserved \\customers}& \makecell{Partial\\ unserved \\customers}&  \makecell{$RP$\\(\euro/d)} & \makecell{CPU\\(s)} & Gap\\
    \midrule
    20\_2\_100\_1\_RED  & 3     & 1     & 35,0554 & 1,5415 & 25 & 60 & 0 & 4 &  33,1316 & 43921,8 & 0\% \\
    20\_2\_100\_2\_RED  & 3     & 1     & 35,0353 & 1,5415 & 8 & 63 & 0 & 2 &  32,9840 & 58721,3 & 0\% \\
    20\_2\_100\_3\_RED  & 3     & 1     & 28,9750 & 9,3979 & 13 & 70 & 0 & 4 &  32,8903 & 49628,5 & 0\% \\
    20\_2\_100\_4\_RED  & 3     & 1     & 34,3780 & 3,5608 & 12 & 79 & 0 & 1 &  32,9113 & 44862,3 & 0\% \\
    20\_2\_100\_5\_RED  & 3     & 1     & 30,4693 & 10,0905 & 8 & 85 & 0 & 4 &  33,0521 & 43892,6 & 0\% \\
    20\_2\_100\_6\_RED  & 3     & 1     & 33,5473 & 3,5608 & 37 & 79 & 0 & 1 &  32,9034 & 55762,7 & 0\% \\
    20\_2\_100\_7\_RED  & 3     & 1     & 30,5905 & 10,0905 & 13 & 76 & 0 & 3 &  33,1873 & 62478,1 & 0\% \\
    20\_2\_100\_8\_RED  & 3     & 1     & 30,1152 & 9,3979 & 16 & 62 & 0 & 4 &  33,0014 & 38291,6 & 0\% \\
    20\_2\_100\_9\_RED  & 3     & 1     & 33,1886 & 3,1139 & 39 & 90 & 0 & 5 &  33,3484 & 68923,0 & 0\% \\
    20\_2\_100\_10\_RED & 3     & 1     & 35,6967 & 5,4020 & 19 & 72 & 0 & 1 &  32,7652 & 58273,7 & 0\% \\
    \midrule
     Average  & 3 & 1 & 32,6967 & 5,4020 & 19 & 74 & 0 & 3 & 33,0175 & 52475,6 & 0\%\\

    \bottomrule
    \end{tabular}}
    \caption{Computational results of $\mathcal{M}$ for reduced stochastic instances with increasing number of scenarios.}
    \label{stochastic_computational_results_scen_red}
    \end{table}

Results in Table \ref{stochastic_computational_results_scen_red} indicate that the fleet size and composition remain stable across instances. However, different routing plans are produced for the fleet, with ECBs performing, on average, longer routes. This leads to a higher number of split recourse actions compared with complete recourse actions.

The instances generated through scenario reduction are then used to evaluate the impact of customer demand uncertainty. First, we compare the stochastic formulation ($RP$) with the perfect-information case, i.e., the \textit{Wait and See} solution ($WS$), through the \textit{Expected Value of Perfect Information} ($EVPI$). We then quantify the benefits of explicitly accounting for uncertainty by comparing the stochastic model with its deterministic counterpart - the \textit{Expected Value} problem $(EV)$ - through the \textit{Value of the Stochastic Solution} ($VSS$) (see \cite{birge2011introduction}).  To compute $VSS$, we determine the \textit{Expectation of using the EV solution} ($EEV$). The $EEV$ is the objective value of the $RP$ model when first-stage decisions are fixed to those of the $EV$ solution. 
In addition, we compute the \textit{Loss of Upgrading the Deterministic Solution} ($LUDS$) (see \cite{maggioni2012analyzing}). With the $LUDS$, we seek to determine
whether the solution to the $EV$ is upgradeable, i.e., it can be used
as a starting point for generating a high-quality solution to the $RP$. To
do so, we solve the $RP$, albeit with additional constraints ensuring that
the values of the first-stage variables are at least as large as their values
in the optimal solution to the $EV$. We refer to the objective function
value of the optimal solution to this restricted $RP$ as the \emph{Expected Input
Value} ($EIV$) and compute the $LUDS$ as a matheuristic
in which a mathematical model, based on information from the deterministic
solution, is solved. 
 Since first-stage decisions include both fleet sizing and routing, we evaluate two approaches when computing $EEV$ and $EIV$. In the first, all first-stage decisions are fixed to the $EV$ solution, yielding the model $EEV_{FR}$ (Fleet and Routes fixed). In the second, only the fleet decisions are fixed, allowing the routing plan to adjust; we denote this by $EEV_{F}$ (Fleet fixed). Analogously, we refer to $EIV_{FR}$ and $EIV_{F}$ when applying the same logic to the computation of $EIV$. 
In Table \ref{tab:stoch_measures_1} and Table \ref{tab:stoch_measures_2} values of the classical stochastic measures are reported (see \hyperref[sec:appendix_c]{Appendix C} for extended results). 

\begin{table}[h!]
    \centering
    \resizebox{\textwidth}{!}{
    \begin{tabular}{l|c|c|c|c|c|c|c|c|c|c|c}
    \toprule
          Instance & $RP$ & $EV$ &$WS$ & $EVPI$ & $EEV_{FR}$ & $EEV_{F}$ &$VSS_{FR}$ & $VSS_{F}$ & $\%EVPI$ & $\%VSS_{FR}$ & $\%VSS_{F}$ \\
         \midrule
         20\_2\_100\_1\_RED & 33,1316 & 27,1944 & 26,7279 & 6,4037 & 60,3214 & 60,2146 & 27,1898 & 27,0830 & 19,33\% & 82,07\% & 81,74\%  \\
         20\_2\_100\_2\_RED & 32,9840 & 27,0981& 26,9011 & 6,0829 & 53,6357 & 53,6245 & 20,6517 & 20,6405&  18,44\% & 62,61\% & 62,58\%\\
         20\_2\_100\_3\_RED  & 32,8903  & 27,1070 & 26,8278 &  6,0625 & 61,7699 & 61,6935 & 28,8796 & 28,8032 & 18,43\% & 87,81\% & 87,57\%\\
         20\_2\_100\_4\_RED & 32,9113 & 27,9745
 & 27,3106 & 5,6007 & 36,2036 & 36,1716 & 3,2923
 & 3,2603 & 17,02\% & 10,00\% & 9,91\% \\
         20\_2\_100\_5\_RED & 33,0521 & 27,5756 & 26,8904
 & 6,1617 & 38,7154 & 38,6831 & 5,6633 & 5,6633 & 18,64\% & 17,13\% & 17,04\% \\
         20\_2\_100\_6\_RED & 32,9034 & 27,9745 & 26,8989 & 6,0045 & 37,0564 & 37,0564 & 4,1530 & 4,1530 & 18,25\% & 12,62\% &  12,62\%\\
         20\_2\_100\_7\_RED & 33,1873 & 27,6445 &  26,7044 & 6,4829 & 35,8342 & 35,8247 &  2,6469 & 2,6374  & 19,53\% & 7,98\% & 7,95\% \\
         20\_2\_100\_8\_RED & 33,0014 & 27,2800 & 26,8584&  6,1430 & 61,4908 & 61,4353 & 28,4894
          & 28,4339 & 18,61\% & 86,33\% & 86,16\%\\
         20\_2\_100\_9\_RED & 33,3484 & 27,8658 & 27,0660 &  6,2824 & 37,0132 & 37,0132 & 3,6648 & 3,6648 & 18,84\% & 10,99\% & 10,99\% \\
         20\_2\_100\_10\_RED & 32,7652 & 27,2421 & 26,8165 & 5,9487 & 61,6554 & 61,5263 & 28,8902 & 28,7611 & 18,16\% & 88,17\% & 87,78\%\\
         \hline
         Average & 33,0175 & 27,4956 & 26,9091 & 6,9091 & 48,5030 & 48,4555 & 15,5221 & 15,4746 & 18,41\% & 47,18\% & 47,04\% \\
         \bottomrule
    \end{tabular}}
    \caption{Summary results of stochastic measures $EVPI$, $VSS_{FR}$ and $VSS_{F}$, also expressed in 
    percentage gap to the corresponding $RP$ problem.}
    \label{tab:stoch_measures_1}
\end{table}

According to Table \ref{tab:stoch_measures_1}, the average value of $\%EVPI$ is 18,41\%, meaning that a decision maker would be ready to pay at most 18,41\% of the profit to get the perfect information about the customers demand. As for the $VSS$, average values of $\%VSS_{FR}$ and $\%VSS_{F}$ are 47,18\% and 47,04\%, respectively. This means that the cost for ignoring uncertainty is almost half of the $RP$. Additionally, since there is no significant difference between $\%VSS_{FR}$ and $\%VSS_{F}$, fleet decisions are more relevant than routing decisions in the stochastic framework. Even if only the fleet is fixed, indeed, it is still difficult to cover all the customers demand without an external delivery service. The results discussed so far justify the adoption of a stochastic model when dealing with a Last Mile delivery problem, since the deterministic solutions in a stochastic framework leads to significantly high additional costs. As for the upgradeability of the $EV$ solution, values of $LUDS_{FR}$ and $LUDS_{F}$ are reported in Table \ref{tab:stoch_measures_2}. First we note that $LUDS_{FR}=VSS_{FR}$ for all the instances, underlying the \textit{non upgradeability} of the deterministic solution if $EV$ fleet and routes are set as a lower bound in the $RP$ program. However, the $EV$ solution can be used as a valuable starting point in the stochastic framework if a lower bound is set only for fleet decisions. $EIV_{F}$, indeed, provides updated solutions characterized by a different fleet size and composition, allowing to reduce the number of recourse actions and, consequently, the total costs. Furthermore, from the extended results in \hyperref[sec:appendix_c]{Appendix C}, we observe that there is \textit{perfect upgradeability} when the $EV$ solution provides a delivery fleet smaller than the one provided by the $RP$. On the other hand, when the $EV$  delivery fleet has the same size of the $RP$ one, but with a different composition, there is only \textit{partial upgradeability}. This happens since the $EIV$ optimal solution involves a higher number of vehicles in the fleet to reduce the number of unserved customers.

\begin{table}[h!]
    \centering
    \resizebox{0.9\textwidth}{!}{
    \begin{tabular}{l|c|c|c|c|c|c|c}
    \toprule
          Instance & $RP$ &  $EIV_{FR}$ & $EIV_{F}$ & $LUDS_{FR}$ & $LUDS_{F}$ & $\%LUDS_{FR}$ & $\%LUDS_{F}$ \\
         \midrule
         20\_2\_100\_1\_RED & 33,1316 & 60,3214 & 34,7305 & 27,1898 & 1,5989 & 82,07\% & 4,83\%  \\
         20\_2\_100\_2\_RED & 32,9840 & 53,6357 & 34,6758 & 20,6517 & 1,6918 & 62,61\% & 5,13\%\\
         20\_2\_100\_3\_RED  & 32,8903  & 61,7699 &  34,4339 & 28,8796 & 1,5436 & 87,81\% & 4,69\% \\
         20\_2\_100\_4\_RED & 32,9113
 & 36,2036 & 32,9113 
 & 3,2923 & 0,0000 
& 10,00\% & 0,00\%\\
         20\_2\_100\_5\_RED & 33,0521 & 38,7154
 & 33,0521  & 5,6310
 & 0,0000 & 17,13\% & 0,00\%\\
         20\_2\_100\_6\_RED & 32,9034 & 37,0564 & 32,9304 & 4,1530 & 0,0000 & 12,62\% & 0,00\%\\
         20\_2\_100\_7\_RED &33,1873 & 35,8342 & 33,1873 & 2,6469
 & 0,0000 & 7,98\% & 0,00\% \\
         20\_2\_100\_8\_RED & 33,0014 & 61,4908&  34,5814 & 28,2894 & 1,5800 & 86,33\% & 4,79\%\\
         20\_2\_100\_9\_RED & 33,3484 & 37,0132 & 33,3484 & 3,6648 & 0,0000 & 10,99\% & 0,00\% \\
         20\_2\_100\_10\_RED & 32,7652 & 61,6554 & 34,0193 & 28,8902 & 1,2541 & 88,17\% & 3,83\%\\
         \midrule
         Average & 33,0175 & 48,5030 & 33,7843 & 15,5221 & 0,7324 & 47,18\% & 2,23\%\\
         \bottomrule
    \end{tabular}}
    \caption{Summary results of $LUDS_{FR}$ and $LUDS_F$, also expressed in 
    percentage gap to the corresponding $RP$ problem.}
    \label{tab:stoch_measures_2}
\end{table}

As reported in Table \ref{tab:stoch_measures_2}, $\%LUDS_{F}$ is, on average, 2,23\% of the $RP$ corresponding value. 

Therefore, an important takeaway for decision makers is that the $EV$ solution can serve as a valuable starting point within the stochastic framework, as it can be upgraded without leading to substantially higher delivery costs. As discussed above, the $LUDS$ provides a means to assess whether the $EV$ solution offers useful guidance in the stochastic setting. More refined stochastic measures, which further evaluate how the structure of the $EV$ optimal solution can be exploited, can be found in \cite{crainic2018reduced}.

\subsubsection{Real-World Case Study Results}\label{sec:large_results}

In this section, we discuss the computational results obtained on large real-world instances. 

We first evaluate the performance of $\mathcal{M}$ and $\mathcal{M}_{PB}$ on the same deterministic instances. Five deterministic instances, each representing a different working day, are selected based on the customer demand distribution. All instances include the same number of demand nodes, $|\mathcal{I}| = 67$, but differ in the total demand volume. Table \ref{results_det_model_real_data} summarizes the results of both $\mathcal{M}$ and $\mathcal{M}_{PB}$ on these deterministic instances. Results show that $\mathcal{M}_{PB}$ significantly outperforms $\mathcal{M}$, achieving notably smaller optimality gaps within the time limit of two hours. In addition, objective function values provided by $\mathcal{M}_{PB}$ is smaller in all the five instances, due to the fact that within the time limit the path-based formulation is able to select less expensive routes to be assigned to the vehicles included in the delivery fleet. Specifically, these routes allow CMs to cover longer distances and serve more customers, hence reducing the number of demand externalization.

\begin{table}[h!]
  \centering
  \resizebox{\textwidth}{!}{
  \begin{tabular}{l|c|c|c|c|c|c|c|c|c}
  \toprule
    Instance & $s\in\mathcal{S}$ & $\sum_{i\in\mathcal{I}} d_{is}$ & CMs & ECBs & \makecell{CMs distance\\ (km)} & \makecell{ECBs distance\\(km)} & \makecell{Objective\\ Function\\ (\euro/d)} & \makecell{CPU\\(s)} & Gap \\
    \midrule
    67\_1\_$\mathcal{M}$ & 71 & 979 & 3   & 5     & 54,5385 & 53,3917 & 91,1166 & 7200,0 & 15,38\% \\
    67\_2\_$\mathcal{M}$ & 83 & 1060 & 3  & 6   & 59,7699
 & 49,7761 & 99,3337
 & 7200,0  & 13,20\% \\
    67\_3\_$\mathcal{M}$ & 5 & 1154 & 4 & 5 & 69,5475 & 40,6694 & 105,2104 & 7200,0 & 13,01\% \\
    67\_4\_$\mathcal{M}$ & 21 & 1221 & 2  & 9  & 34,1926		
 & 92,0676 & 106,5276 & 7200,0 & 8,03\% \\
    67\_5\_$\mathcal{M}$ & 69 & 1281 & 3     & 8   		
  & 52,4664 & 73,7491 & 113,2559 & 7200,0  & 12,71\% \\
  \midrule
$67\_1\_\mathcal{M}_{PB}$ & 71 & 979 & 3   & 5 & 61,3284  & 39,9072  & 90,1865 & 7200,00 & 2,50\% \\
$67\_2\_\mathcal{M}_{PB}$ & 83 & 1060 & 3   & 6 & 59,1669  & 45,5170  & 96,7485 & 7200,00 & 3,26\% \\
$67\_3\_\mathcal{M}_{PB}$ & 5 & 1154 & 4   & 5 & 65,9975  & 40,2662  & 104,4049 & 7200,00 & 3,92\% \\
$67\_4\_\mathcal{M}_{PB}$ & 21 & 1221 & 2   & 9 & 39,6204  & 81,4506  & 105,3670 & 1915,76 & 0,00\% \\
$67\_5\_\mathcal{M}_{PB}$ & 69 & 1281 & 3   & 8 & 54,6874  & 65,3792  & 112,7720 & 7200,00 & 1,07\% \\
    \midrule
    \end{tabular}}
    \caption{Computational results of $\mathcal{M}$ and $\mathcal{M}_{PB}$ for the deterministic instances generated from the real-world data. The notation $67\_n\_\mathcal{M}$ and $67\_n\_\mathcal{M}_{PB}$ indicate the $n$-th instance on which model $\mathcal{M}$ or $\mathcal{M}_{PB}$ is tested, with $n=1,2,\dots,5$.}
\label{results_det_model_real_data}
\end{table}

Since solving the deterministic formulation using $\mathcal{M}$ or $\mathcal{M}_{PB}$ cannot be solved within a time limit of 2 hours, the full stochastic formulation of $\mathcal{M}_{PB}$ is solved through the adoption of the KS-based approach, described in Section \ref{sec:KS_approach}, enhancing the scalability for large real-world instances. However, since both $\mathcal{M}_{PB}$ and the KS-based approach are heuristic methods, in the next paragraph we will assess the quality of these approaches compared with $\mathcal{M}$.

\paragraph{Assessing the Quality of Path-Based and Kernel Search-based Approaches}\mbox{}\\

\noindent To quantify the near-optimality of the heuristic solutions to the optimal one, we compute the relative optimality gap as follows:
\begin{equation*}
\text{Gap}_{PB-\mathcal{M}}(\%) = 
\frac{M_{PB}^* - M_{\mathcal{M}}^*}{M_{\mathcal{M}}^*} \cdot 100,\qquad \text{Gap}_{\text{KS}(N)-\mathcal{M}}(\%) = 
\frac{M_{\text{KS}(N)}^* - M_{\mathcal{M}}^*}{M_{\mathcal{M}}^*} \cdot 100,
\end{equation*}
where $M_{PB}^*$ denotes the objective value obtained by approach $\mathcal{M}_{PB}$, $M_{\text{KS}(N)}^*$ denotes the objective value obtained by approach $\mathcal{M}_{\text{KS}(N)}$ with bucket size $N$, and $M_{\mathcal{M}}^*$ is the optimal objective value returned by model $\mathcal{M}$. Since our optimization problem is a minimization one, a positive gap indicates that the heuristic solution is worse than the optimal one, while a value equal to zero implies optimality.
In addition, to directly compare the performance of the two heuristic approaches, we compute the relative gap between their objective values:
\begin{equation*}
\text{Gap}_{\text{KS}(N)-PB}(\%) =
\frac{M_{\text{KS}(N)}^*-M^*_{\mathcal{M}_{PB}} }{M^*_{\mathcal{M}_{PB}}} \cdot 100.
\end{equation*}
This measure evaluates the relative difference between the solutions obtained by the two heuristic methods. In the context of a minimization problem, a negative value indicates that the solution provided by $\mathcal{M}_{PB}$ is worse than that obtained by the KS-based approach, whereas a positive value indicates the opposite.

To evaluate the near-optimality of the heuristic solutions we consider a reduced real-world case study instance with $|\mathcal{I}| = 15$ customers and the full set of scenarios ($|\mathcal{S}| = 84$). Under this reduced setting, model $\mathcal{M}$ is computationally tractable and can determine an optimal solution, which serves as a benchmark for evaluating the performance of the heuristic approaches. Given the reduced number of demand nodes, a new tailored set of preprocessed routes $\mathcal{R}$ was generated, comprising approximately 320 routes. In addition, we consider three different values of the penalty parameter $\beta$, namely $\beta = \{1, 2, 3\}$, while keeping all other problem parameters unchanged. For the KS-based method, three bucket sizes are tested, i.e. $N \in \{20, 50, 100\}$, and \textit{OT} is set to $0.5\%$. The time limit for these experiments is set to 12 hours (i.e. $T_{max}=43600 s$). Table~\ref{tab:perform_comp} reports the computational results for this reduced real-world case study instance. Results show that both heuristic approaches provide solutions within a limited near-optimality range with respect to the optimal benchmark obtained by $\mathcal{M}$. In particular, the optimality gaps remain relatively small across all tested values of the penalty parameter $\beta$, confirming that the implemented heuristics are able to identify high-quality solutions. Even in the most challenging configurations, the deviation from optimality is contained within a few percentage points, which supports the practical applicability of these approaches in the large-scale setting where solving $\mathcal{M}$ to optimality becomes computationally intractable.

\begin{table}[h!]
	\centering
	\resizebox{0.6\textwidth}{!}{
		\begin{tabular}{c|c|c|c|c|c}
			\toprule
			Approach & $\beta$ & \makecell{Objective\\Function\\(\euro/d)} &  \makecell{CPU\\(s)} & $\text{Gap}_{(\cdot) -\mathcal{M}}(\%)$ &$\text{Gap}_{\text{KS}(N)-PB}(\%) $\\
			\midrule
			$\mathcal{M}$ & & 48.8928 & 436.57 & - & -\\
			$\mathcal{M}_{PB}$ &   & 50.6984 & 6110.33 & 3.69 & - \\
			$\mathcal{M}_{\text{KS}(20)}$& 1 & 49.9542 & 1931.8 & 2.34 & 1.47\\
			$\mathcal{M}_{\text{KS}(50)}$ & & 49.9912 & 6312.6 & 2.27 & 1.39\\
			$\mathcal{M}_{\text{KS}(100)}$ &  & 49.9541 & 5890.2 & 2.34 & 1.47\\
			\midrule
			$\mathcal{M}$ & & 51.1508 & 492.27	& - & -\\
			$\mathcal{M}_{PB}$ &  & 52.0525 & 5109.49	 & 1.76 & - \\
			$\mathcal{M}_{\text{KS}(20)}$ & 2 & 52.1364 & 3502.40 & 1.93 & \textbf{-0.16} \\
			$\mathcal{M}_{\text{KS}(50)}$ &  & 52.1359 & 2575.80 & 1.93 & \textbf{-0.16}\\
			$\mathcal{M}_{\text{KS}(100)}$ &  &52.0831 & 4386.61 & 1.82 & \textbf{-0.06}\\
			\midrule
			$\mathcal{M}$ & &53.0742 & 634.48 & - & -\\
			$\mathcal{M}_{PB}$ &  & 53.1210 & 5109.49 & 0.09 & - \\
			$\mathcal{M}_{\text{KS}(20)}$ & 3 & 53.1056 & 2122.43 & 0.06 & 0.03 \\
			$\mathcal{M}_{\text{KS}(50)}$ &  & 53.1667 & 4447.81 & 0.17 & \textbf{-0.09} \\
			$\mathcal{M}_{\text{KS}(100)}$ & & 53.2466 & 4560.27 & 0.32 & \textbf{-0.24}\\
			\bottomrule
	\end{tabular}}
	\caption{Computational results provided by $\mathcal{M}$, $\mathcal{M}_{PB}$ and the KS$(N)$-based approach on reduced real-world case study instances with a time limit of 12 hours. $\text{Gap}_{(\cdot) -\mathcal{M}}(\%)$ denotes the percentage gap between $PB$ or KS$(N)$ and model $\mathcal{M}$.}
	\label{tab:perform_comp}
\end{table}

A comparison between the two heuristic methods further highlights that the KS-based approach achieves solution values that are very close to those obtained by $\mathcal{M}_{PB}$, and in some instances even slightly improves them, as indicated by the small (and occasionally negative) relative gaps between the two heuristics. The best performances are reported in bold. At the same time, the KS-based method exhibits a more structured and scalable solution process, particularly when varying the bucket size parameter. These results suggest that the KS procedure represents a valid and effective strategy to enhance the scalability of the problem while preserving a satisfactory level of solution quality, thus making it particularly suitable for larger real-world instances.

Finally, we test the effectiveness of the KS-based approach against $\mathcal{M}$ and $\mathcal{M}_{PB}$ on the full real-world case study instances, with increasing values of $|\mathcal{S}|$, namely $|\mathcal{S}| \in  \{10, 30, 60, 84\}$. Bucket sizes, \textit{OT}, and time limit are set as before. The corresponding results are reported in Table \ref{tab:KS_results}.

\begin{table}[h!]
    \centering
    \resizebox{\textwidth}{!}{
    \begin{tabular}{c|c|c|c|c|c|c|c|c|c|c|c|c|c|c|c}
        \toprule
        Instance & $|\mathcal{S}|$ & \multicolumn{4}{c|}{$\mathcal{M}$} & \multicolumn{4}{c|}{$\mathcal{M}_{PB}$} & \multicolumn{3}{c|}{$\text{KS}(N)$} & \makecell{$\text{Gap}_{PB-\mathcal{M}}$\\$(\%)$} & \makecell{$\text{Gap}_{\text{KS}(N)-\mathcal{M}}$\\$(\%)$}  &\makecell{$\text{Gap}_{\text{KS}(N)-PB}$\\$(\%)$}\\
        \cline{3-13}
        & & \makecell{CPU\\(h)} & \makecell{Incumbent\\(\euro/d)} & \makecell{Lower Bound\\(\euro/d)} & Gap & \makecell{CPU\\(h)} & \makecell{Incumbent\\(\euro/d)} & \makecell{Lower Bound\\(\euro/d)} & Gap & \makecell{CPU\\(h)} & \makecell{Objective \\Function\\ (\euro/d)}  & $N$& & & \\
        \midrule
        1 & 10 & 12,00 & 114,7662 & 105,1791 & 8,35\% & 12,00 & 115,4832 & 114,6547 & 0,72\% & 4,63 & 116,9417& 20 & & 1,86\% & 1,25\% \\
        & & & & & & & & & & 2,50 & 116,2527 & 50 & 0,62\% & 1,28\% & 0,66\% \\
        & & & & & & & & & &  2,27 & 116,0443  & 100 & & 1,10\% & 0,48\% \\
        \midrule
        2 & 30 & 12,00 & 119,1763 & 106,5978 & 10,55\% & 12,00 & 120,8643 & 116,9983 & 3,20\% & 5,43 & 120,1705 & 20 & & 0,83\% & \textbf{-0,58\%} \\
        & & & & & & & & & &  5,51 & 120,0772 & 50 & 1,40\% & 0,75\% & \textbf{-0,66\%} \\
        & & & & & & & & & &  9,61 & 119,7524 & 100 & & 0,48\% & \textbf{-0,93\%} \\
        \midrule
        3 & 60 & 12,00 & 120,7208 & 107,6772 & 10,80\% & 12,00 & 124,0341 & 116,7434 & 5,88\% & 11,52 & 121,5753 & 20 & & 0,70\% & \textbf{-2,02\%} \\
        & & & & & & & & & &  12,00  & 121,6983  & 50 & 2,67\% & 0,80\% & \textbf{-1,92\%} \\
        & & & & & & & & & &  12,00  & 122,0204  & 100 & & 1,07\% & \textbf{-1,65\%} \\
        \midrule
        4 & 84 & 12,00 & 129,3415 & 112,6587 & 12,90\% & 12,00 & 126,1372 & 118,3908 & 6,14\% & 12,00 & 124,0163  & 20 & & \textbf{-4,29\%} & \textbf{-1,71\%} \\
        & & & & & & & & & &  12,00 & 126,7430 & 50 & \textbf{-2,54\%} & \textbf{-2,05\%} & 0,48\% \\
        & & & & & & & & & &  12,00 & 124,9363 & 100 & & \textbf{-3,53\%} & \textbf{-0,96\%} \\
        \bottomrule
    \end{tabular}}
    \caption{Computational results provided by $\mathcal{M}$, $\mathcal{M}_{PB}$ and the $KS$-based approach on the full real-world case study instances.}
    \label{tab:KS_results}
\end{table}

\begin{figure}[h!]
\centering
\includegraphics[width=0.5\textwidth]{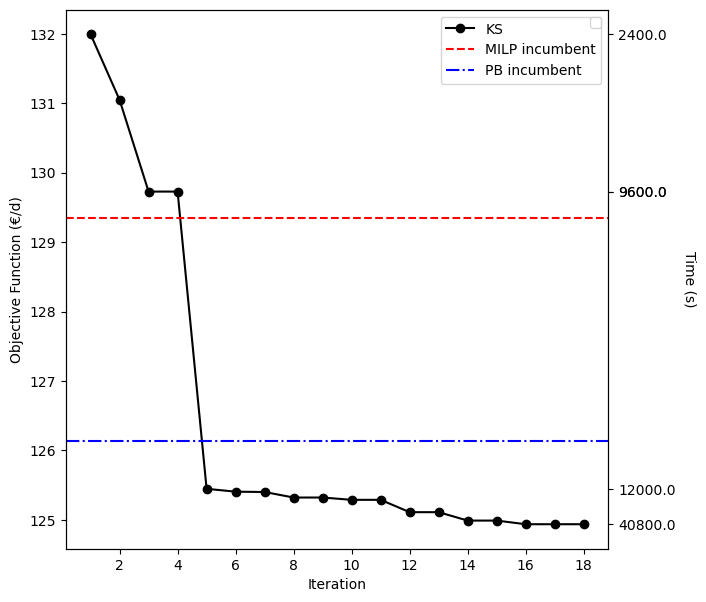}
\caption{Performance of KS-based approach, with bucket size $N=100$, compared to performances of $\mathcal{M}$ and $\mathcal{M}_{PB}$ when $|\mathcal{S}|=84$.}
\label{fig:KS_performances} 
\end{figure}

According to Table \ref{tab:KS_results}, neither $\mathcal{M}$ nor $\mathcal{M}_{PB}$ finds an optimal solution within the time limit, with optimality gaps increasing as the number of scenarios grows. In particular, $\mathcal{M}$ exhibits larger gaps than $\mathcal{M}_{PB}$, while $\mathcal{M}_{PB}$ yields a better incumbent for $|\mathcal{S}| = 84$, confirming the advantages of its reformulation. In contrast, the KS-based approach produces better solutions than $\mathcal{M}_{PB}$ for $|\mathcal{S}| > 10$, and also outperforms $\mathcal{M}$ for $|\mathcal{S}| = 84$. More importantly, the KS-based approach with $N = 100$ generates high-quality incumbent solutions than both $\mathcal{M}$ and $\mathcal{M}_{PB}$ after only five iterations, as shown in Figure \ref{fig:KS_performances}. It identifies a better incumbent 

In terms of CPU time, the incumbents obtained by $\mathcal{M}$ and $\mathcal{M}_{PB}$ require the full 12-hour limit, with optimality gaps of $12.9\%$ and $6.14\%$, respectively. In contrast, the KS-based approach identifies a better incumbent in approximately 3.5 hours, outperforming both models. Therefore, we derive  derive managerial insights from the implementation of the KS-based approach with $N = 100$.

\subsection{Practical Insights}\label{sec:practical_insights}

\begin{figure}[h!]
    \centering
    \subfloat[$WS$ fleets compared to the stochastic fleet provided by $\mathcal{M}_{PB}$ when $|\mathcal{S}|=84$.\label{fig:PB_fleet_size}]{\includegraphics[width=0.415\textwidth]{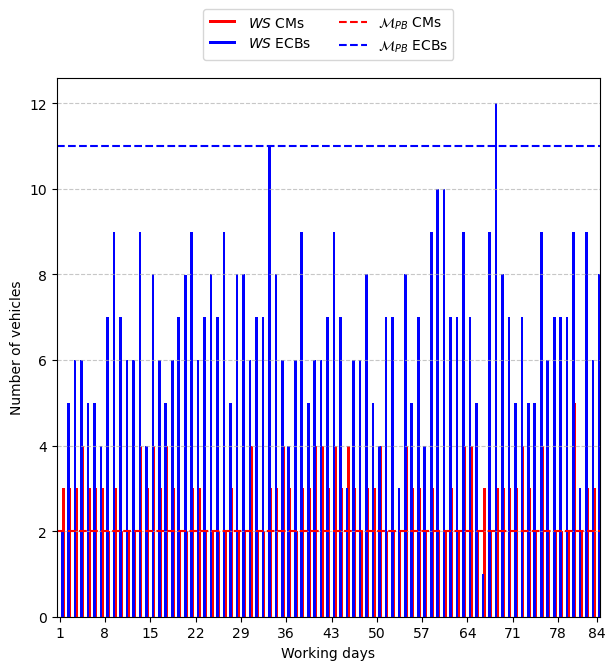}}\quad
    \subfloat[$WS$ traveled distances compared to the stochastic traveled distance provided by $\mathcal{M}_{PB}$ when $|\mathcal{S}|=84$.\label{fig:PB_covered_distances_beta_2}]{\includegraphics[width=0.45\textwidth]{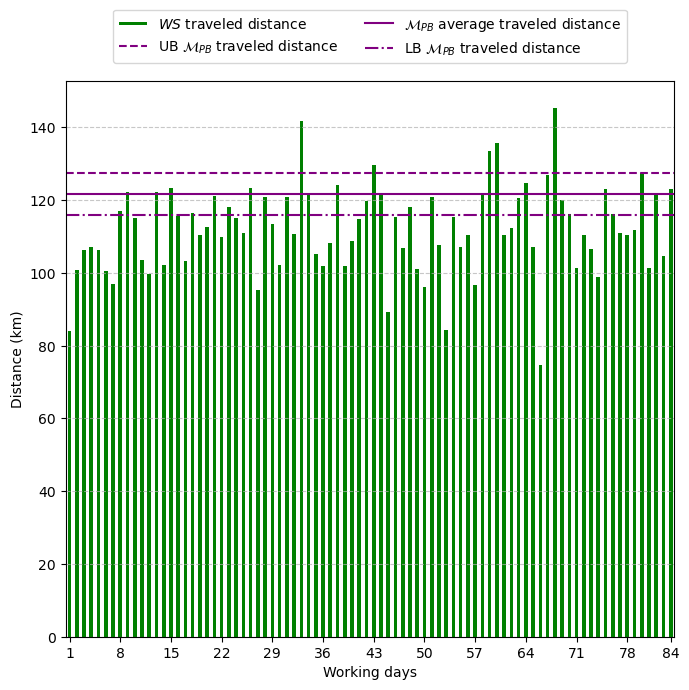}}
    \caption{Comparison between fleet composition and traveled distances of the scenario-specific deterministic $WS$ solution and stochastic $\mathcal{M}_{PB}$ solution.}
    \label{fig:comparison_stoch_det_real}
\end{figure}

In this section, we discuss several managerial insights derived from the solution of the large real-world instance obtained with the $\mathcal{M}_{PB}$ model solved with the KS-based heuristic. We first compare the scenario-specific deterministic solutions (i.e., the \textit{WS} solutions), obtained by solving $\mathcal{M}_{PB}$ for each working day considered as an individual scenario, with the stochastic solution, in which the entire set of working days is represented by the scenario set $\mathcal{S}$. Figure~\ref{fig:comparison_stoch_det_real} compares the fleet size, fleet composition (Figure~\ref{fig:PB_fleet_size}), and travel distances (Figure~\ref{fig:PB_covered_distances_beta_2}) obtained from the scenario-specific deterministic solutions with those resulting from the stochastic solution.

As depicted in Figure \ref{fig:PB_fleet_size}, while the \textit{WS} solutions determine the fleet composition separately for each working day, leading to substantial day-to-day variations driven by the realized customer demand, the stochastic $\mathcal{M}_{PB}$ solution identifies a single fleet configuration that remains unchanged across all scenarios. In the considered case study, the stochastic fleet consists of 2 CMs and 11 ECBs. The vertical bars represent the daily fleet composition obtained from the \textit{WS} solutions, whereas the dashed horizontal lines indicate the fixed fleet size selected by the stochastic model. Compared with the deterministic solutions, the larger number of ECBs in the stochastic fleet provides sufficient flexibility to absorb demand uncertainty, thereby reducing fleet variability and limiting the need for external delivery services.

Figure~\ref{fig:PB_covered_distances_beta_2} further illustrates the reduced variability achieved by the stochastic $\mathcal{M}_{PB}$ solution compared with \textit{WS} solutions. Each vertical green bar represents the total distance traveled on a given working day by the corresponding \textit{WS} solution, highlighting the significant day-to-day fluctuations induced by demand variability. In contrast, the stochastic solution yields much more stable routing decisions. The solid horizontal line denotes the average total traveled distance associated with the stochastic $\mathcal{M}_{PB}$ solution. Since the actual customer demand is revealed only during operations, the realized travel distance depends on the recourse actions activated after the first-stage routing decisions. Consequently, the total traveled distance varies only within the interval delimited by the Lower Bound (LB) and Upper Bound (UB), represented by the dashed horizontal lines. These bounds capture the range of travel distances resulting from the combination of the fixed first-stage routes and the scenario-dependent recourse actions. Notably, this interval is considerably narrower than the variability exhibited by the \textit{WS} solutions across working days, demonstrating that the stochastic approach effectively mitigates the impact of demand uncertainty.

\begin{figure}[h!]
    \centering
    \subfloat[Consistent routing with $\beta=1$.\label{fig:KS_consistent_route_1}]{\includegraphics[width=0.33\textwidth]{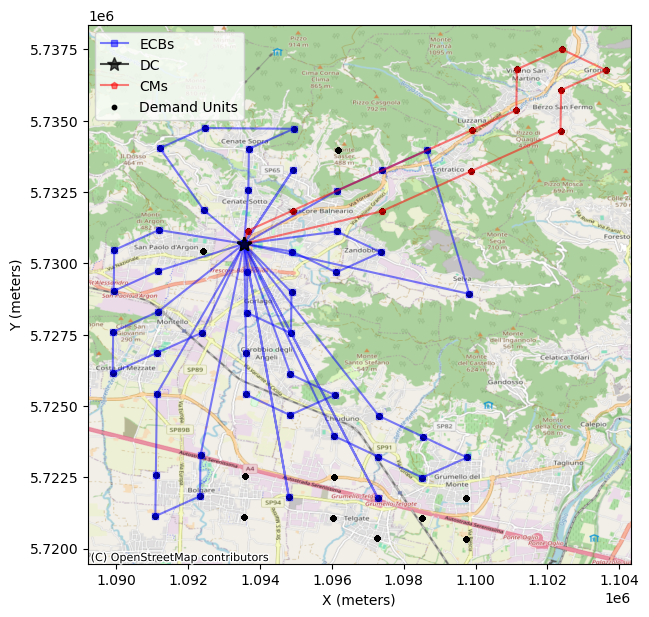}}
    \subfloat[Consistent routing with $\beta=2$.\label{fig:KS_consistent_route_2}]{\includegraphics[width=0.33\textwidth]{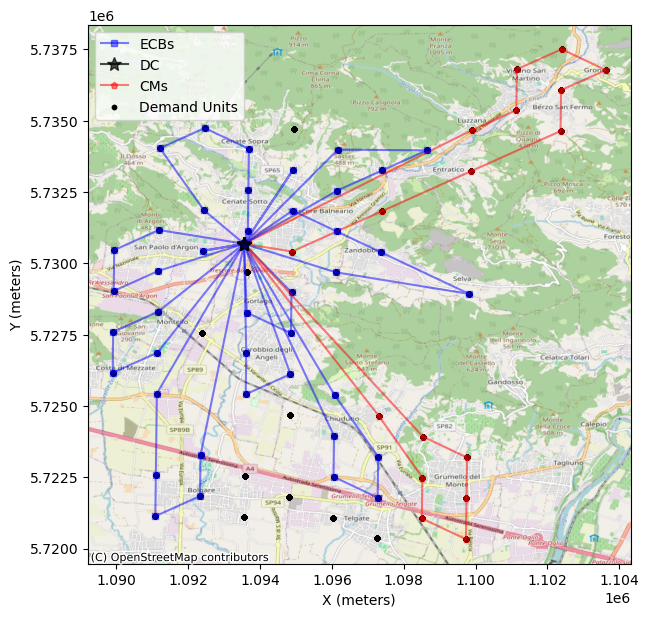}}
	\subfloat[Consistent routing with $\beta=3$.\label{fig:KS_consistent_route_3}]{\includegraphics[width=0.33\textwidth]{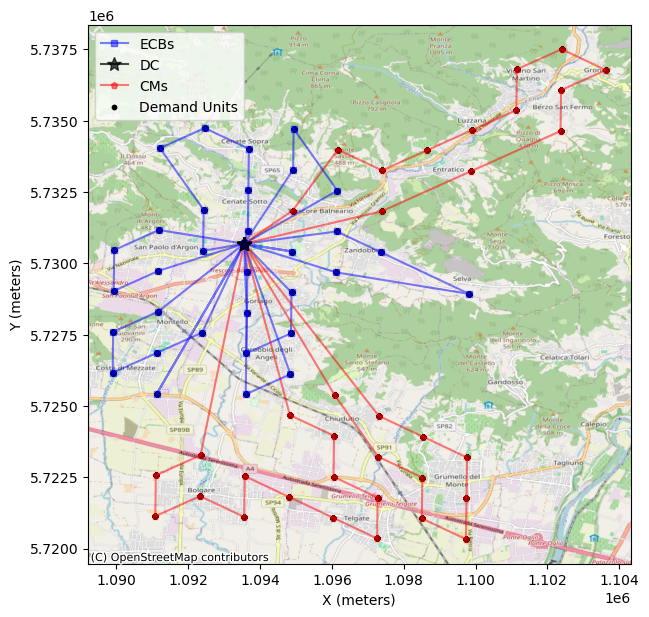}}\\
    \vspace{0.25cm}
    \subfloat[Recourse actions with $\beta=1$.\label{fig:KS_recourse_actions_1}]{\includegraphics[width=0.33\textwidth]{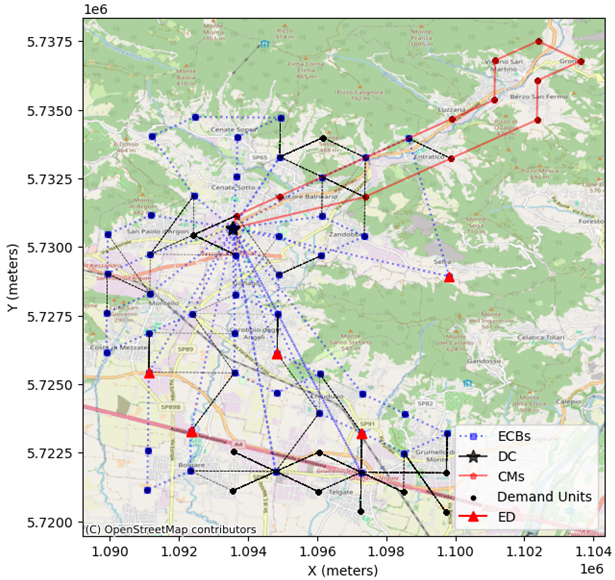}}
    \subfloat[Recourse actions with $\beta=2$.\label{fig:KS_recourse_actions_2}]{\includegraphics[width=0.33\textwidth]{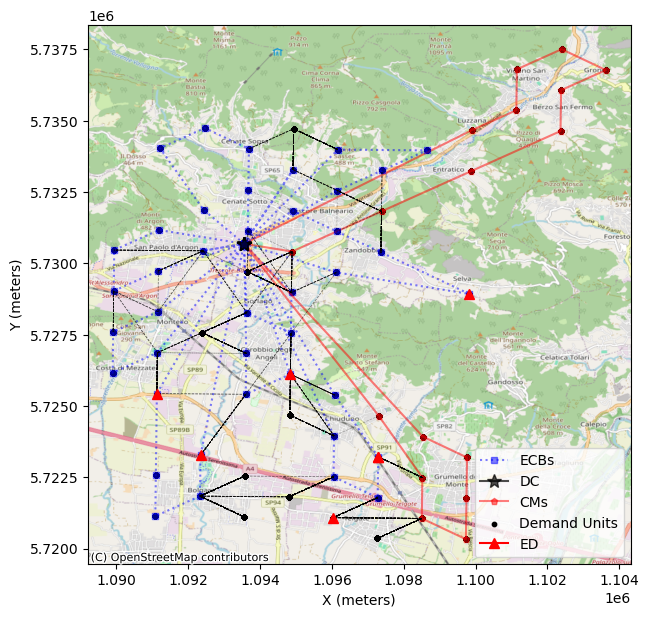}}
    \subfloat[Recourse actions with $\beta=3$.\label{fig:KS_recourse_actions_3}]{\includegraphics[width=0.33\textwidth]{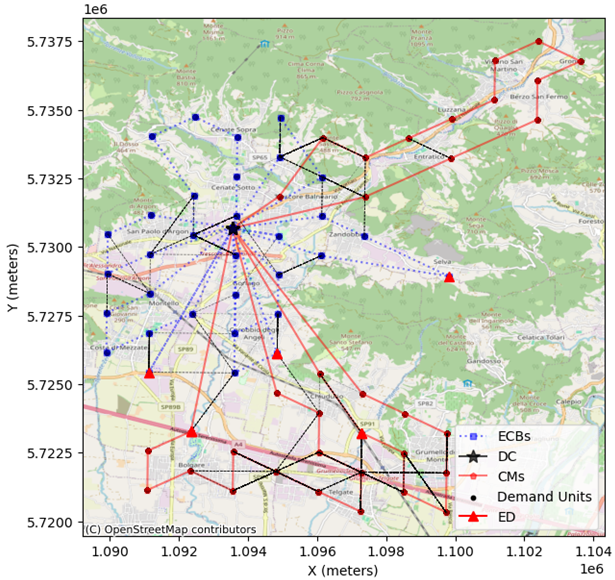}}    
    \caption{Routing plans for the real-world case study obtained with $\mathcal{M}_{PB}$ solved using the KS-based heuristic ($N=100$) for different values of $\beta=1,2,3$. Black dots represent customer demand locations, while red triangles identify customers assigned to Externalized Demand (ED), i.e., customers not included in the consistent first-stage routing plan.}\label{fig:routing_comparisons}
\end{figure}

The capability of the stochastic model $\mathcal{M}_{PB}$ to provide consistent first-stage delivery plans, with a flexibility to adapt according to demand realization, is further illustrated in Figure~\ref{fig:routing_comparisons}. 

Specifically, Figures~\ref{fig:KS_consistent_route_1}-\ref{fig:KS_consistent_route_3} show the consistent routing plans determined by $\mathcal{M}_{PB}$, whereas Figures~\ref{fig:KS_recourse_actions_1}-\ref{fig:KS_recourse_actions_3} illustrate how these plans are adjusted through approximate recourse actions according to the demand realized on each working day for different values of penalty parameter $\beta$.

As shown in Figures~\ref{fig:KS_consistent_route_1}-\ref{fig:KS_consistent_route_3}, increasing $\beta$ makes recourse actions more expensive and therefore encourages the model to construct more conservative delivery plans. When $\beta=1$ (Figure~\ref{fig:KS_consistent_route_1}), recourse actions are only mildly penalized, and the model deliberately leaves some customers (shown as red triangles) outside the consistent routing plan, expecting them to be served through recourse. As $\beta$ increases, relying on recourse becomes progressively less attractive. Consequently, for $\beta=3$ (Figure~\ref{fig:KS_consistent_route_3}), all customers are included in the consistent routing plan, significantly reducing the need for corrective actions during operations.

This trade-off also affects fleet composition. For $\beta=1$, the model selects a fleet consisting of 1 CM and 12 ECBs, as the additional capacity provided by a second CM is not justified by the relatively low recourse costs. As $\beta$ increases, the fleet becomes progressively more robust, evolving to 2 CMs and 11 ECBs for $\beta=2$, and 3 CMs and 9 ECBs for $\beta=3$. These larger and more diversified fleets increase the capability of absorbing demand uncertainty within the first-stage plan, thereby reducing the expected reliance on costly recourse actions.
 
The penalty parameter $\beta$ also influences the extent to which the consistent routing plans are adjusted through recourse actions. Figures~\ref{fig:KS_recourse_actions_1}, \ref{fig:KS_recourse_actions_2}, and \ref{fig:KS_recourse_actions_3} illustrate the approximate recourse actions (represented by black edges) activated for each scenario $s\in\mathcal{S}$ to satisfy additional scenario-specific demand. The darker the black edges, the more frequently the corresponding recourse actions are required across the scenarios. When $\beta=1$ (Figure~\ref{fig:KS_recourse_actions_1}), recourse actions are used more extensively because their contribution to the objective function is relatively small. As $\beta$ increases to 2 and 3 (Figures~\ref{fig:KS_recourse_actions_2} and \ref{fig:KS_recourse_actions_3}), recourse becomes progressively more expensive, encouraging the model to anticipate demand uncertainty within the first-stage routing plan and consequently reducing both the number and the frequency of recourse interventions.

Overall, these results demonstrate that the stochastic $\mathcal{M}_{PB}$ model successfully combines consistency with operational flexibility. It generates a common set of delivery routes that can be repeatedly implemented over the planning horizon while requiring only limited adjustments after demand realization. The parameter $\beta$ provides a direct mechanism for controlling the trade-off between first-stage consistency and operational flexibility, allowing decision makers to calibrate the desired level of conservatism. From a practical perspective, such consistent routing plans simplify daily operations, facilitate driver familiarization with delivery routes, improve planning predictability, and reduce the operational disruptions caused by day-to-day demand fluctuations.

\section{Conclusions}\label{sec:conclusion}

In this paper, we studied a Stochastic Fleet Size ands Mix Consistent Vehicle Routing Problem for a Last Mile delivery context with uncertain demand. To address this challenge, we proposed a two-stage stochastic programming formulation that minimizes delivery costs arising from tactical decisions (fleet composition and consistent routing) and from operational decisions taken after demand realization. We introduced valid inequalities to strengthen the formulation and characterized key properties of the approximate recourse actions. In addition, we developed a path-based reformulation to enhance tractability, complemented by a Kernel Search-based heuristic.

Numerical experiments were carried out on instances derived from a real-world case study. Small instances were synthetically generated through a population-density estimation approach and used to validate the formulation. We also computed classical stochastic measures, after applying a scenario reduction procedure, to assess the impact of uncertainty. The results not only highlight the advantages of the stochastic methodology over its deterministic counterpart but also provide managerial insights, particularly regarding the use of deterministic solutions as effective starting points in the stochastic setting. Subsequently, computational tests on large-scale instances, generated from a real-world dataset provided by the Italian postal company, were conducted to evaluate the effectiveness of our approach. These results demonstrate the improvements achieved through the path-based reformulation and the additional gains obtained with the Kernel Search-based approach. Overall, the proposed formulations generate consistent routing plans that can be adapted to demand variability and enhanced through the use of light freight vehicles such as ECBs.

Several avenues for future research emerge from this work. First, implementing exact solution methods, such as Benders decomposition, may be promising for improving solution quality. Moreover, incorporating additional sources of uncertainty relevant to Last Mile delivery, such as travel times and energy consumption, represents an important direction for investigation. Finally, the integration of real-time data offers opportunities to move from static toward dynamic routing models, enhancing the flexibility and adaptability of the solutions in real-world applications.

\bibliography{biblio}

\newpage
\section*{Appendix A: Split Approximate Recourse Actions - example}\label{sec:appendix_a}

\begin{figure}[h!]
    \centering
    \begin{tikzpicture}[roundnode/.style={circle, draw=black!60, fill=black!5, very thick, minimum size=0.3mm},squarednode/.style={rectangle, draw=black!60, fill=black!5, very thick, minimum size=8mm}]
    \begin{scope}[scale=0.7]
        \node[squarednode] (0) at (0,0) {0};
        \node[roundnode] (1) at (5,1.5) {\tiny 1};
        \node[roundnode] (2) at (1,5) {\tiny 2};
        \node[roundnode] (3) at (-3,3) {\tiny 3};
        \node[roundnode] (4) at (-1.5,-3) {\tiny 4};
        \node[roundnode] (5) at (3.5,-2) {\tiny 5};
        \node[roundnode] (6) at (2,6.2) {\tiny 6};
        \node[roundnode] (7) at (-2,-4.2) {\tiny 7};
        \node[roundnode] (8) at (-1.5,6.2) {\tiny 8};
        \node[roundnode] (9) at (1,-1.8) {\tiny 9};
        \node[roundnode] (10) at (-3.5,0.5) {\tiny 10};
        \node[roundnode] (11) at (3.5,3.5) {\tiny 11};

        \path (0) edge[->,ultra thick, dotted, bend right,blue]node[right]{} (2);
        \path (2) edge[->,ultra thick, dotted, bend right,blue]node[above]{}  (3);
        \path (3) edge[->,ultra thick, dotted, bend right,blue]node[right]{} (0);
        \path (3) edge[<->,ultra thick,orange] node[above left]{} (10);

        \path (0) edge[->,ultra thick, dash dot, bend left,red]node[below]{}(1);
        \path (1) edge[->,ultra thick, dash dot, bend left,red]node[right]{} (5);
        \path (5) edge[->,ultra thick, dash dot, bend left,red]node[below]{} (4);
        \path (4) edge[->,ultra thick, dash dot, bend left,red] node[above left]{} (0);
        \path (4) edge[<->,ultra thick,orange]node[left]{} (10);
        \path (1) edge[<->,ultra thick,orange]node[above right]{} (11);


    \end{scope}
    \begin{scope}[shift={(6,3.5)}] 
        \draw[](-1.5,0.5) node[right, black]{\footnotesize\textbf{Vehicle types:}};
        \draw[ultra thick, dash dot, red] (-1,0) -- (-0.5,0) node[right, black]{: \footnotesize CMs};
        \draw[ultra thick, dotted, blue] (-1,-0.5) -- (-0.5,-0.5) node[right, black]{: \footnotesize ECBs};
        \draw[](-1.5,-1.25) node[right, black]{\footnotesize \textbf{Load capacity:}};
        \draw[](-1.5,-1.75) node[right, black]{\footnotesize $l_{\text{CMs}}=10$ units};
        \draw[](-1.5,-2.25) node[right, black]{\footnotesize $l_{\text{ECBs}}=5$ units};
        \draw[](-1.5,-3) node[right, black]{\footnotesize \textbf{Demand of customers in $1^{st}$-stage routing:}};
        \draw[](-1.5,-3.5) node[right, black]{\footnotesize $d_{1s}=d_{3s}=d_{4s}=2$ units};
        \draw[](-1.5,-4) node[right, black]{\footnotesize $d_{2s}=d_{5s}=1$ unit};
        \draw[](-1.5,-4.75) node[right, black]{\footnotesize \textbf{Additional customers demand:}};
        \draw[](-1.5,-5.25) node[right, black]{\footnotesize $d_{6s}=d_{7s}=d_{8s}=d_{9s}=0$ units};
        \draw[](-1.5,-5.75) node[right, black]{\footnotesize $d_{10s}=3$ units};
        \draw[](-1.5,-6.25) node[right, black]{\footnotesize $d_{11s}=4$ units};
    \end{scope}
        
    \end{tikzpicture}
    \caption{\small Example of delivery network with split approximate recourse actions (orange double arrows).}
    \label{fig:example_split_delivery}
\end{figure}
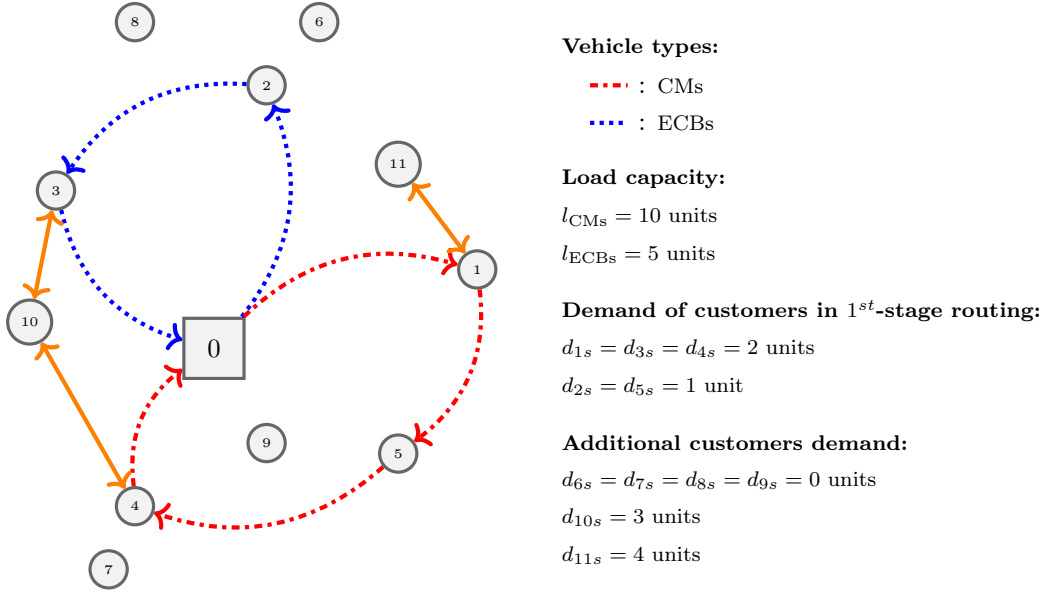

Consider a delivery network composed of 11 customers served by a fleet that includes two vehicle types. Figure \ref{fig:example_split_delivery} displays the tactical routes for both vehicles, along with their load capacities and the customer demand under a given scenario $s \in \mathcal{S}$.

According to Figure \ref{fig:example_split_delivery}, the values of $u_{is}$ for each customer $i$ appearing in the tactical routes are reported in Table \ref{tab:values_uis}.

\begin{table}[h!]
    \centering
    \scalebox{0.8}{
    \begin{tabular}{l|l|l}
    \toprule
    Vehicle & \makecell{$1^{st}$-stage\\distributed\\ demand} & \makecell{$2^{nd}$-stage\\distributed\\ demand} \\
    \midrule
         ECB & $(u_{2s},u_{3s}) = (1,3)$ & \makecell{$(u_{2s},u_{11s}, u_{3s}) = (1,5,\textbf{7})$\\ $(u_{2s},u_{3s}, u_{10s}) = (1,3,\textbf{6})$} \\
         \hline
         CM & $(u_{1s}, u_{5s}, u_{4s})  = (2, 3 ,5)$ & \makecell{$(u_{1s},u_{11s}, u_{5s},u_{4s},u_{10s}  ) = (2,6,7,9,\textbf{12})$\\ $(u_{1s},u_{5s}, u_{4s},u_{10s} ) = (2,3,5,8)$}\\
    \bottomrule
    \end{tabular}}
    \caption{\small Values of the accumulated satisfied demand $u_{is}$ including additional customers under scenario $s$ without split recourse actions.}
    \label{tab:values_uis}
\end{table}

In this specific scenario, neither vehicle is able to fully satisfy the additional customer demand. As shown in the second column of Table \ref{tab:values_uis}, the condition $u_{hs} \leq l_p$ is violated for some customers, for all vehicles $p \in \mathcal{P}$ and some customers $h \in \mathcal{I}$. In particular, the ECB cannot serve both additional customers because its load capacity is insufficient, as highlighted by the bold values in Table \ref{tab:values_uis}. Similarly, the CM can serve at most one of the two additional customers in full; serving both would exceed its load capacity, and the demand of at least one customer would need to be split. Only when one of the additional customers is excluded does the CM avoid split recourse actions.
\begin{table}[h!]
    \centering
    \scalebox{0.8}{
    \begin{tabular}{l|l|l}
    \toprule
    Vehicle & \makecell{$1^{st}$-stage\\distributed\\ demand} & \makecell{$2^{nd}$-stage\\distributed\\demand} \\
    \midrule
         ECB & $(u_{2s},u_{3s}) = (1,3)$ & \makecell{$(u_{2s},u_{11s}, u_{3s}) = (1,3,5)$\\ $(u_{2s},u_{3s}, u_{10s}) = (1,3,5)$\\ $(u_{2s},u_{11s},u_{3s}, u_{10s}) = (1,2,4, 5)$}\\
         \hline
         CM & $(u_{1s}, u_{5s}, u_{4s})  = (2, 4 ,5)$ & \makecell{$(u_{1s},u_{11s}, u_{5s},u_{4s},u_{10s}  ) = (2,4,5,7,10)$\\ $(u_{1s},u_{11s}, u_{5s},u_{4s}, u_{10s}) = (2,6,7,9,10)$\\$(u_{1s},u_{11s}, u_{5s},u_{4s}, u_{10s}) = (2,5,6,8,10)$}\\
    \bottomrule
    \end{tabular}}
    \caption{\small Values of the accumulated satisfied demand $u_{is}$ including the additional customers under scenario $s$ with split approximate recourse actions.}
    \label{tab:correct_values_uis}
\end{table}

Thus, in the absence of split approximate recourse actions, the only way to fully satisfy all delivery requests would be to rely on an external delivery service, incurring higher operational costs. A feasible adjustment to the first-stage routing is illustrated by the orange double arrows in Figure \ref{fig:example_split_delivery}. In this case, all customer demands can be satisfied, as summarized in Table \ref{tab:correct_values_uis}, which provides all possible distributions of the delivered demand.

As shown in Table \ref{tab:correct_values_uis}, split approximate recourse actions are not only necessary under certain demand realizations but also offer multiple feasible configurations that enable full demand satisfaction without resorting to external delivery services, thereby reducing overall operating costs.

\newpage
\section*{Appendix B: Small Instances Extended Results}\label{sec:appendix_b}

\begin{table}[h!]
  \centering
  \resizebox{\textwidth}{!}{
    \begin{tabular}
    {l|c|c|c|c|c|c|c|c|c|c|c|c|c}
    \toprule
    Instance & $\mid \mathcal{S} \mid $ & CMs & ECBs &  \makecell{CMs distance\\ (km)} & \makecell{ECBs distance\\(km)}
 & \makecell{Complete\\ recourse \\actions} & \makecell{Split \\recourse \\actions} &  \makecell{Complete\\ unserved \\customers}& \makecell{Partial\\ unserved \\customers}& \makecell{Average\\recourse\\distance (km)} &  \makecell{$RP$\\ (\euro/d)} & \makecell{CPU\\(s)} & Gap\\
    \midrule
    20\_2\_5\_1     &  & 3     & 1     & 35,0688 & 1,5415 & 0 & 2 & 0 & 1 & 1,07 &  34,2889 & 3251,2 & 0\% \\
    20\_2\_5\_2     & & 3     & 1     & 29,0034 & 10.0840 & 0 & 2 & 0 & 0 & 0,45 &  31,3914 & 2764,1 & 0\% \\
    20\_2\_5\_3     & & 3     & 1     & 33,4385 & 1,5415 & 1 & 4 & 0 & 0 & 2,57 & 31,0586 & 503,1 & 0\% \\
    20\_2\_5\_4     & & 3     & 1     & 26,6817 & 10,1679 & 1 & 5 & 0 & 1 &  2,65 & 34,0130 & 797,6 & 0\% \\
    20\_2\_5\_5    &  & 3     & 1     & 33,3191 & 1,5415 & 0 & 5 & 0 & 0 & 2,37 & 31,3103 & 1133,8 & 0\% \\
    20\_2\_5\_6     &5   & 3     & 1     & 32,4141 & 1,5415 & 0 & 3 & 0 & 0 & 1,75 & 31,0206 & 203,1 & 0\% \\
    20\_2\_5\_7     & & 3     & 0     & 33,6741 & 0,0000 & 0 & 3 & 0 & 0 & 1,63 & 28,0204 & 312,9 & 0\% \\
    20\_2\_5\_8     & & 3     & 1     & 34,6233 & 1,5415 & 0 & 5 & 0 & 0 & 1,49 & 31,4162 & 517,8 & 0\% \\
    20\_2\_5\_9     & & 3     & 0     & 34,8994 & 0,0000 & 0 & 3 & 0 & 0 & 1,38 & 28,2220 & 434,3 & 0\% \\
    20\_2\_5\_10     & & 2     & 2     & 24,5220 & 11,5675 & 0 & 3 & 0 & 0 & 1,81 & 26,9562 & 818,4 & 0\% \\
    \midrule

    20\_2\_10\_1    &   & 3     & 1     & 35,4982 & 1,5415 & 1 & 9 & 0 & 1 & 1,07 & 32,7672 & 3411,6 & 0\% \\
    20\_2\_10\_2    &  & 3     & 2     & 29,1347 & 10,9394 & 0 & 4 & 0 & 0 & 0,64 & 34,5793 & 2252,4 & 0\% \\
    20\_2\_10\_3    &  & 3     & 2     & 34,1143 & 3,9569 & 1 & 6 & 0 & 0 & 1,45 & 34,6703 & 884,8 & 0\% \\
    20\_2\_10\_4    &  & 3     & 0     & 34,8791 & 0,0000 & 1 & 6 & 0 & 1 & 0,97 & 30,6459 & 751,9 & 0\% \\
    20\_2\_10\_5    &  & 3     & 0     & 35,8044 & 0,0000 & 3 & 3 & 0 & 0 & 1,70 & 28,4586 & 2254,6 & 0\% \\
    20\_2\_10\_6    & 10  & 3     & 1     & 35,0954 & 2,4155 & 1 & 5 & 0 & 1 & 1,60 & 34,1613 & 2037,3 & 0\% \\
    20\_2\_10\_7    &  & 3     & 1     & 34,5066 & 2,4155 & 0 & 10 & 0 & 0 & 5,82 & 31,3761 & 4754,2 & 0\% \\
    20\_2\_10\_8    &  & 3     & 1     & 26,5928 & 10,1679 & 3 & 6 & 0 & 0 & 1,32 & 31,3058 & 3497,9 & 0\% \\
    20\_2\_10\_9    &   & 3     & 0     & 34,8792 & 0,0000 & 2 & 3 & 0 & 1 & 0,71 & 29,5282 & 739,8 & 0\% \\
    20\_2\_10\_10    &  & 3     & 2     & 28,2788 & 11,6320 & 2 & 8 & 0 & 0 & 1,47 & 34,6581 & 1734,5 & 0\% \\

    \midrule
    20\_2\_20\_1    &  & 3     & 1     & 34,8634 & 1,5415 & 2 & 17 & 0 & 1 & 1,81 & 33,3599 & 7501,1 & 0\% \\
    20\_2\_20\_2    &  & 3    & 1     & 32,4141 & 1,5415 & 3 & 22 & 0 & 0 & 3,71 & 31,3628 & 4090,1 & 0\% \\
    20\_2\_20\_3    &  & 3     & 1     & 33,4385 & 1,5415 & 4 & 15 & 0 & 1 & 2,76 & 34,2589 & 12898,8 & 0\% \\
    20\_2\_20\_4    & & 3     & 2 & 26,5127 & 11,7094 & 5 & 23 & 0 & 1 & 3,26 & 35,3446 & 5992,3 & 0\% \\
    20\_2\_20\_5    & & 3     & 1     & 34,7440 & 1,5415 & 2 & 21 & 0 & 0 & 1,71 & 31,4793 & 7231,8 & 0\% \\
    20\_2\_20\_6    & 20 & 3     & 1     & 33,4557 & 1,7243 & 5 & 13 & 0 & 0 & 2,75 & 31,4304 & 3764,0 & 0\% \\
    20\_2\_20\_7    & & 3     & 1     & 33,3191 & 1,5415 & 1 & 18 & 0 & 1 & 2,15 & 34,1279 & 2474,1 & 0\% \\
    20\_2\_20\_8    &  & 3     & 1   & 35,3887 & 1,5415 & 2 & 18 & 0 & 1 & 1,25 & 34,0273 & 7750,3 & 0\% \\
    20\_2\_20\_9    &  & 3     & 1     & 34,8634 & 1,5415 & 2 & 15 & 0 & 0 & 1,18 & 31,4109 & 10894,8 & 0\% \\
    20\_2\_20\_10    & & 3     & 1     & 34,8634 & 1,5415 &  2 & 20 & 0 & 1 & 1,26 & 33,9252 & 6838,2 & 0\% \\

    \midrule
    20\_2\_50\_1    &   & 3    &   1    &   34,7440   &  1,5415  & 6 & 41 & 0 & 2 & 1,63 & 32,0011  &    32397,4   & 0\% \\
    20\_2\_50\_2    &   & 3    &   1    &   34,8634   &  1,5415  &  6 & 42 & 0 & 4 & 1,89 & 34,4871  &   69239,9   & 0\% \\
    20\_2\_50\_3    &  & 3    &   2    &   32,7061   &  5,1023  &  4 & 46 & 0 & 1 & 1,87 & 34,8846  &    14878,6   & 0\% \\
    20\_2\_50\_4    &  & 3    &   1    &   35,0744   &  1,5415 & 10 & 36 & 0 & 2 & 1,52 & 33,6551  & 63521,9 & 0\% \\
    20\_2\_50\_5    &   & 3    &   1    &   34,7440   &  1,5415  & 5 & 31 & 1 & 1 & 1,20 & 34,0569  & 14349,1 & 0\% \\
    20\_2\_50\_6    &   50 & 4    &   0    &   36,5047   &  0,0000  &  5 & 16 & 0 & 0 & 1,18 & 35,5073  & 19313,2 & 0\% \\
    20\_2\_50\_7    &  & 3    &   1    &   34,8634   &  1,5415  & 3 & 33 & 0 & 3 & 0,99 & 32,4277  & 19114,5 & 0\% \\
    20\_2\_50\_8    &   & 3   &   1    &   34,8634   &  1,5415 &  5 & 35 & 0 & 2 & 1,43 &  33,7030  & 6542,8 & 0\% \\
    20\_2\_50\_9    &   & 3    &   1    &   26,8845   &  10,1679  & 13 & 46 & 0 & 1 & 2,79 & 32,1396  & 49785,3 & 0\% \\
    20\_2\_50\_10   &    & 3    &   1    &   34,8634   &  1,5415  & 6 & 33 & 0 & 1 & 1,24 & 31,7059  & 18580,2 & 0\% \\

    \midrule
    20\_2\_100\_1   &  & 3 & 1 &  35,0743 & 1,5415 & 28 & 78 & 0 & 4 & 1,69 & 34,0300 & 78813,2  & 0\% \\
    20\_2\_100\_2   & & 3 & 1 &  34,8643 & 1,5415 & 11 & 83 & 0 & 1 & 1,46 & 33,2102 & 50054,2  & 0\% \\
    20\_2\_100\_3   & & 3 & 1 &  34,7616 & 1,5415 & 46 & 90 & 0 & 5 & 2,08 & 34,4217 & 56786,2  & 0\% \\
    20\_2\_100\_4   &  & 3 & 1 &  34,7440 & 1,5415 & 12 & 82 & 0 & 4 & 1,50 & 33,5315 & 40166,8  & 0\% \\
    20\_2\_100\_5   & & 3 & 1 &  34,7616 & 1,5415 & 24 & 75 & 0 & 4 & 1,67 & 33,4149 & 41016,3  & 0\% \\
    20\_2\_100\_6   & 100  & 3 & 1 & 35,0763 & 1,5415 & 14 & 71 & 0 & 4 & 1,38 & 32,6634 & 43201,7  & 0\% \\
    20\_2\_100\_7   &   & 3 & 1 & 33,0324 & 1,5415 & 12 & 93 & 0 & 5 & 2,25 & 33,7571 & 68983,1  & 0\% \\
    20\_2\_100\_8   &  & 3 & 1 &  35,4919 & 1,5415 & 22 & 76 & 0 & 3 & 1,28 & 33,5757 & 73324,5  & 0\% \\
    20\_2\_100\_9   &   & 3 & 1 &  33,3191 & 1,5415 & 38 & 105 & 0 & 2 & 3,84 & 32,2269 & 50071,9  & 0\% \\
    20\_2\_100\_10   & & 3 & 1 &  34,8822 & 1,5415 & 51 & 92 & 1 & 4 & 2,26 & 33,5384 & 66931,8  & 0\% \\

    \midrule
    20\_2\_150\_1  &  & 3 & 1 &  36,0683 & 1,5415 & 30  & 93 & 0 & 6 & 1,47 & 33,5974 & $86400,0^\dagger$  & $1,24\%$ \\
    20\_2\_150\_2   & & 3 & 1 &  36,3974 & 1,5415 & 33 & 101 & 0 & 7 & 1,58 & 33,4026 & 73206,1  & 0\% \\
    20\_2\_150\_3   & & 3 & 1 &  34,0628 & 3,5608 & 14 & 74 & 0 & 4 & 1,26 & 33,5069 & $86400,0^\dagger$  & $1,35\%$ \\
    20\_2\_150\_4   & &  3 & 1 &  35,2592 & 1,5415 & 32 & 96 & 0 & 6 & 1,74 & 33,0513 & 65402,1 & 0\% \\
    20\_2\_150\_5   &  & 3 & 1 &  34,8229 & 3,5608 & 27 & 109 & 0 & 3 & 1,15 & 33,3025 & 72391,4 & 0\% \\
    20\_2\_150\_6  & 150 & 3 & 1 &  30,4064 & 10,0901 & 42 & 105 & 0 & 12 & 3,65 & 34,6471 & $86400,0^\dagger$   & $1,03\%$ \\
    20\_2\_150\_7   & & 3 & 1 &  35,4622 & 1,5415 & 34 & 101 & 0 & 4 & 1,43 & 33,3391 & 68946,2  & 0\% \\
    20\_2\_150\_8   & & 3 & 1 &  34,4123 & 1,5415 & 31 & 94 & 0 & 6 & 3,42 & 33,3349 & $86400,0^\dagger$  & $0,94\%$ \\
    20\_2\_150\_9   & & 3 & 1 &  35,4981 & 1,5415 & 29 & 89 & 0 & 3 & 1,71 & 32,8219 & 68920,6  & 0\% \\
    20\_2\_150\_10   & & 3 & 1 &  33,9565 & 1,5415 & 29 & 97 & 0 & 7 & 3,36 & 33,3681 & 70239,6  & 0\% \\
    \bottomrule
    \end{tabular}}
  \caption{\small Computational results of the model for stochastic instances with increasing number of scenarios. The superscript $^\dagger$
 on the CPU time indicates that the time limit was reached.}
\label{stochastic_computational_results}
\end{table}

\newpage

\section*{Appendix C: Stochastic Measures - Small Instances Extended Results}\label{sec:appendix_c}

\begin{table}[h!]
    \centering
    \resizebox{\textwidth}{!}{
    \begin{tabular}{l|c|c|c|c|c|c|c|c|c|c|c|c}
    \toprule
    Instance & CMs & ECBs & \makecell{CMs distance\\ (km)} & \makecell{ECBs distance\\(km)}
 & \makecell{Complete\\ recourse \\actions} & \makecell{Split \\recourse \\actions} &  \makecell{Complete\\ unserved \\customers}& \makecell{Partial\\ unserved \\customers}& \makecell{Average recourse\\distance (km)} &  \makecell{$EV$\\(\euro/d)} & \makecell{CPU\\(s)} & Gap\\
    \midrule
        20\_2\_1\_EV & 2 & 2 & 23,6581 & 13,6438 & 0 & 0 & 0 & 0 & 0,00 & 27,1944 & 2594,1 & 0,00\% \\
        20\_2\_2\_EV & 2 & 2 & 23,7254 & 14,9705 & 0 & 0 & 0 & 0 & 0,00 & 27,0981 & 7641,1 & 0,00\% \\
        20\_2\_3\_EV & 2 & 2 & 23,6466 & 13,6438 & 0 & 0 & 0 & 0 & 0,00 & 27,1070 & 6473,1 & 0,00\%\\
        20\_2\_4\_EV & 3 & 0 & 34,8727 & 0,0000 & 0 & 0 & 0 & 0 & 0,00 & 27,9745 & 247,5 & 0,00\%\\
        20\_2\_5\_EV & 3 & 0 & 34,8792 & 0,0000 & 0 & 0 & 0 & 0 & 0,00 & 27,5756 & 875,3 & 0,00\%\\
        20\_2\_6\_EV & 3 & 0 & 34,8727 & 0,0000 & 0 & 0 & 0 & 0 & 0,00 & 27,9745 & 333,7 & 0,00\%\\
        20\_2\_7\_EV & 3 & 0 & 34,8423 & 0,0000 & 0 & 0 & 0 & 0 & 0,00 & 27,6445 & 3782,4 & 0,00\%\\
        20\_2\_8\_EV  & 2 & 2 & 23,6581 & 13,6438 & 0 & 0 & 0 & 0 & 0,00 & 27,2800 & 3591,4 & 0,00\%\\
        20\_2\_9\_EV  & 3 & 0 & 34,8592 & 0,0000 & 0 & 0 & 0 & 0 & 0,00 & 27,8658 & 1342,7 & 0,00\%\\
        20\_2\_10\_EV & 2 & 2 & 23,2997 & 14,6699 & 0 & 0 & 0 & 0 & 0,00 & 27,2421 & 3572,4 & 0,00\%\\
        \midrule
    Average & & & 29,2364 & 7,0572 & & & & & 0,00 & 27,4956 & 3045,6 & 0,00\%\\
    \bottomrule
    \end{tabular}}
    \caption{\small Computational results provided by the $EV$ on each stochastic small instance.}
    \label{tab:EV_computational_results}
\end{table}

\begin{table}[h!]
    \centering
    \resizebox{\textwidth}{!}{
    \begin{tabular}{l|c|c|c|c|c|c|c|c|c|c|c|c}
    \toprule
    Instance & CMs & ECBs & \makecell{CMs distance\\ (km)} & \makecell{ECBs distance\\(km)} 
 & \makecell{Complete\\ recourse \\actions} & \makecell{Split \\recourse \\actions} &  \makecell{Complete\\ unserved \\customers}& \makecell{Partial\\ unserved \\customers}& \makecell{Average recourse\\distance (km)} &  \makecell{$EEV$\\(\euro/d)} & \makecell{CPU\\(s)} & Gap\\
    \midrule
        $20\_2\_1\_EEV_{FR}$ & 2 & 2 & 23,6581 & 13,6438 & 52 & 143 & 14 & 40 & 4,06 & 60,3214 & 0,67 & 0,00\% \\
        $20\_2\_2\_EEV_{FR}$ & 2 & 2 & 23,7254 & 14,9705 & 46 & 159 & 10 & 43 & 3,47 & 53,6357 & 0,43 & 0,00\% \\
        $20\_2\_3\_EEV_{FR}$  & 2 & 2 & 23,6466 & 13,6438 & 59 & 157 & 14 & 49 & 4,41 & 61,7699 & 0,55 & 0,00\%\\
        $20\_2\_4\_EEV_{FR}$  & 3 & 0 & 34,8727 & 0,0000 & 17 & 70 & 0 & 20 & 2,32 & 36,2036 & 0,49 & 0,00\%\\
        $20\_2\_5\_EEV_{FR}$ & 3 & 0 & 34,8792 & 0,0000 & 16 & 68 & 1 & 20 & 2,19 & 38,7154 & 0,47 & 0,00\%\\
        $20\_2\_6\_EEV_{FR}$  & 3 & 0 & 34,8727 & 0,0000 & 14 & 76 & 0 & 24 & 2,45 & 37,0564 & 0,52 & 0,00\%\\ 
        $20\_2\_7\_EEV_{FR}$  & 3 & 0 & 34,8423 & 0,0000 & 18 & 69 & 1 & 16 & 2,66 & 35,8342 & 0,48 & 0,00\%\\
        $20\_2\_8\_EEV_{FR}$ & 2 & 2 & 23,6581 & 13,6438 & 56 & 161 & 13 & 50 & 4,60 & 61,4908 & 0,56 & 0,00\%\\
        $20\_2\_9\_EEV_{FR}$ & 3 & 0 & 34,8592 & 0,0000 & 31 & 64 & 3 & 18 & 2,33 & 37,0132 & 0,63 & 0,00\%\\
        $20\_2\_10\_EEV_{FR}$  & 2 & 2 & 23,2997 & 14,6699 & 52 & 131 & 10 & 46 & 3,69 & 61,6554 & 0,59 & 0,00\%\\
    \midrule
    Average & & & 29,2364 & 7,0572 & & & & & 3,22 & 48,3696 & 0,54 & 0,00\%\\
    \midrule
        $20\_2\_1\_EEV_{F}$ & 2 & 2 & 22,2144 & 13,6438 & 35 & 167 & 14 & 41 & 4,06 & 60,2146 & 8212,2 & 0,00\% \\
        $20\_2\_2\_EEV_{F}$  & 2 & 2 & 21,9017 & 14,9705 & 40 & 186 & 10 & 42 & 5,37 & 53,6245 & 25681,2 & 0,00\% \\
        $20\_2\_3\_EEV_{F}$   & 2 & 2 & 23,2997 & 14,6699 & 75 & 173 & 15 & 48 & 5,45 & 61,6935 & 2347,5 & 0,00\%\\
        $20\_2\_4\_EEV_{F}$ & 3 & 0 & 35,2289 & 0,0000 & 20 & 78 & 0 & 20 & 1,70 & 36,1716 & 6087,1 & 0,00\%\\
        $20\_2\_5\_EEV_{F}$  & 3 & 0 & 35,3405 & 0,0000 & 28 & 74 & 0 & 21 & 1,49 & 38,6831 & 3901,7 & 0,00\%\\
        $20\_2\_6\_EEV_{F}$  & 3 & 0 & 34,8727 & 0,0000 & 14 & 76 & 0 & 24 & 2,45 & 37,0564 & 7138,4 & 0,00\%\\
        $20\_2\_7\_EEV_{F}$ & 3 & 0 & 35,0076 & 0,0000 & 15 & 70 & 2 & 14 & 2,52 & 35,8247 & 548,9 & 0,00\%\\
        $20\_2\_8\_EEV_{F}$   & 2 & 2 & 22,0976 & 13,6438 & 92 & 172 & 14 & 47 & 6,01 & 61,4353 & 7932,6 & 0,00\%\\
        $20\_2\_9\_EEV_{F}$  & 3 & 0 & 34,8792 & 0,0000 & 31 & 64 & 3 & 18 & 2,33 & 37,0132 & 1837,1 & 0,00\%\\
        $20\_2\_10\_EEV_{F}$  & 2 & 2 & 21,8676 & 14,6699 & 32 & 161 & 9 & 47 & 4,44 & 61,5263 & 5498,2 & 0,00\%\\
        \midrule
    Average  & &  & 28,6710 & 7,1598 & & & &  &3,58 & 48,3243 & 6918,5 & 0,00\%\\
    \bottomrule
    \end{tabular}}
    \caption{\small Computational results provided by the $EEV_{FR}$ and $EEV_{F}$ on each stochastic instance with $\mid \mathcal{S}\mid = 100$.}
    \label{tab:EEV_computational_results}
\end{table}

\begin{table}[h!]
    \centering
    \resizebox{\textwidth}{!}{
    \begin{tabular}{l|c|c|c|c|c|c|c|c|c|c|c|c|c}
    \toprule
    Instance & CMs & ECBs & \makecell{CMs distance\\ (km)} & \makecell{ECBs distance\\(km)}  
 & \makecell{Complete\\ recourse \\actions} & \makecell{Split \\recourse \\actions} &  \makecell{Complete\\ unserved \\customers}& \makecell{Partial\\ unserved \\customers}& \makecell{Average recourse\\distance (km)} &  \makecell{$EIV$\\(\euro/d)} & \makecell{CPU\\(s)} & Gap\\
    \midrule
        $20\_2\_1\_EIV_{FR}$  & 2 & 2 & 23,6581 & 13,6438 & 52 & 143 & 14 & 40 & 4,06 & 60,3214 & 0,77 & 0,00\% \\
        $20\_2\_2\_EIV_{FR}$  & 2 & 2 & 23,7254 & 14,9705 & 46 & 159 & 10 & 43 & 3,47 & 53,6357 & 0,58 & 0,00\% \\
        $20\_2\_3\_EIV_{FR}$  & 2 & 2 & 23,6466 & 13,6438 & 59 & 157 & 14 & 49 & 4,41 & 61,7699 & 0,51 & 0,00\%\\
        $20\_2\_4\_EIV_{FR}$ & 3 & 0 & 34,8727 & 0,0000 & 17 & 70 & 0 & 20 & 2,32 & 36,2036 & 0,45 & 0,00\%\\
        $20\_2\_5\_EIV_{FR}$ & 3 & 0 & 34,8792 & 0,0000 & 16 & 68 & 1 & 20 & 2,19 & 38,7154 & 0,87 & 0,00\%\\
        $20\_2\_6\_EIV_{FR}$  & 3 & 0 & 34,8727 & 0,0000 & 14 & 76 & 0 & 24 & 2,45 & 37,0564 & 0,72 & 0,00\%\\ 
	$20\_2\_7\_EIV_{FR}$  & 3 & 0 & 34,8423 & 0,0000 & 18 & 69 & 1 & 16 & 2,66 & 35,8342 & 0,81 & 0,00\%\\
        $20\_2\_8\_EIV_{FR}$ & 2 & 2 & 23,6581 & 13,6438 & 56 & 161 & 13 & 50 & 4,60 & 61,4908 & 0,39 & 0,00\%\\
        $20\_2\_9\_EIV_{FR}$  & 3 & 0 & 34,8592 & 0,0000 & 31 & 64 & 3 & 18 & 2,33 & 37,0132 & 0,73 & 0,00\%\\
        $20\_2\_10\_EIV_{FR}$  & 2 & 2 & 23,2997 & 14,6699 & 52 & 131 & 10 & 46 & 3,69 & 61,6554 & 0,52 & 0,00\%\\
        \midrule
    Average & & &  29,2364 & 7,0572 & & & & & 3,22 & 48,3696 & 0,64 & 0,00\%\\
    \midrule
        $20\_2\_1\_EIV_{F}$  & 3 & 2 & 24,2558 & 12,7200 & 41 & 136 & 0 & 0 & 5,48 & 34,7305 & 11960,6 & 0,00\% \\
        $20\_2\_2\_EIV_{F}$  & 2 & 4 & 22,4414 & 19,3319 & 13 & 106 & 0 & 2 & 2,07 & 34,6758 & 18873,7 & 0,00\% \\
        $20\_2\_3\_EIV_{F}$ & 2 & 4 & 22,1566 & 19,5777 & 21 & 118 & 0 & 4 & 2,52 & 34,4339 & 4374,1 & 0,00\%\\
        $20\_2\_4\_EIV_{F}$  & 3 & 1 & 34,3780 & 3,5608 & 12 & 79 & 0 & 1 & 1,31 & 32,9113 & 9408,6 & 0,00\%\\
        $20\_2\_5\_EIV_{F}$  & 3 & 1 & 30,4693 & 10,0905 & 8 & 65 & 0 & 4 & 1,63 & 33,0521 & 4382,6 & 0,00\%\\
        $20\_2\_6\_EIV_{F}$  & 3 & 1 & 33,5473 & 3,5608 & 37 & 79 & 0 & 1 & 2,29 & 32,9304 & 8475,9 & 0,00\%\\
        $20\_2\_7\_EIV_{F}$  & 3 & 1 & 30,5905 & 10,0905 & 13 & 76 & 0 & 3 & 1,07 & 33,1873 & 927,8 & 0,00\%\\
        $20\_2\_8\_EIV_{F}$  & 2 & 4 & 22,3107 & 18,0052 & 24 & 143 & 0 & 4 & 4,28 & 34,5814 & 4187,4 & 0,00\%\\
        $20\_2\_9\_EIV_{F}$ & 3 & 1 & 33,1886 & 3,1139 & 39 & 90 & 0 & 5 & 3,56 & 33,3484 & 6892,3 & 0,00\%\\
        $20\_2\_8\_EIV_{F}$ & 2 & 4 & 21,8488 & 17,5583 & 25 & 146 & 0 & 1 & 4,19 & 34,0193 & 5029,6 & 0,00\%\\
        \midrule
    Average &  & & 27,5187 & 11,7610 & & & & & 2,84 & 33,7843 & 7448,6 & 0,00\%\\
    \bottomrule
    \end{tabular}}
    \caption{\small Computational results provided by the $EIV_{FR}$ and $EIV_{F}$ on each stochastic small instance with $\mid \mathcal{S}\mid = 100$.}
    \label{tab:EIV_computational_results}
\end{table}

\newpage

\end{document}